\PassOptionsToPackage{dvipsnames}{xcolor}
\documentclass{lmcs}
\pdfoutput=1

\usepackage{lastpage}
\lmcsdoi{22}{2}{13}
\lmcsheading{}{\pageref{LastPage}}{}{}%
{Jun.~10,~2025}{May~05,~2026}{}

\usepackage[T2A,T1]{fontenc}			% for cyrillic characters in Golubtsov ref

\usepackage{amsmath,amsfonts,amssymb}
\usepackage{bm}
\usepackage{enumitem}
\setlist[enumerate]{label=(\roman*)}  
\setlist[enumerate,2]{label=(\alph*)} 
\usepackage{mathtools}
\usepackage{tikz}
	\usetikzlibrary{external,arrows}
	\usetikzlibrary{decorations.pathreplacing}	% for curly brace

\usepackage{tikz-cd}
\usepackage{tikzit}				% String diagrams
\tikzstyle{morphism}=[fill=white, draw=black, shape=rectangle]
\tikzstyle{medium box}=[fill=white, draw=black, shape=rectangle, minimum width=0.7cm, minimum height=0.7cm]
\tikzstyle{large morphism}=[fill=white, draw=black, shape=rectangle, minimum width=1.7cm, minimum height=1cm]
\tikzstyle{bn}=[fill=black, draw=black, shape=circle, inner sep=1.5pt]
\tikzstyle{state}=[fill=white, draw=black, regular polygon, regular polygon sides=3, minimum width=0.8cm, shape border rotate=180, inner sep=0pt]
\tikzstyle{medium state}=[fill=white, draw=black, regular polygon, regular polygon sides=3, minimum width=1.3cm, inner sep=0pt, shape border rotate=180]
\tikzstyle{large state}=[fill=white, draw=black, regular polygon, regular polygon sides=3, minimum width=2.2cm, shape border rotate=180, inner sep=0pt]
\tikzstyle{wide state}=[fill=white, draw=black, shape=isosceles triangle, minimum width=0.8cm, shape border rotate=270, inner sep=1.4pt, minimum height=0.5cm, isosceles triangle apex angle=80]
\tikzstyle{wn}=[fill=white, draw=black, shape=circle, inner sep=1.5pt]
\tikzstyle{blue morphism}=[fill=white, draw={rgb,255: red,15; green,0; blue,150}, shape=rectangle, text={rgb,255: red,15; green,0; blue,150}, tikzit category=blue]
\tikzstyle{red morphism}=[fill=white, draw={rgb,255: red,150; green,0; blue,2}, shape=rectangle, text={rgb,255: red,150; green,0; blue,2}, tikzit category=red]
\tikzstyle{blue state}=[fill=white, draw={rgb,255: red,15; green,0; blue,150}, shape=circle, regular polygon, regular polygon sides=3, minimum width=0.8cm, shape border rotate=180, inner sep=0pt, text={rgb,255: red,15; green,0; blue,150}, tikzit category=blue]
\tikzstyle{blue node}=[fill={rgb,255: red,15; green,0; blue,150}, draw={rgb,255: red,15; green,0; blue,150}, shape=circle, tikzit category=blue, inner sep=1.5pt]
\tikzstyle{blue}=[text={rgb,255: red,15; green,0; blue,150}, tikzit draw={rgb,255: red,191; green,191; blue,191}, tikzit category=blue, tikzit fill=white, inner sep=0mm]
\tikzstyle{blue wide state}=[fill=white, draw={rgb,255: red,15; green,0; blue,150}, text={rgb,255: red,15; green,0; blue,150}, shape=isosceles triangle, minimum width=0.8cm, shape border rotate=270, inner sep=1.4pt, minimum height=0.5cm, isosceles triangle apex angle=80]
\tikzstyle{red node}=[fill={rgb,255: red,150; green,0; blue,2}, draw={rgb,255: red,150; green,0; blue,2}, shape=circle, inner sep=1.5pt]
\tikzstyle{Purple node}=[fill={rgb,255: red,120; green,0; blue,120}, draw={rgb,255: red,120; green,0; blue,120}, text={rgb,255: red,120; green,0; blue,120}, shape=circle, inner sep=1.5pt]
\tikzstyle{red}=[text={rgb,255: red,150; green,0; blue,2}, inner sep=0mm, tikzit fill=white, tikzit draw={rgb,255: red,191; green,191; blue,191}]
\tikzstyle{purple}=[text={rgb,255: red,150; green,0; blue,150}, inner sep=0mm, tikzit fill=white, tikzit draw={rgb,255: red,191; green,191; blue,191}]
\tikzstyle{white morphism}=[fill=white, draw=white, shape=rectangle, tikzit draw={rgb,255: red,139; green,139; blue,139}]
\tikzstyle{leak morphism}=[fill=white, draw={rgb,255: red,120; green,0; blue,85}, shape=rectangle, text={rgb,255: red,120; green,0; blue,85}, tikzit category=leak]
\tikzstyle{leak}=[text={rgb,255: red,120; green,0; blue,85}, inner sep=0mm, tikzit fill=white, tikzit draw={rgb,255: red,191; green,191; blue,191}, tikzit category=leak]
\tikzstyle{leak node}=[fill={rgb,255: red,120; green,0; blue,85}, draw={rgb,255: red,120; green,0; blue,85}, shape=circle, inner sep=1.5pt, tikzit category=leak]
\tikzstyle{horiz state}=[fill=white, draw=black, regular polygon, regular polygon sides=3, minimum width=1cm, shape border rotate=90, inner sep=0pt]

\tikzstyle{arrow}=[->]
\tikzstyle{dashed box}=[-, dashed]
\tikzstyle{blue arrow}=[-, draw={rgb,255: red,15; green,0; blue,150}, tikzit category=blue]
\tikzstyle{red arrow}=[-, draw={rgb,255: red,150; green,0; blue,2}, tikzit category=red]
\tikzstyle{purple arrow}=[->, draw={rgb,255: red,120; green,0; blue,120}, >=stealth, shorten <=2pt, shorten >=2pt]
\tikzstyle{protected purple arrow}=[->, draw={rgb,255: red,120; green,0; blue,120}, >=stealth, shorten <=2pt, shorten >=2pt, preaction={line width=1.8pt, white, draw}]
\tikzstyle{mapsto}=[{|->}]
\tikzstyle{double wire}=[-, double]
\tikzstyle{protected double wire}=[-, double, preaction={line width=2.5pt,white,draw}]
\tikzstyle{triple wire}=[-, draw, line width=0.4pt, preaction={-, draw, line width=1.4pt, white, preaction={-, draw, line width=2.2pt}}]
\tikzstyle{protected}=[-, preaction={line width=1.8pt,white,draw}]
\tikzstyle{leak arrow}=[-, tikzit draw={rgb,255: red,150; green,0; blue,120}]
\tikzstyle{protected leak arrow}=[-, tikzit draw={rgb,255: red,150; green,0; blue,120}]
\tikzstyle{hollow arrow}=[-, very thin, white, preaction={line width=0.7pt,draw={rgb,255: red,120; green,0; blue,85}}, tikzit category=leak, tikzit draw={rgb,255: red,150; green,0; blue,120}]
\tikzstyle{protected hollow arrow}=[-, very thin, white, preaction={line width=0.7pt,draw={rgb,255: red,120; green,0; blue,85},preaction={line width=2.1pt,white,draw}}, tikzit category=leak, tikzit draw={rgb,255: red,150; green,0; blue,120}]
\tikzstyle{over arrow}=[-, black, preaction={draw=white, double}]
\tikzstyle{purple line}=[-, draw={rgb,255: red,120; green,0; blue,120}]
\tikzstyle{orange line}=[-, draw={rgb,255: red,255; green,100; blue,0}]
\tikzstyle{red line}=[-, draw={rgb,255: red,150; green,0; blue,2}]
\tikzstyle{blue line}=[-, draw={rgb,255: red,15; green,0; blue,150}]
\tikzstyle{blue double arrow}=[-, double, draw={rgb,255: red,15; green,0; blue,150}, tikzit category=blue]
\tikzstyle{red double arrow}=[-, double, draw={rgb,255: red,150; green,0; blue,2}, tikzit category=red]
\tikzstyle{purple double line}=[-, double, draw={rgb,255: red,120; green,0; blue,120}]
\tikzstyle{orange double line}=[-, double, draw={rgb,255: red,255; green,100; blue,0}]
\tikzstyle{curly brace}=[decorate, decoration={brace,amplitude=5pt}]

\tikzstyle{d-wire1 plate}=[-, double=red!20!white]
\tikzstyle{d-wire2 plate}=[-, double=red!32!white]
\tikzstyle{dotted_plate}=[-, densely dotted, draw=red, fill opacity=0.4, fill=red!50!white, rounded corners]
\tikzstyle{twire1}=[-, draw, line width=0.4pt, preaction={-, draw, line width=1.4pt, red!20!white, preaction={-, draw, line width=2.2pt}}]
\tikzstyle{twire2}=[-, draw, line width=0.4pt, preaction={-, draw, line width=1.4pt, red!32!white, preaction={-, draw, line width=2.2pt}}]

\usepackage[colorlinks=true,linkcolor=blue,citecolor=blue,urlcolor=blue]{hyperref}
\usepackage[capitalize,noabbrev]{cleveref}
\usepackage[utf8]{inputenc}

\makeatletter
\def\Cref@thmoptarg[#1]#2#3#4{%
	    \ifhmode\unskip\unskip\par\fi%
	    \normalfont%
	    \trivlist%
	    \let\thmheadnl\relax%
	    \let\thm@swap\@gobble%
	    \thm@notefont{\fontseries\mddefault\upshape}%
	    \thm@headpunct{.}%
	    \thm@headsep 5\p@ plus\p@ minus\p@\relax%
	    \thm@space@setup%
	    #2%
	    \@topsep \thm@preskip
	    \@topsepadd \thm@postskip
	    \def\@tempa{#3}\ifx\@empty\@tempa%
	      \def\@tempa{\@oparg{\@begintheorem{#4}{}}[]}%
	    \else%
	      \refstepcounter[#1]{#3}%
	      \@namedef{Cref@#3@alias}{#1}%
	      \def\@tempa{\@oparg{\@begintheorem{#4}{\csname the#3\endcsname}}[]}%
	    \fi%
	    \@tempa}%
\makeatother

\Crefname{thm}{theorem}{theorems}
\Crefname{thm}{Theorem}{Theorems}
\Crefname{lem}{lemma}{lemmas}
\Crefname{lem}{Lemma}{Lemmas}
\Crefname{cor}{corollary}{corollaries}
\Crefname{cor}{Corollary}{Corollaries}
\Crefname{prop}{proposition}{propositions}
\Crefname{prop}{Proposition}{Propositions}
\Crefname{defi}{definition}{definitions}
\Crefname{defi}{Definition}{Definitions}
\Crefname{exa}{example}{examples}
\Crefname{exa}{Example}{Examples}
\Crefname{rem}{remark}{remarks}
\Crefname{rem}{Remark}{Remarks}
\Crefname{asm}{assumption}{assumptions}
\Crefname{asm}{Assumption}{Assumptions}
\Crefname{nota}{notation}{notations}
\Crefname{nota}{Notation}{Notations}

\newcommand{\N}{\mathbb{N}}

\newcommand{\Q}{\mathbb{Q}}
\newcommand{\R}{\mathbb{R}}
\newcommand{\eps}{\varepsilon}
\newcommand{\perms}[1]{S_{#1}}	% group of finite permutations on a set
\newcommand{\permact}[2]{{#2}^{#1}}	% having a permutation (1st argument) act on a Kolmogorov power of an object (2nd argument)
\let\P\undefined
\DeclarePairedDelimiterXPP{\P}[1]{\operatorname{\mathbf{P}}}[]{}{#1}
\DeclarePairedDelimiterXPP{\E}[1]{\operatorname{\mathbf{E}}}[]{}{#1}
\newcommand{\cat}[1]{{\mathsf{#1}}} 
\newcommand{\op}{\mathrm{op}}

\newcommand{\id}{\mathrm{id}} 		% identity
\newcommand{\tensor}{\otimes}
\newcommand{\Par}[1]{\mathsf{Partial}\left({#1}\right)}	% category of partial maps, as in \cite{cockettlack2002partialmaps}
\newcommand{\dom}[1]{\mathrm{dom}\left(#1\right)}		% domain

\tikzset{pullback/.style={minimum size=1.2ex,path picture={	% pullback symbol in diagram
			\draw[opacity=1,black,-,#1] (-0.5ex,-0.5ex) -- (0.5ex,-0.5ex) -- (0.5ex,0.5ex);%
}}}
\newcommand{\comp}{ 		% Command for sequential composition
	\mathchoice{\,}{\,}{}{} 	% First two are for displaystyle and text style, the remaining two for smaller math.
}
\newcommand{\domext}{\sqsupseteq}

\newcommand{\ph}{{\kern0.06em}\mathord{\rule[-0.05em]{0.6em}{0.05em}}{\kern0.06em}}		% Argument placeholder
\newcommand{\phsm}{{\kern0.04em}\mathord{\rule[-0.035em]{0.4em}{0.035em}}{\kern0.04em}}		% Scriptstyle argument placeholder
\newcommand{\cC}{\mathsf{C}}		% quasi-Markov cat
\newcommand{\cD}{\mathsf{D}}		% Markov cat
\renewcommand{\det}{\mathrm{det}}	% deterministic morphisms
\newcommand{\samp}{\mathsf{samp}}	% sampling map

\newcommand{\FinStoch}{\mathsf{FinStoch}}

\newcommand{\BorelStoch}{\mathsf{BorelStoch}}
\newcommand{\BorelMeas}{\mathsf{BorelMeas}}

\newcommand{\as}[1]{% 					almost surely
		\def\relstate{#1}%
		\ifx\relstate\empty
		  \text{a.s.}%
		\else
		  {#1\text{-a.s.}}%
		\fi
	}
\newcommand{\ase}[1]{=_{#1\text{-a.s.}}}					% almost surely equal
\newcommand{\es}{\mathsf{es}}			% empirical sampling map
\newcommand{\esav}{\overline{\es}}

\newcommand{\cop}{\mathrm{copy}}
\newcommand{\discard}{\mathrm{del}}

	\DeclarePairedDelimiter{\abs}{\lvert}{\rvert}
	\DeclarePairedDelimiter{\norm}{\lVert}{\rVert}
	\DeclarePairedDelimiterXPP{\pnorm}[2]{}{\lVert}{\rVert}{_{#1}}{#2}
	\makeatletter
		\let\oldabs\abs
		\def\abs{\@ifstar{\oldabs}{\oldabs*}}
		\let\oldnorm\norm
		\def\norm{\@ifstar{\oldnorm}{\oldnorm*}}
		\let\oldpnorm\pnorm
		\def\pnorm{\@ifstar{\oldpnorm}{\oldpnorm*}}
	\makeatother

	\providecommand{\given}{}			% Just to make sure the \given command exists.
	\newcommand{\SetSymbol}[1][]{%
		\nonscript\;\,#1\vert
		\allowbreak
		\nonscript\;\,
		\mathopen{}
	}
	\DeclarePairedDelimiterX{\Set}[1]{\{}{\}}{%
		\renewcommand{\given}{\SetSymbol[\delimsize]}
		#1
	}
	\makeatletter
		\let\oldSet\Set
		\def\Set{\@ifstar{\oldSet}{\oldSet*}}
	\makeatother

	\DeclarePairedDelimiterX{\Family}[1]{(}{)}{%
		\renewcommand{\given}{\SetSymbol[\delimsize]}
		#1
	}
	\makeatletter
		\let\oldFamily\Family
		\def\Family{\@ifstar{\oldFamily}{\oldFamily*}}

\newcommand{\newterm}[1]{\textbf{#1}}

\let\originalleft\left
\let\originalright\right
\renewcommand{\left}{\mathopen{}\mathclose\bgroup\originalleft}
\renewcommand{\right}{\aftergroup\egroup\originalright}

\title[Empirical Measures and Strong Laws in Categorical Probability]{Empirical Measures and \texorpdfstring{\\}{} Strong Laws of Large Numbers \texorpdfstring{\\}{} in Categorical Probability}

\author[T.~Fritz]{Tobias Fritz\lmcsorcid{0000-0002-7081-2635}}[a]
\author[T.~Gonda]{Tom\'{a}\v{s} Gonda\lmcsorcid{0000-0002-1531-0058}}[a]
\author[A.~Lorenzin]{Antonio Lorenzin\lmcsorcid{0000-0002-2415-4261}}[b]
\author[P.~Perrone]{Paolo Perrone\lmcsorcid{0000-0002-9123-9089}}[c]
\author[A.~Shah~Mohammed]{Areeb Shah Mohammed\lmcsorcid{0009-0003-6830-9098}}[a]

\address{Department of Mathematics, University of Innsbruck, Austria}
\email{tobias.fritz@uibk.ac.at, tomas.gonda@uibk.ac.at, areeb.shah-mohammed@uibk.ac.at}

\address{Independent researcher, Trento, Italy}
\email{a.lorenzin.95@gmail.com}

\address{Department of Computer Science, University of Oxford, United Kingdom}
\email{paolo.perrone.math@gmail.com}

\makeatletter
\def\shortauthors{T.~Fritz et al.}
\makeatother

\keywords{categorical probability, Markov categories, quasi-Markov categories, empirical measures, Glivenko--Cantelli theorem, strong law of large numbers, de Finetti theorem}
\amsclass{18M05, 60B10, 60F15, 60G09}

\begin{document}

\begin{abstract}
	The Glivenko--Cantelli theorem is a uniform version of the strong law of large numbers.
	It states that for every IID sequence of random variables, the empirical measure converges to the underlying distribution (in the sense of uniform convergence of the CDF).
	In this work, we provide tools to study such limits of empirical measures in categorical probability.
	We propose two axioms, namely permutation invariance and empirical adequacy, that a morphism of type $X^\N \to X$ should satisfy to be interpretable as taking an infinite sequence as input and producing a sample from its empirical measure as output.
	Since not all sequences have a well-defined empirical measure, such \emph{empirical sampling morphisms} live in quasi-Markov categories, which, unlike Markov categories, allow for partial morphisms.
	
	Given an empirical sampling morphism and a few other properties, we prove representability
		as well as abstract versions of the de Finetti theorem, the Glivenko--Cantelli theorem and the strong law of large numbers.
	
	We provide several concrete constructions of empirical sampling morphisms as partially defined Markov kernels on standard Borel spaces.
	Instantiating our abstract results then recovers the standard Glivenko--Cantelli theorem and the strong law of large numbers for random variables with finite first moment.
	Our work thus provides a joint proof of these two theorems in conjunction with the de Finetti theorem from first principles.
\end{abstract}

\maketitle
\tableofcontents

\section{Introduction}

Laws of large numbers are widely regarded among the most important results in probability theory both in technical and conceptual terms.
Conceptually, they provide self-consistency to the idea of a probability measure: 
Although the actual expectation value of a function $f$ under an unknown distribution $\mu$ cannot be inferred from a finite sequence of samples from $\mu$, its \newterm{empirical averages} converge with probability $1$ to this expectation value as the number of observations grows.
In formulas, we have
\begin{equation}
	\label{eq:LLN_intro}
	\lim_{n \to \infty} \frac{f(x_1) + \dots + f(x_n)}{n} \;\ase{\mu}\; \int f(x) \, \mu(\mathrm{d}x),
\end{equation}
where $\ase{\mu}$ indicates equality $\mu$-almost surely.
This is why, given a finite number of samples $(x_1, \dots, x_n)$, one can reasonably approximate the expectation value on the right by the empirical average, which is the fraction on the left.

In particular, consider $f$ to be the indicator function of a measurable set $T \subseteq X$.
We obtain that the relative frequency of the event $T$ in the sequence $(x_i)_{i \in \N}$ converges almost surely to its probability:
\begin{equation}
	\label{eq:rel_freq_conv}
	\lim_{n \to \infty} \frac{|\{ i \le n \mid x_i \in T \}|}{n} \;\ase{\mu}\; \mu(T).
\end{equation}
As $T$ varies, one can interpret the fraction on the left as the probability of $T$ under the \newterm{empirical measure} 
\begin{equation}
	\label{eq:finite_empirical_distribution}
	\frac1n \left( \delta_{x_1} + \dots + \delta_{x_n} \right) 
\end{equation}
of the first $n$ elements of the sequence, where $\delta_{x_i}$ stands for the Dirac measure at $x_i$.
One might therefore expect that the sequence of empirical measures in some sense converges to $\mu$.
The details of this convergence are subtle, and we will see that there are several reasons why it may fail.
However, there also are many positive results that avoid these pitfalls and provide conditions that ensure the desired convergence~\cite{vaartwellner1996empiricalprocesses}.
One such result that is of primary interest to us is the \newterm{Glivenko--Cantelli theorem}.
It applies to $x_i$ valued in $\R$ and establishes the uniform convergence of the relative frequencies of~\Cref{eq:rel_freq_conv} as $T$ ranges over all intervals in $\R$.

We regard such theorems on the convergence of empirical measures as the conceptually most fundamental laws of large numbers.
In this vein, we can think of a standard strong law of large numbers as a consequence obtained by taking the expectation value of $f$ on both sides of \Cref{eq:rel_freq_conv}.\footnote{However, making this precise for unbounded $f$ requires dealing with a number of subtleties, which we return to in \Cref{sec:es_integration}.}

The main goal of this paper is to develop results like this in the framework of \newterm{categorical probability theory}.
This recent approach to probability theory
strives to redevelop the main results of probability theory in a structured manner.
It is based on a more abstract starting point than measure theory.
Namely, the idea is to capture the essential features of reasoning in probability theory as properties of certain categories (often\ \newterm{Markov categories}) that generalize the category of measurable spaces and Markov kernels.
One benefit of this higher-level language is its focus on fundamental properties while leaving irrelevant details aside.

In order to develop laws of large numbers in this framework, we need to express limits mentioned above in the categorical language.
One approach would be to equip the hom-sets with a topology and to consider convergence in this topology~\cite{perrone2024markov}.
Here, we pursue a more abstract approach instead:

\begin{itemize}
	\item For every standard Borel space $X$, we construct a Markov kernel $X^\N \to X$ that takes an \emph{infinite} sequence $(x_i)_{i \in \N}$ as input and returns a sample from its empirical measure as output.
		Since the limit of empirical measures (of the first $n$ elements) as $n \to \infty$ need not exist, as we discuss on the next page,
		this kernel is only \emph{partially defined}.
		Instead of specifying a topology that would determine when these limits converge, we thus have the notion of the domain of the Markov kernel that tells us which infinite sequences have a well-defined empirical measure.

	\item In Markov categories, we \emph{axiomatize} the properties that a partial morphism ${\es_X : X^\N \to X}$ should satisfy in order to carry the interpretation of an \newterm{empirical sampling morphism}.
		Our two axioms are permutation invariance and empirical adequacy.
		The former states that $\es$ should be invariant under finite permutations of the input sequence.
		Empirical adequacy roughly states that sampling from an exchangeable distribution $\mu$ on $X^\N$ is the same thing as iterated empirical sampling applied to a sampled sequence.

\end{itemize}
Permutation invariance and empirical adequacy are not only important properties to give $\es_X$ the intended interpretation, but they are also useful in synthetic proofs of results involving empirical measures.

Given a Markov category\footnote{Technically, we need a positive quasi-Markov category in which balanced idempotents split.} with empirical sampling morphisms, we derive abstract versions of the de Finetti theorem (\Cref{cor:dF}), 
the Glivenko--Cantelli theorem (\Cref{thm:lln}) and the strong law of large numbers (\Cref{cor:genericLLN}).
Another one of our key technical contributions is a constructions of concrete empirical sampling morphisms in measure-theoretic probability and, in particular, proofs that they satisfy the two axioms mentioned above.
To set the stage for these developments, we also prove a number of technical results on Kolmogorov products in quasi-Markov categories and on the passage to partial morphisms, which are of independent interest.
Particularly relevant here is the assumption of $\sigma$-continuity (\Cref{def:count_meets}), by which we mean the existence of countable directed meets in the hom-sets and their preservation by composition and tensor. 
In measure-theoretic probability, this turns out to be related to the fact that probability measures are $\sigma$-continuous as functions on the $\sigma$-algebra (\Cref{prop:sigma_parborelstoch}).

Together, our abstract theorems and the concrete constructions recover the de Finetti theorem for standard Borel spaces~\cite{fritz2021definetti}, the Glivenko--Cantelli theorem (\Cref{thm:glivenko}) and the strong law of large numbers for real-valued random variables with finite first moment (\Cref{cor:strong_lln}).

Our framework thus provides a structured proof of these three fundamental results from first principles.
We believe that this joint proof is short compared to the traditional purely measure-theoretic treatment, and moreover displays enhanced conceptual clarity.

Translating probabilistic concepts and results into a categorical framework is often not straightforward,
	but it can be very valuable.
The above summary demonstrates that empirical sampling morphisms are an invaluable tool when translating laws of large numbers and related theorems.
Furthermore, we suggest that they are a key concept in probability theory that has been underappreciated so far.
We hope that this realization will facilitate further progress on categorical probability theory.
For example, we expect empirical sampling morphisms to be useful for the development of ergodic theorems in the categorical framework.

The rest of the Introduction provides a more detailed overview of our results.

\paragraph{Empirical sampling morphisms in measure-theoretic probability.}

For a finite standard Borel space $F$, we can construct an empirical sampling morphism (in this case, a partial Markov kernel) according to the intuition above.
The transition probability of an event $T \subseteq F$ given a sequence $(x_i)_{i \in \N}$ is the limiting relative frequency of the occurence of the event within the sequence, i.e.\
\begin{equation}\label{eq:desired_es}
	\es_F \bigl( T \,|\, (x_i) \bigr) \;\coloneqq\; \lim_{n\to\infty} \frac{ \abs{ \Set{ i \le n \given x_i \in T } } }{n}
\end{equation}
whenever the limit exists.

There are three types of subtleties with this construction, and a substantial part of our work is devoted to addressing them.
\begin{enumerate}
	\item Indeed, not all sequences $(x_i)$ admit a limit as in \Cref{eq:desired_es} above.  
	For example, consider $X = \{0,1\}$ and take $x_i = 0$ whenever $\lceil \log_2{i} \rceil$ is even and $x_i = 1$ otherwise.
	Then the relative frequency of either outcome oscillates indefinitely and does not have a limit.
\end{enumerate}
One could try to get around this by forcing every sequence to have well-defined limits, for example by using ultrafilter convergence or Banach limits (with respect to which every bounded sequence has a limit). 
However, this is not what we want, because we would not recover the convergence in the $\eps$-$\delta$ sense as part of the conclusion of our results.
Therefore, we formalize empirical sampling as a \emph{partial} Markov kernel whose domain is a subset of $X^\N$. 
The question of which sequences should be included in the domain is itself subtle and mathematically interesting.
There is a certain flexibility in the choice of domain which is closely related to the assumption of finite first moment in the strong law of large numbers (see \Cref{sec:es_integration}).

\begin{enumerate}[resume]
	\item For infinite $X$, even if the limits exist, they may not define a probability measure $\es_X(\ph|(x_i))$. 
	For example, the sequence $(1,2,3,\dots)$ over $X = \N$ has limiting relative frequency \emph{zero} for each singleton subset $T$ of $\N$, so it does not satisfy $\sigma$-additivity.
\end{enumerate}
To address this second issue, we require the limits in~\Cref{eq:desired_es} to be uniform on sets of the form $T = \{1,\ldots,t\}$,\footnote{In this article, we use $\N = \{1,2,3,\ldots\}$ as the conventionally chosen countably infinite set.\label{foot:N}} and take $\es_\N$ to be undefined on all sequences that do not satisfy this requirement (see \Cref{prop:eq_uniform} for a classification of such sequences).

\begin{enumerate}[resume]
		\item For uncountable $X$, we cannot expect \Cref{eq:desired_es} to hold for every measurable set $T \subseteq X$.
		Requiring this would entail no sequence with mutually distinct elements has an empirical measure, which is clearly undesirable.
		
		Indeed if $(x_i)$ is such a sequence and \Cref{eq:desired_es} is postulated for all measurable sets, then the empirical measure will assign probability $1$ to the (countable) support of the sequence, but probability $0$ to every individual point in the support, contradicting $\sigma$-additivity.
\end{enumerate}
Therefore for uncountable $X$, we need to restrict the validity of \Cref{eq:desired_es} to a certain class of measurable sets.
We do this in \Cref{sec:es_BorelStoch} for $X = \R$ by restricting to intervals and requiring the uniform existence of the limit on these.

As already indicated, we identify two structural properties satisfied by the above constructions that serve as axioms to define an abstract empirical sampling morphism.
\begin{itemize}
	\item \emph{Permutation invariance:} For every measurable $T\subseteq X$, the probability $\es_X (T|(x_i))$ is invariant under finite permutations of the sequence $(x_i)$.
		
		This encapsulates the fact that the relative frequencies do not depend on the order of outcomes appearing in the sequence $(x_i)$, provided that one keeps almost all of them fixed.
	\item \emph{Empirical adequacy:}\footnote{The formal empirical adequacy axiom in \Cref{def:es} is slightly stronger than this, but here we focus on the case $Y = I$ for simplicity, which contains the essence of the idea.}
			If $\mu$ is an exchangeable probability measure on $X^\N$, then we have
	\begin{equation}\label{eq:invariance}
		\int_{X^n} \es(T_1|(x_i)) \, \cdots \, \es(T_n|(x_i)) \, \mu_{1\ldots n}(\mathrm{d}x_1 \times \dots \times \mathrm{d}x_n) \;=\; \mu_{1\ldots n} \bigl( T_1\times\dots\times T_n  \bigr),
	\end{equation}
		where $\mu_{1\ldots n}$ denotes the marginal of $\mu$ on the first $n$ factors.
		For $\mu = \mu_1^{\otimes \N}$ a product measure and $n = 1$, this says that applying empirical sampling to a sequence of IID samples from $\mu_1$ produces another sample from $\mu_1$ as expected.

		For our particular constructions in measure-theoretic probability, this property follows from \Cref{eq:es_random_permutation}, where we show that the \emph{resampling} kernel $X^\N \to X^\N$ defined as on the left of \Cref{eq:invariance},
	\begin{equation}
		\mathsf{resamp} \bigl( T_1\times\dots\times T_n \times X \times \dots \,|\, (x_i) \bigr) \; \coloneqq \; \es(T_1|(x_i)) \, \cdots \, \es(T_n|(x_i)),
	\end{equation}
	can be written as the application of a uniformly random permutation to the sequence elements.
\end{itemize}
\noindent %FIXED indent
As we then show in \Cref{sec:theorems}, the existence of a kernel with these two properties features as the main common ingredient in proofs of the de Finetti theorem, the Glivenko--Cantelli theorem and the strong law of large numbers.

\paragraph{The main ideas, category-theoretically.}

All of the above takes a simpler and more intuitive form in the categorical framework, which we now briefly recap.

\emph{Markov categories}~\cite{chojacobs2019strings,fritz2019synthetic} are an abstract generalization of the category of measurable spaces and Markov kernels.
Research on categorical probability via Markov categories rests on the idea that the main concepts of probability theory, such as statistical (in)dependence, determinism, conditioning, invariance under symmetries, etc., can be abstractly generalized from categories of Markov kernels to more general Markov categories.
Many of these notions can be captured in terms of universal properties.
To illustrate this, consider the Markov category $\BorelStoch$, which has standard Borel spaces as objects and Markov kernels between them as morphisms and is the Markov category of primary interest for measure-theoretic probability.
Its subcategory of deterministic morphisms is $\BorelMeas$, consisting of standard Borel spaces and measurable functions between them.
Now given a standard Borel space $X$, we can characterize the measurable space $PX$ of probability measures on $X$~\cite{giry} in terms of a universal property: 
It is equipped with a natural bijection
\begin{equation}
	\BorelStoch(A,X) \;\cong\; \BorelMeas(A,PX)
\end{equation}
between Markov kernels $A \to X$ and measurable functions $A \to PX$ defined for every standard Borel space $A$ \cite{fritz2023representable}.
Additionally, the de Finetti theorem can be formulated as the fact that the space $PX$ is the limit{\,---\,}in the categorical sense{\,---\,}of the diagram of all finite permutations acting on the space $X^\N$ \cite{fritz2021definetti,moss2022probability}. 

Abstractions like these allow us to formulate and prove theorems of probability and related fields by means of diagram manipulation, where measure theory is needed only for showing that the relevant axioms are satisfied for the Markov category $\BorelStoch$.
Results that have been developed like this over the past few years feature the Kolmogorov and Hewitt--Savage zero-one laws~\cite{fritzrischel2019zeroone}, the de Finetti theorem~\cite{fritz2021definetti,moss2022probability}, the d-separation theorem for graphical models~\cite{fritz2022dseparation}, the ergodic decomposition theorem~\cite{moss2022ergodic,ensarguet2023ergodic}, the Blackwell--Sherman--Stein theorem~\cite{fritz2023representable} and the Aldous--Hoover theorem~\cite{chen2024aldoushoover}.\footnote{This fresh perspective has also gained popularity in the computer science community, as it allows for a more systematic treatment of the logic underlying probability theory~\cite{jacobs2020logical,stein2021structural} as well as applications to topics including active inference~\cite{tull2024active}, random graphs \cite{ackerman2024randomgraphs}, causal inference~\cite{jacobs2019causal_surgery} and evidential decision theory~\cite{dilavore2023evidential}, to name a few.}
In addition to the categorical proofs often being simpler and more intuitive than the measure-theoretic ones, they also allow for greater generality, as the categorical results can often be instantiated in Markov categories that model other kinds of uncertainty than measure-theoretic probability.

In this work we follow the same approach, focusing on the study of \emph{empirical sampling} as introduced above.
The axioms on empirical sampling morphisms, in the form of the permutation invariance and empirical adequacy discussed above, have simple and natural categorical formulations (\Cref{def:es}).

As already indicated, empirical sampling morphisms are \emph{partial} maps that need not be defined for every input sequence. 
Formulating this in categorical probability therefore requires going beyond the standard theory of Markov categories.
Notions of partiality for Markov categories
have been first given in \cite{dilavore2023evidential,dilavore2024partial,dilavore2025OrderPartialMarkov}.
For our purposes, we assume quasi-totality \cite{dilavore2023evidential} and don't need conditionals, and thus our approach differs from the earlier ones.
To this end, we develop the theory of \newterm{quasi-Markov categories} in \Cref{sec:markovcats}, where the main technical contribution is a functorial version of Kolmogorov products in this setting.
This framework is adequate for our purposes and general enough to include any category constructed from any ``partializable'' Markov category $\cC$ by adding formal partial versions of morphism from $\cC$, as shown by one of the authors in \cite{shahmohammed2025partial}.

\paragraph{Related work.}

In measure-theoretic probability, there is ample literature on empirical processes, including textbooks such as~\cite{shorack1986empirical,vaartwellner1996empiricalprocesses}.
This literature generally focuses on the empirical measure of a finite sequence $(x_1,\ldots,x_n)$, and then studies this as a stochastic process as $n$ varies, proving various results on when this process converges to the underlying distribution.

On the other hand, we are not aware of any prior literature on the question of when a general infinite sequence $(x_i)_{i \in \N}$ has a well-defined empirical measure.
There seems to be remarkably little literature on empirical measures of infinite sequences, even in the measure-theoretic setting.
One reference where these have been considered, in a similar spirit as in the present manuscript, is in a 2014 paper by Austin and Panchenko~\cite{austin2014hierarchical}.
We will comment on this further in \Cref{ex:es_domain}.

\paragraph{Brief outline.}
In \Cref{sec:markovcats}, we introduce our framework of quasi-Markov categories and its features.
\Cref{sec:emp_samp} then concerns empirical sampling morphisms.
We first give the categorical definition in \Cref{sec:es_def} and subsequently construct concrete instances thereof in measure-theoretic probability.
Proving that our constructions satisfy the definition is somewhat technical, which is why we include the proofs in the appendix (\Cref{sec:measure_theory,sec:proofs}).
Finally, in \Cref{sec:theorems} we prove our main theorems, most notably the synthetic Glivenko--Cantelli theorem, and show how they recover the standard measure-theoretic results.

\section*{Acknowledgements}

We thank Alexander Glazman, Andreas Klingler and Christian Wei\ss{} for helpful discussions.
Research for Paolo Perrone is funded, at the time of writing, by Sam Staton's Consolidator Grant ``BLaSt -- a Better Language for Statistics'' from the European Research Council.
Tobias Fritz, Antonio Lorenzin and Areeb Shah Mohammed acknowledge funding by the Austrian Science Fund (FWF) [doi:\href{https://www.doi.org/10.55776/P35992}{10.55776/P35992}].
Tomáš Gonda has been funded in whole or in part by the Austrian Science Fund (FWF) 10.55776/ESP3451824 and the Start Prize Y1261-N.
Areeb Shah Mohammed has also been supported by the Doctoral Scholarship of the University of Innsbruck.
Additionally, Antonio Lorenzin has received support from the ARIA Safeguarded AI TA1.1 programme.

\section{Quasi-Markov Categories}
\label{sec:markovcats}

Before delving into our main contributions, the study of empirical sampling morphisms (\Cref{sec:emp_samp}) and synthetic laws of large numbers (\Cref{sec:theorems}), let us set the stage for this discussion by introducing quasi-Markov categories (\Cref{def:quasi-total}).
In \Cref{sec:basic_definitions} we also introduce our key instance that supports measure-theoretic interpretation (\Cref{ex:par_borelstoch}).
More general construction of quasi-Markov categories through which one can obtain further examples can be found in \cite{shahmohammed2025partial}.

In \Cref{sec:qM_defs}, we spell out a few basic definitions that are familiar to readers accustomed to the theory of Markov categories.
Unlike Markov categories, quasi-Markov ones come with a canonical (and non-trivial) ordering among morphisms that reminds one of other types of categories to model partiality \cite{cockettlack2002partialmaps}.
We discuss this partial order in \Cref{sec:poset_enrichment} and use it to develop tools for working with infinite tensor products in quasi-Markov categories in \Cref{sec:products}, which were initially introduced for Markov categories in~\cite{fritzrischel2019zeroone}.
\Cref{sec:representable} then introduces distribution objects in quasi-Markov categories.
Finally, following the arguments laid out in \cite{fritz2023involutive}, we show how a specific equalizer of permutations of a countable sequence of objects (a de Finetti object) is in fact also a distribution object (\Cref{thm:defin_obs}).

\subsection{Basic Definitions}
\label{sec:basic_definitions}

We first sketch the definition of CD categories, of which quasi-Markov categories are a special case. 
For more details, we refer the reader to \cite{chojacobs2019strings,fritz2019synthetic}. 

\begin{defi}
 A \newterm{copy-discard} (\newterm{CD}) \newterm{category} is a symmetric monoidal category in which every object $X$ is equipped with a distinguished commutative comonoid structure, compatible with the tensor product. 
 The comonoid structure maps
 \begin{equation}
	 \cop_X : X\to X\tensor X, \qquad \discard_X : X\to I
 \end{equation}
 are called \newterm{copy} and \newterm{delete}, respectively.
\end{defi}

It is helpful to draw morphisms in CD categories in terms of string diagrams. 
In our conventions, we draw diagrams from bottom to top, and copy and delete like this:
\begin{equation} \input{comultiplication.tikz} \qquad \qquad \begin{tikzpicture}
	\begin{pgfonlayer}{nodelayer}
		\node [style=bn] (0) at (0, 0.5) {};
		\node [style=none] (1) at (0, -0.5) {};
		\node [style=none] (2) at (0, -1) {$X$};
		\node [style=none] (3) at (-3, -0) {$\discard_{X} \quad =$};
	\end{pgfonlayer}
	\begin{pgfonlayer}{edgelayer}
		\draw (0) to (1.center);
	\end{pgfonlayer}
\end{tikzpicture} \end{equation}
Morphisms from the monoidal unit $I$, referred to as \newterm{states},
	are represented by triangles:
\begin{equation} \begin{tikzpicture}
	\begin{pgfonlayer}{nodelayer}
		\node [style=none] (2) at (0, 1) {};
		\node [style=none] (4) at (0, 1.5) {$X$};
		\node [style=state] (5) at (0, 0) {$m$};
	\end{pgfonlayer}
	\begin{pgfonlayer}{edgelayer}
		\draw (5) to (2.center);
	\end{pgfonlayer}
\end{tikzpicture}
 \end{equation}

\begin{defi}\label{def:CDTotal}
	A morphism $f : A \to X$ in a CD category is \newterm{total} if it satisfies
	\begin{equation}\label{equation:CDTotalDef}
		\input{CDTotalDef.tikz}
	\end{equation}
	A CD category is a \newterm{Markov category} if every morphism is total, or equivalently if the monoidal unit is terminal.
\end{defi}
	
	We can interpret Markov categories as models of information flow that may involve randomness or nondeterminism.
	A morphism $X \to Y$ is then viewed as a (potentially noisy) channel from $X$ to $Y$. 
	The copy map duplicates the information perfectly, and the delete map discards it. 

	Here are two of the most prominent examples of Markov categories.

	\begin{exa}[A Markov category for discrete probability]\label{ex:finstoch}
	The category $\FinStoch$ has:
	\begin{itemize}
	\item As objects, finite sets;
	\item As morphisms $X\to Y$, stochastic matrices with entries denoted by $f(y|x)$ and the usual composition: 
	For $f : X\to Y$ and $g  :Y\to Z$, we have
	\begin{equation} (g\comp f)(z|x) = \sum_{y\in Y} g(z|y)\cdot f(y|x) ; \end{equation}
	\item The tensor product is given by the cartesian product of finite sets and the tensor (or Kronecker) product of matrices, with entries given by $f(x|a)\cdot g(y|b)$ for $f : A\to X$ and $g : B\to Y$;
	\item The copy and delete maps are given by the following matrices:
	\begin{equation}
		\cop_X (x_1,x_2|x) = \begin{cases}
								1 & \mbox{if } x_1=x_2=x ;\\
								0 & \mbox{else};
							\end{cases}
		\qquad \qquad
		\discard_X( \ast | x ) = 1;
	\end{equation}
	where the monoidal unit $I$ is the singleton set $\{\ast\}$.
	\end{itemize}
	This category is the prototype for discrete probability theory.
	\end{exa}

	%FIXED IMPORTANT Swapped the position of ''$\BorelStoch$'' and ''category'' for margins, NOTE PLEASE REVIEW IF IT CHANGES MEANING AND LET ME KNOW
	\begin{exa}[A Markov category for measure-theoretic probability]\label{ex:borelstoch}
	The $\BorelStoch$ category  has:
	\begin{itemize}
	\item As objects, standard Borel spaces, i.e.\ measurable spaces whose $\sigma$-algebra can be expressed as the Borel $\sigma$-algebra of a Polish space;
	\item As morphisms $X\to Y$, Markov kernels, with their usual (Chapman--Kolmogorov) composition: 
	For $f : X\to Y$ and $g : Y\to Z$, and given a measurable set $A \subseteq Z$, we have
	\begin{equation}
		\label{eq:ChapmanKolmogorov}
		(g\comp f)(A|x) = \int_Y g(A|y)\cdot f(\mathrm{d}y|x) ;
	\end{equation}
	
	\item The tensor product is given by the cartesian product of measurable spaces and the product of measures;
	
	\item The copy and delete maps are given by the following matrices:
	\begin{equation}
		\cop_X ( \ph | x) = \delta_{(x,x)} 
		\qquad\qquad 
		\discard_X(\{\ast\}| x) = 1.
	\end{equation}
	\end{itemize}
	This category is the prototype for measure-theoretic probability theory.
	It contains $\FinStoch$ as a full subcategory.
	\end{exa}
	
	In order to model partial Markov kernels, we need to relax \Cref{equation:CDTotalDef}, since it forces a Markov kernel $f$ to produce an output for each input.
	Partial Markov kernels satisfy the following weaker condition, which will play an important role in our proofs.

	\begin{defiC}[{\cite[Definition~3.1]{dilavore2023evidential}}] %FIXED switched to defiC
		\label{def:quasi-total}
		A morphism $f : A \to X$ in a CD category is \newterm{quasi-total} if it satisfies
		\begin{equation}\label{eq:quasi-total}
			\input{quasi-total.tikz}
		\end{equation}
		A \newterm{quasi-Markov category} is then a CD category in which every morphism is quasi-total. 
	\end{defiC}
	
	Trivially, every total morphism is quasi-total and every Markov category is quasi-Markov.
	We can interpret quasi-Markov categories as a variant of Markov categories in which a morphism is allowed to fail, i.e.\ it may not produce any output for some of its input values.
	This failure is, however, deterministic in a sense that we illustrate for $\cat{FinSubStoch}$, the variant of $\FinStoch$ without the normalization condition.
	A morphism $f : X \to Y$ in $\cat{FinSubStoch}$ (i.e.\ a substochastic matrix) is quasi-total if and only if for every input value $a \in A$, we have either
	\begin{equation}
		f(x | a) = 0 \quad \forall x \in X
	\end{equation}
	or
	\begin{equation}
		\sum_{x \in X} f(x | a) = 1.
	\end{equation}
	In the former case, we say that $f$ is undefined on input $a$, i.e.\ $f$ does not produce any output.
	In the latter case, $a$ is said to be in the domain of $f$ (cf.\ \Cref{def:domain}), i.e.\ $f$ applied to $a$ produces an output (with probability 1).
	
	The quasi-totality condition therefore forces $f$ to be normalized for every input where it does not fail.
	Contrast this with the totality condition, which requires, in addition, that $f$ produces an output for every input, i.e.\ that $f$ is a stochastic matrix.

	Both $\FinStoch$ and $\BorelStoch$ are Markov categories, i.e.\ all of their morphisms are total.
	For each, we can consider a corresponding quasi-Markov category which additionally includes partial morphisms.
	Here, let us merely introduce the version of $\BorelStoch$ with partial morphisms added.
	This is the key example for us, in which we will instantiate our categorical results from \Cref{sec:theorems}.

	\begin{exa}[A quasi-Markov category for measure-theoretic probability]\label{ex:par_borelstoch}
		Specifically, we denote the resulting category $\Par{\BorelStoch}$ and call it the \newterm{partialization} of $\BorelStoch$~\cite{shahmohammed2025partial}. It has:
		\begin{itemize}
			\item As objects, standard Borel spaces;
			\item As morphisms $X\to Y$, partial Markov kernels, which are pairs $(D_f,f)$ consisting of a measurable set $D_f \subseteq X$ called \newterm{domain} and a Markov kernel $f : D_f \to Y$.
		\end{itemize}
		We often write $f : X \to Y$ by abuse of notation for both the partial Markov kernel with domain $D_f$ and the ordinary Markov kernel $f : D_f \to Y$.

		\begin{itemize}
			\item Composition of $f : X \to Y$ and $g : Y \to Z$ is $g \comp f : X \to Z$ with domain
				\begin{equation}\label{eq:dom_comp}
					D_{g \comp f} \coloneqq \Set{ x \in D_f \given f(D_g|x) = 1 }.
				\end{equation}
				and with the Markov kernel $g \comp f : D_{g \comp f} \to Z$ defined as in~\Cref{eq:ChapmanKolmogorov}.\footnote{As an alternative to $f(\ph|x)$ being undefined for $x \notin D_f$, we could also set $f(T|x) = 0$ for all measurable sets $T \subseteq Y$.
				However, this should be used with caution, since we cannot use the Chapman--Kolmogorov \Cref{eq:ChapmanKolmogorov} for those $x \in D_f$ that lie outside $D_{g \comp f}$ as given by \Cref{eq:dom_comp}. 
				Namely if $x \in D_f$ and we have $0 < f(D_g|x) < 1$, then Chapman--Kolmogorov equation would give $0 < (g \comp f)(Z|x) < 1$ which contradicts quasi-totality. 
				In other words, it is important that our composite fails either if the first part $f$ fails or if the second part $g$ fails with non-zero probability.}
			\item The tensor product of $f : A \to X$ and $g : B \to Y$ is given by the partial Markov kernel ${f \tensor g : A \times B \to X \times Y}$ with domain
				\begin{equation}
					D_{f \tensor g} = (D_f \times B) \cap (A \times D_g),
				\end{equation}
				and defined on this domain as the product measure, just as in $\BorelStoch$.
			\item Copy and delete morphisms are as in $\BorelStoch$.
		\end{itemize}
		\noindent %FIXED indent
		We leave it to the reader to verify that this indeed forms a quasi-Markov category, either by direct verification or by noting that this category is the partialization of $\BorelStoch$ in the sense of~\cite{shahmohammed2025partial}.
	\end{exa}

	\begin{rem}[Disclaimer on the use of ``partialization'']
		The term ``partialization'' used above follows a tradition in category theory, where it refers to categories whose morphisms are spans with monic left legs \cite{cockettlack2002partialmaps}.
		Indeed, in $\Par{\BorelStoch}$, morphisms can be represented as spans
		\[
			X \hookleftarrow D_f \xrightarrow{f} Y,
		\]
		where the left leg is a deterministic monomorphism (see \cite{shahmohammed2025partial} for further details on this connection).

		In contrast, in category-theoretic approaches to probability, the term ``partial Markov categories'' has been used to refer to categories such as $\cat{FinSubStoch}$ \cite{dilavore2023evidential}, where not all morphisms are quasi-total. 
		This also motivates our use of the term ``quasi-Markov categories'', which avoids a clash of terminology with the existing literature.
		Partial Markov categories will not be discussed further in this paper.
	\end{rem}

	An important feature of quasi-Markov categories is the following property, which will be particularly relevant in the discussion of representability (\Cref{sec:representable,sec:dF_obs_rep}).
	\begin{lem}\label{lem:mono_total}
		Every monomorphism $m : A \to X$ in a quasi-Markov category is total.
	\end{lem}
	\begin{proof}
		Under these assumptions, we have
		\begin{equation}
			\input{quasi-total-mono.tikz}\qquad \implies\qquad \input{quasi-total-mono-cancelled.tikz}\qquad \implies\qquad \input{quasi-total-mono-total.tikz}
		\end{equation}
		where the first implication holds because $m$ is a monomorphism and the second follows by applying $\discard_X$ to the output.
		As a result, quasi-totality of $m$ (the antecedent) implies its totality (the final consequent).
	\end{proof}

\subsection{Determinism, Almost Sure Equality and Positivity}
\label{sec:qM_defs}

To get a suitable language for discussing probability, we need some natural generalizations of existing notions for Markov categories. 
Namely, we define what it means for a morphism to be deterministic, almost sure equalities, and the positivity axiom.

\begin{defiC}[{\cite[Section 2]{carboniwalters1987bicartesian}}]\label{def:copyable} %FIXED defC
	Let $\cC$ be a quasi-Markov category.
	\begin{enumerate}
		\item A morphism $f : X\to Y$ in $\cC$ is \newterm{copyable} if it commutes with copying:
 			\begin{equation}\label{eq:copyable}
 				\input{multiplication_natural.tikz}
 			\end{equation}
 			$f$ is \newterm{deterministic} if it is both total and copyable. 
 		\item The wide subcategory of $\cC$ consisting of copyable morphisms is denoted by $\cC_{\rm cop}$.
	\end{enumerate}
\end{defiC}
\noindent %FIXED indent
Clearly all copyable morphisms are quasi-total.
In a Markov category, a morphism is deterministic if and only if it is copyable.
We think of a copyable $f$ as a channel which gives the same result (with probability $1$) when applied twice to the same input. 

\begin{exa}[Copyable morphisms in concrete categories]
In many examples, copyable morphisms correspond indeed to what one would expect.
 \begin{itemize}
  \item In $\FinStoch$, copyable morphisms are precisely the deterministic stochastic matrices, i.e.\ those that contain only zero and one entries;
  \item In $\BorelStoch$, a Markov kernel $f : X\to Y$ is copyable if and only if $f(T|x) \in \{0,1\}$ for all ${x \in X}$ and all measurable $T \subseteq Y$.
	  These are precisely those Markov kernels of the form ${f(\ph ,x)=\delta_{g(x)}}$, where $g : X\to Y$ is a measurable function (necessarily unique).
	  Thus $\BorelStoch_{\rm cop}$ is isomorphic to $\BorelMeas$, the category of standard Borel spaces and measurable functions between them.
  \item In $\Par{\BorelStoch}$, the copyable morphisms are those partial Markov kernels that are deterministic on their domain, i.e.~of the form $f(\ph|x) = \delta_{g(x)}$ for some measurable $g : D_f \to Y$.
  In particular, $\Par{\BorelStoch}_{\rm cop}$ coincides with $\Par{\BorelMeas}$ (this is a general fact of partialization of Markov categories, see \cite[Proposition 3.22]{shahmohammed2025partial}).
	  Such an $f$ is deterministic if $D_f = X$ in addition.
 \end{itemize}
\end{exa}
\noindent %FIXED indent
Throughout this paper, \emph{almost sure equalities} play an important role. 
The definition is the same as for Markov categories~\cite{fritz2019synthetic,fritz2023supports} or more generally CD categories \cite{chojacobs2019strings}. 

\begin{defi}
\label{def:ase}
In a quasi-Markov category $\cC$, consider a morphism $p : A \to X$ and two  parallel morphisms $f,g : X \to Y$. 
We say that $f$ and $g$ are \newterm{\as{$\bm{p}$} equal}, and write $f \ase{p} g$, if the following equality holds:
\begin{equation}\label{eq:as}
	\input{as_eq.tikz}
\end{equation}
\end{defi}

The following additional axiom is useful in the context of Markov categories~\cite[Definition~11.22]{fritz2019synthetic}.
Its name is motivated by the fact that its satisfaction is related to the non-negativity of measures~\cite[Example~11.25]{fritz2019synthetic}.

\begin{defi}\label{def:pos}
	A quasi-Markov category is \newterm{positive} if for each $f : X\to Y$ and $g : Y\to Z$, whose composite $g\comp f : X\to Z$ is copyable, the following equation holds:
	\begin{equation}\label{positivity}
		\input{positivity.tikz}
	\end{equation}
\end{defi}
To get an intuition, consider the case of $X$ being the monoidal unit $I$.
Then positivity means that coarse-graining a state $m : I\to Y$ to a copyable state $g \comp m : I\to Z$ makes $Y$ and $Z$ necessarily independent in the sense that 
\begin{equation}
	\input{positivity_unit.tikz}
\end{equation}
holds.
Thus the positivity axiom can be also seen as a requirement that a deterministic variable cannot display correlation with other variables~\cite[Section~2.2]{fritz2022dilations}.

Our main example, $\Par{\BorelStoch}$, is a positive quasi-Markov category.
One can prove this by adapting the proof for $\BorelStoch$ \cite[Example 11.25]{fritz2019synthetic} or by using the general results of \cite{shahmohammed2025partial} and the fact that $\Par{\BorelStoch}$ is the partialization of $\BorelStoch$.

\begin{rem}
	In categorical probability, \newterm{conditionals}~\cite{fritz2019synthetic,dilavore2024partial} also play a central role.
	However, in this paper, they are only of tangential relevance, and we will merely comment on them for readers familiar with the concept. Nonetheless, $\Par{\BorelStoch}$ does have conditionals \cite[Proposition 4.10]{shahmohammed2025partial}.
\end{rem}

\subsection[Poset Enrichment]{Poset Enrichment of Quasi-Markov Categories}\label{sec:poset_enrichment}

Let us first recall some terminology introduced in earlier works such as \cite[Section B]{lorenz2023causalmodels}, \cite[Definition 2.13]{fritz2023lax}, and \cite[Definition 2.5]{gonda2024framework}.

\begin{defi}\label{def:domain}
	Let $f : A \to X$ be a morphism in a quasi-Markov category.
	Then the \newterm{domain} of $f$, denoted also by $\dom{f}$, is the morphism
	\begin{equation}\label{eq:domain}
		\input{domain.tikz}
	\end{equation}
\end{defi}

By quasi-totality, it is easy to see that the domain of any morphism is copyable and an idempotent.
For a total morphism, the domain is simply the identity.
This matches with its semantics in $\Par{\BorelStoch}$, where the domain of a morphism $f : A \to X$ is the inclusion $D_f \hookrightarrow A$, considered as the partial Markov kernel $A\to A$ that acts as the identity on its domain $D_f$ and is undefined elsewhere.
It is therefore only a slight abuse of terminology to use ``domain'' for both $D_f$ and for \Cref{eq:domain}.

In terms of domains, quasi-totality (\Cref{def:quasi-total}) amounts to every morphism absorbing its domain in the sense that $f \comp \dom{f} = f$ holds.
Consequently, the domain of a quasi-total monomorphism is simply the identity, which matches the totality of \Cref{lem:mono_total}.

Moreover, for arbitrary composable morphisms $f$ and $g$, we have
\begin{equation}\label{eq:dom_repeat}
	\input{dom_repeat.tikz}
\end{equation}
which can also be written as $\dom{f \comp g} = \dom{\dom{f} \comp g}$.
If $g$ is moreover copyable, we also have
\begin{equation}\label{eq:dom_comp_copyable}
	\dom{f} \comp g = g \comp \dom{f \comp g}.
\end{equation}

\begin{defi}\label{def:dom_ext}
	For any two parallel morphisms $f, g : A \to X$ in a quasi-Markov category, we say that $f$ \newterm{extends} $g$, denoted $f \domext g$, if we have
	\begin{equation}\label{eq:eq_on_domain}
		\input{eq_on_domain.tikz}
	\end{equation}
\end{defi}
For example, $\discard_X : X \to I$ extends any morphism of type $X \to I$.
The extension relation $\domext$ is a partial order on each hom-set of a quasi-Markov category \cite[Lemma 98]{lorenz2023causalmodels} (see also \cite{dilavore2025OrderPartialMarkov}).
Moreover, assuming positivity, we can show that it gives rise to an \emph{enrichment} in posets.
\begin{prop}[Enrichment of quasi-Markov categories]\label{prop:enrichment}
	Let $\cC$ be a positive quasi-Markov category.
	Then the extension partial order $\domext$ makes $\cC$ a category monoidally enriched in posets, i.e.~composition and tensor are monotone in each argument.
\end{prop}
\begin{proof}
	The fact that $\domext$ is preserved under tensoring is an immediate consequence of the fact that the copy is compatible with the tensor product.
	We therefore focus on the composition: 
	Consider parallel morphisms $f$ and $g$ such that $f \domext g$ holds.
	We first show compatibility with pre-composition by a morphism $h$, i.e.\ we show ${f \comp h \domext g \comp h}$.
	Since $\discard{} \comp g \comp h$ is copyable by quasi-totality of $g \comp h$, positivity gives
	\begin{equation}
		\input{quasi-MarkovEnrichmentPreCompLem.tikz}
	\end{equation}
	and thus we find
	\begin{equation}
		\input{quasi-MarkovEnrichmentPreComp.tikz}
	\end{equation}
	as desired.
	To show compatibility with post-composition by a morphism $k$, we need to prove $k \comp f \domext k \comp g$.
	Using the quasi-totality of $g$ in combination with $f \domext g$, we get\footnote{This shows that compatibility with post-composition holds also without positivity.}
	\begin{equation}
		\input{quasi-MarkovEnrichmentPostComp.tikz}
	\end{equation}
	and therefore
	\begin{equation}
		\input{quasi-MarkovEnrichmentPostCompConc.tikz}
	\end{equation}
	holds as desired.
\end{proof}

\begin{exa}[Extension relation for partial Markov kernels]\label{ex:domext_PBS}
	In $\Par{\BorelStoch}$, the extension ordering is exactly what one would expect: $f \domext g$ holds if and only if we have both $D_f \supseteq D_g$ and that $f$ and $g$ agree when restricted to $D_g$.
	Specifically, we have $\dom{f} \domext \dom{g}$ if and only if $D_f \supseteq D_g$ holds, so that the meet of $\dom{f}$ and $\dom{g}$ in $\Par{\BorelStoch}$ is the partial Markov kernel $(D_f \cap D_g, \id)$, which also equals $\dom{f} \comp \dom{g}$ in this category (see \Cref{ex:par_borelstoch}).
\end{exa}

This fact holds more generally, as we now show.
\begin{lem}[Meet of domains]\label{lem:meet_of_domains}
	Let $\cC$ be a quasi-Markov category, and $f : A \to X$, $g: A \to Y$ two morphisms in $\cC$.
	The meet of $\dom{f}$ and $\dom{g}$ is given by
	\begin{equation}\label{eq:meet_of_domains}
		\input{meet_of_domains.tikz}
	\end{equation}
\end{lem}
\begin{proof}
	The fact that this morphism is a lower bound follows from
	\begin{equation}
		\input{meet_lb.tikz}
	\end{equation}
	which establishes that $\dom{f}$ extends it, and an analogous argument shows that so does $\dom{g}$.
	
	Now consider another lower bound $h : A \to A$, which also extends the morphism from \eqref{eq:meet_of_domains}.
	That is, we have $\dom{f} \domext h$, $\dom{g} \domext h$, and $h \domext \dom{f} \comp \dom{g}$, which can be expanded to:
	\begin{gather}
		\label{eq:meet_lb_2}\input{meet_lb_2.tikz}\\
		\label{eq:meet_lb_4}\input{meet_lb_4.tikz}
	\end{gather}
	Combining the two equations from \eqref{eq:meet_lb_2} gives the first two equations in
	\begin{equation}
		\input{meet_lb_5.tikz}
	\end{equation}
	while the last one follows directly from \eqref{eq:meet_lb_4} upon discarding the output.
	In conclusion, the morphism from \eqref{eq:meet_of_domains} is indeed the greatest common lower bound of $\dom{f}$ and $\dom{g}$.
\end{proof}

One of the basic properties of the extension partial order that we need in the following section is that copyability is preserved by restricting the domain of a morphism.
\begin{lem}\label{lem:copyable_dom_rest}
	Consider two morphisms in a quasi-Markov category that satisfy $f \domext g$.
	\begin{enumerate}
		\item\label{it:copyable_dom_rest} If $f$ is copyable, then so is $g$.
		\item\label{it:dom_monotone} We have $\dom{f} \domext \dom{g}$.
		\item\label{it:idempotent_dom_rest} If $f$ is an idempotent, then $f \comp g = g$ holds.
		\item\label{it:domain_dom_rest} If $f$ is an identity, then $\dom{g} = g$ holds.
	\end{enumerate}
\end{lem}
\begin{proof}
	The condition $f \domext g$ means that $g = f \comp \dom{g}$.
	Quasi-totality ensures that $\dom{g}$ is copyable, so the composite $g$ is copyable provided that $f$ is.
	Explicitly,
	\begin{equation}
		\input{copyable_dom_rest.tikz}
	\end{equation}
	which is precisely saying that $g$ is copyable.

	\noindent %FIXED indent
	For the second part, we need to prove $\dom{f} \comp \dom{\dom{g}} = \dom{g}$, which amounts to
	\begin{equation}\label{eq:dom_monotone}
		\input{dom_monotone.tikz}
	\end{equation}
	where we use the definition of $f \domext g$ in the second equation.
	
	For the third part, we conclude directly from the definition of the extension partial order:
	\begin{equation}
		\input{idempotent_dom_rest.tikz}
	\end{equation}
	\noindent %FIXED indent
	Part \ref{it:domain_dom_rest} follows directly from \Cref{eq:eq_on_domain} when we replace $f$ by the identity.
\end{proof}

\subsection[Kolmogorov Products]{Kolmogorov Products in Quasi-Markov Categories}\label{sec:products}

\subsubsection{The Basic Definition}

We now turn to infinite tensor products in quasi-Markov categories.
These are particularly relevant to us, since we want to study empirical measures of \emph{infinite} sequences of samples.
The definition of Kolmogorov products has been introduced for Markov categories in~\cite{fritzrischel2019zeroone} and generalized to CD categories in~\cite{moss2022probability}.
However, in the presence of non-total morphisms we will need to exercise care with respect to domains of definition while performing certain constructions that are standard in Markov categories.
While a lax notion of Kolmogorov product suitable for quasi-Markov categories is introduced in~\cite[Section 5]{shahmohammed2025partial}, for the purposes of this article we will work with Kolmogorov products and domains of definition directly.
Here we restrict ourselves to countable families, but it is worth noting that the definition (as well as the results below) generalize straightforwardly to arbitrary families of objects.

\begin{defi}[{cf.~\cite[Definitions 3.1 and 4.1]{fritzrischel2019zeroone}}]\label{def:kolmogorov}
	Let $(X_i)_{i \in \N}$ be a sequence of objects in a quasi-Markov category $\cC$.
	Then a limit cone in $\cC$ of the diagram
	\begin{equation}
		\label{eq:kolmogorov_diagram}
		\begin{tikzcd}
			\dots \ar{r} & X_1 \otimes \dots \otimes X_n \ar{r} & \dots \ar{r} & X_1 \otimes X_2 \ar{r} & X_1,
		\end{tikzcd}
	\end{equation}
	in which the arrows are given by deletion of the last tensor factor, is a \newterm{Kolmogorov product} $\bigotimes_{i \in \N} X_i$ if the following conditions hold:
	\begin{enumerate}
		\item The cone components $\pi_n : \bigotimes_{i \in \N} X_i \to \bigotimes_{i \le n} X_i$ are deterministic;
		\item The limit cone is preserved by $\ph \otimes \id_Y$ for every $Y \in \cC$.

	\end{enumerate}
\end{defi}
\noindent %FIXED indent
In $\Par{\BorelStoch}$, countable Kolmogorov products exist and are (as in $\BorelStoch$) given by products of measurable spaces with their product $\sigma$-algebras~\cite[Proposition 5.5]{shahmohammed2025partial}.

\begin{nota}\label{not:products_of_objects}
	In order to simplify the exposition, for positive integers $m < n$ we write
	\begin{equation}\label{eq:products_of_objects}
		X_m^n \coloneq \bigotimes_{i = m}^n X_i  \qquad \qquad  X_m^\N \coloneq \bigotimes_{i \ge m} X_i,
	\end{equation}
	where the latter is a countable Kolmogorov product as defined above.
	Similarly, for a family of morphisms $(f_i : X_i \to Y_i)_{i=m}^n$, we write
	\begin{equation}\label{eq:products_of_morphisms}
		f_m^n \coloneq \bigotimes_{i = m}^n f_i 
	\end{equation}
	for their parallel composite, which is of type $X_m^n \to Y_m^n$.
\end{nota}

\begin{lem}\label{lem:KP_dom_cop}
	Consider a family ${\left(g_n : A \to X_1^n \right)}_{n \in \N}$ that forms a cone over Diagram (\ref{eq:kolmogorov_diagram}), meaning that its elements satisfy $\discard_{X_{n+1}} \comp g_{n+1} = g_n$ for all $n \in \N$.
	Then:
	\begin{enumerate}
		\item\label{item:Kolmogorov_cone_domains} The morphisms $g_n$ in the family have identical domains.
		\item\label{item:Kolmogorov_map_domain} The domain of $g : A \to X_1^\N$, the morphism induced by the universal property, is equal to the common domain of morphisms in the family.
		\item\label{item:Kolmogorov_map_copyablity} The morphism $g : A \to X_1^\N$ is copyable if and only if each $g_n$ is copyable.
	\end{enumerate}
\end{lem}
\begin{proof} %FIXED Split for Margins
	Claim~\ref{item:Kolmogorov_cone_domains} can be proven by induction.
	Since every marginalization map ${\discard_{X_{n+1}} : X_1^{n+1}}$
	$ \to X_1^{n}$ is total, we have
	\begin{equation}
		\dom{g_{n}} = \dom{\discard_{X_{n+1}} \comp g_{n+1}} = \dom{g_{n+1}}
	\end{equation}
	for every $n \in \N$, which proves the claim.

	Claim~\ref{item:Kolmogorov_map_domain} follows from~\ref{item:Kolmogorov_cone_domains} and the assumption that the Kolmogorov cone projections are total, since we have $g_n = \pi_n \comp g$ by the universal property.

	Claim~\ref{item:Kolmogorov_map_copyablity} is a consequence of the assumption that the Kolmogorov cone legs $\pi_n$ are copyable by \Cref{def:kolmogorov}.
	In particular, if $g : A \to X_1^\N$ is copyable, then each $g_n = \pi_n \comp g$ is copyable as a composite of copyable morphisms.
	
	Conversely, if each $g_n$ is copyable, then we have, for any $n \in \N$,
	\begin{equation}\label{eq:copyable_n}
		\input{copyable_n.tikz}
	\end{equation}
	Thanks to the fact that Kolmogorov products are preserved by tensor products, applying the universal property to the first output yields
	\begin{equation}\label{eq:copyable_n_2}
		\input{copyable_n_2.tikz}
	\end{equation}
	and doing the same for the second output shows that $g$ itself is copyable.
\end{proof}

\subsubsection{Infinite Parallel Composites and \texorpdfstring{$\sigma$}{σ}-continuity}

Unfortunately, this direct instantiation of the notion of Kolmogorov products from Markov categories is insufficient to perform some of the constructions we need in quasi-Markov categories.
In particular, we would often like to have a version of the universal property that can be applied to cones that do not have identical domains, which by \Cref{lem:KP_dom_cop}~\ref{item:Kolmogorov_cone_domains} is not possible with plain Kolmogorov products.

For instance, in a Markov category one can construct an infinite parallel composite of morphisms ${f_i : X_i \to Y_i}$ as the unique morphism $X_1^{\N} \to Y_1^{\N}$ induced by the cone given by the composites ${f_1^n \comp \pi_n : X_1^{\N} \to Y_1^{n}}$.
However, in a quasi-Markov category, these maps need \emph{not} have the same domain; and hence they cannot form a cone over Diagram (\ref{eq:kolmogorov_diagram}).
Intuitively, the morphism $f_1^n \comp \pi_n$ is only defined on those sequences $(x_i)_{i \in \N}$ for which each the first $n$ elements is in the domain of the corresponding $f_i$.
This is typically a smaller space as $n$ increases, and hence the domains of the $f^n \comp \pi_n$ need not agree.

For this reason, when using the Kolmogorov product to construct a morphism $X_1^{\N} \to Y_1^{\N}$ from a family of morphisms $f_i : X_i \to Y_i$, the induced morphism should only be defined on those input sequences ${\left(x_i\right)}_{i \in \N}$ for which \emph{all} $f_i(x_i)$ are defined.
This domain is the meet of the domains in the family $\left(f^n\comp \pi_n\right)_{n \in \N}$.
Moreover, this meet is additionally directed, because the sequence of domains is decreasing.
We will see that this is semantically related to $\sigma$-additivity of probability measures, since in $\Par{\BorelStoch}$, countable meets of domains correspond to countable intersections of the associated measurable sets.

To ensure that Kolmogorov products behave in a manner reflecting their intended behavior in probability theory, we need the following more explicit handle on intersections of domains of morphisms.

\begin{defi}\label{def:count_meets}
	A positive quasi-Markov category $\cC$ with poset enrichment given by the extension relation $\domext$ is termed $\bm{\sigma}$\textbf{-continuous} if:
\begin{enumerate}
		\item\label{it:count_meets_exist} The hom-sets of $\cC$ have countable directed meets in the sense that for any descending sequence of morphisms $\left(f_n : X \to Y \right)_{n \in \N}$, the meet $\bigwedge_{j \in \N} f_{j}$ exists;
		\item\label{it:count_meets_comp} This meet is preserved by sequential and parallel composition. 
	\end{enumerate}
\end{defi}
\noindent %FIXED indent
Property \ref{it:count_meets_comp} means that we get an enrichment in the category of posets with countable directed meets and monotone maps which preserve those (with the Cartesian product as the monoidal structure).

Before we turn to proving that $\sigma$-continuity holds in $\Par{\BorelStoch}$ in \Cref{prop:sigma_parborelstoch}, let us illustrate its consequences.
The common domain of two morphisms $f, g : X \to Y$, i.e.\ the meet of their domains, is given by the composite $\dom{f} \comp \dom{g}$ as we saw in \Cref{lem:meet_of_domains}.
Moreover, this lemma also implies that the meet of domains is preserved both by sequential and parallel composition.
\Cref{def:count_meets} extends this property to arbitrary morphisms and to countable directed meets. 
Based on this, $\sigma$-continuity allows us to construct infinite parallel composites of morphisms as follows.

\begin{lem}[Infinite parallel composites of morphisms]\label{lem:inf_tensor}
	Let $\cC$ be a positive and $\sigma$-continuous quasi-Markov category with countable Kolmogorov products.
	Consider a countable sequence of morphisms $f_i : X_i \to Y_i$ in $\cC$.
	Then the composite morphisms (one for each $n \in \N$)
	\begin{equation}\label{eq:finite_tensor_morphisms}
		X_1^{\N} \xrightarrow{\bigwedge_{m \in \N} \dom{f_1^m \pi_m}} X_1^{\N} \xrightarrow{\pi_n} X_1^{n} \xrightarrow{f_1^n} Y_1^{n},
	\end{equation}
	expressed using \Cref{not:products_of_objects}, form a cone over the Kolmogorov diagram from (\ref{eq:kolmogorov_diagram}).
\end{lem}
\begin{proof}
	We start by noting that $\discard_{Y_{n+1}} \domext \discard_{Y_{n+1}} \comp f_{n+1}$ holds by the definition of the extension order, which implies 
	\begin{equation}\label{eq:inf_tensor_1}
		\input{inf_tensor_1.tikz}
	\end{equation}
	by the fact that $\domext$ is a poset enrichment (\Cref{prop:enrichment}).
	The left-hand side of \Cref{eq:inf_tensor_1} is equal to $f_1^n \comp \pi_n$ because the $\pi_n$ are the limit cone components and hence satisfy ${\discard_{X_{n+1}} \comp \pi_{n+1} = \pi_n}$.
	This leads to an inclusion of domains,
	\begin{equation}\label{eq:inf_tensor_3}
		\dom{f_1^n \pi_n} \domext \dom{f_1^{n+1} \pi_{n+1}}
	\end{equation}
	by \Cref{lem:copyable_dom_rest}~\ref{it:dom_monotone} and thus also 
	\begin{equation}\label{eq:inf_tensor_2}
		f_1^n \comp \pi_n \comp \dom{f_1^{n+1} \pi_{n+1}} = \discard_{Y_{n+1}} \comp f_1^{n+1} \comp \pi_{n+1}. \qquad
	\end{equation}
	
	\Cref{eq:inf_tensor_3} means that the morphisms $\dom{f_1^m \pi_m}$ form a descending sequence and thus their countable meet exists by $\sigma$-continuity.
	Let us simplify the notation and denote it by
	\begin{equation}\label{eq:inf_tensor_meet}
		\varphi \coloneq \bigwedge_{m \in \N} \dom{f_1^m \comp \pi_m}.
	\end{equation}

	For any $k \in \N$, we then have 
	\begin{equation}\label{eq:inf_tensor_4}
		\begin{split}
			\dom{f_1^k \pi_k} \comp \varphi &=  \bigwedge_{m \in \N} \left[ \dom{f_1^k \comp \pi_k} \comp \dom{f_1^m \comp \pi_m} \right] \\
				&= \bigwedge_{m \ge k} \dom{f_1^m \comp \pi_m}  \\
				&= \varphi ,
		\end{split}
	\end{equation}
	where the first equation is by \Cref{def:count_meets}~\ref{it:count_meets_comp}, the second by \Cref{lem:meet_of_domains,eq:inf_tensor_3}, and the last one by repeated application of \Cref{eq:inf_tensor_3} again.
	
	We can use these ingredients to establish the claim by combining \Cref{eq:inf_tensor_2,eq:inf_tensor_4} to get
	\begin{equation}
		\begin{split}
			\discard_{Y_{n+1}} \comp f_1^{n+1} \comp \pi_{n+1} \comp \varphi &= f_1^n \comp \pi_n \comp \dom{f_1^{n+1} \comp \pi_{n+1}} \comp \varphi \\[2pt]
				&= f_1^n \comp \pi_n \comp \varphi,
		\end{split}
	\end{equation}
	which is the equation required to establish that the morphisms from \eqref{eq:finite_tensor_morphisms} form a cone over the Kolmogorov diagram.
\end{proof}

Consequently, we can define the \newterm{infinite parallel composite} $f_1^\N : X_1^{\N} \to Y_1^{\N}$ of the family $(f_i)_{i \in \N}$ as the unique morphism induced by the cone from \Cref{lem:inf_tensor} via the universal property of the Kolmogorov product $Y_1^\N$.
As we show in \Cref{lem:inf_tensor_domain} below, its domain is given by the meet $\varphi$ from \eqref{eq:inf_tensor_meet}.

\begin{rem}
	\Cref{lem:inf_tensor} demonstrates that, under the assumption of $\sigma$-continuity, Kolmogorov products are not only limits over the diagram~\eqref{eq:kolmogorov_diagram} in the usual sense, but also restriction limits in the sense of Cockett and Lack \cite[Section 4.4]{cockettlackRestrictionCategoriesIII2007}.
	In their terminology, a lax cone over the diagram~\eqref{eq:kolmogorov_diagram} is a family of morphisms $(g_n : A \to X_1^n)_{n \in \N}$ satisfying the lax commutativity conditions $g_n \domext \discard_{X_{n+1}} \comp g_{n+1}$ for all $n \in \N$.
	
	The proof of \Cref{lem:inf_tensor}, in particular \Cref{eq:inf_tensor_2}, demonstrates that the composite morphisms
	\[
		X^{\mathbb N} \xrightarrow{\pi_n} X^n \xrightarrow{f_1^n} Y^n
	\]
	form such a lax cone.
	Moreover, by restricting each of them to their common domain $\varphi$, we obtain a strict cone in the sense of \Cref{def:kolmogorov}.
	
	One could generalize \Cref{lem:inf_tensor} to show that given any lax cone $(A \xrightarrow{g_n} X_1^n)_{n \in \N}$, there is a unique morphism $A \xrightarrow{g} X_1^n$ such that for all $n$, we have $g_n \domext \pi_n \comp g$ and moreover we also have $\dom{g} = \bigwedge_{m \in \N} \dom{g_{m}}$.
	This is precisely the universal property of a restriction limit as formulated in~\cite[Section 2.1]{cockett2012rangeCategoriesII}.
	More on this perspective can be found in~\cite[Section 5]{shahmohammed2025partial}.
	For the purposes of the current work, where only infinite parallel composites arise, it suffices to work with the products and meets explicitly as in \Cref{lem:inf_tensor}.
\end{rem}

\begin{lem}[Finite projections of infinite parallel composites]\label{lem:inf_tensor_proj}
	Given the same assumptions as in \Cref{lem:inf_tensor}, we have
	\begin{equation}\label{eq:inf_tensor_proj}
		\pi_n \comp f_1^\N = \bigwedge_{m \ge n} \left( \discard_{Y_{n+1}^m} \comp f_1^m \comp \pi_m \right),
	\end{equation}
	i.e.\ the projection of the infinite tensor product $f_1^\N$ onto $Y_1^{n}$ is given by the meet of the projections of the finite tensor products $f_1^m$.
\end{lem}
\begin{proof}
	By \Cref{lem:inf_tensor,eq:inf_tensor_3}, we have the first equation in
	\begin{equation}\label{eq:inf_tensor_proj_intermediate}
		\begin{split}
			\pi_n \comp f_1^{\N} &= f_1^n \comp \pi_n \comp \bigwedge_{m \geq n} \dom{f_1^m \pi_m} \\
			&= \bigwedge_{m \ge n} \Bigl(f_1^n \comp \pi_n \comp \dom{f^m \pi_m}\Bigr),
		\end{split}
	\end{equation}
	while the second one is an application of the assumption that countable meets are preserved by sequential composition (\Cref{def:count_meets}).
	For any $m>n$, we can apply \Cref{eq:inf_tensor_2} iteratively to obtain 
	\begin{equation}
		f_1^n \comp \pi_n \comp \dom{f^m \pi_m} = \discard_{Y_{n+1}^m} \comp f_1^m \comp \pi_m,
	\end{equation}
	which also holds for $m = n$ by quasi-totality.
	Plugging this into~\Cref{eq:inf_tensor_proj_intermediate} yields the desired result.
\end{proof}

In the next result, we establish that the domain of the parallel composite of a family $(f_i)_{i \in \N}$ is the infinite parallel composite of the individual domains.
Consequently, it becomes clear that $\dom{f_1^\N}$ is an infinite parallel composite of copyable morphisms, and hence can be studied also from the perspective of restriction products and cartesian bicategories~\cite{cockettlackRestrictionCategoriesIII2007,carboniwalters1987bicartesian}.

\begin{lem}[Domain of infinite parallel composites]\label{lem:inf_tensor_domain}
	Given the same assumptions as in \Cref{lem:inf_tensor}, we have
	\begin{equation}\label{eq:inf_tensor_domain}
		\dom{f_1^\N} = \bigwedge_{m \in \N} \dom{f_1^m \comp \pi_m} = \bigotimes_{i \in \N} \dom{f_i}.
	\end{equation}
\end{lem}
\begin{proof}
	Let us once again use the notation from \Cref{eq:inf_tensor_meet}.
	In order to compute the domain of the infinite parallel composite, according to \Cref{lem:KP_dom_cop}~\ref{item:Kolmogorov_map_domain} we can equivalently compute the domain of either one of the legs of its cone.
	To this end, we find
	\begin{equation}
		\begin{split}
			\dom{f_1^n \comp \pi_n \comp \varphi} &= \dom{\dom{f_1^n \comp \pi_n} \comp \varphi} \\[2pt]
			&= \dom{\varphi}
		\end{split}
	\end{equation}
	where the first equation is by \eqref{eq:dom_repeat} and the second one by \eqref{eq:inf_tensor_4}.
	Since countable meets are preserved by sequential and parallel composition (\Cref{def:count_meets}), we can express $\dom{\varphi}$ as 
	\begin{equation}\label{eq:count_meet_domain}
		\input{count_meet_domain_1.tikz} \;\;\; = \;\, \bigwedge_{m \in \N} \left( \;\input{count_meet_domain_2.tikz} \, \right) \;\; = \;\, \bigwedge_{m \in \N} \varphi_m
	\end{equation}
	where we denote $\dom{f_1^m \pi_m}$ by $\varphi_m$, which satisfies $\dom{\varphi_m} = \varphi_m$.
	Since the right-hand side of \eqref{eq:count_meet_domain} equals $\varphi$, this establishes $\dom{\varphi} = \varphi$ and hence the first equality in \eqref{eq:inf_tensor_domain}.
	
	To prove the second equality, it suffices to consider its composition with $\pi_n$ for each $n \in \N$.
	Using \Cref{lem:inf_tensor_proj}, we can compute $\pi_n \comp \bigotimes_{i \in \N} \dom{f_i}$ as the meet of morphisms
	\begin{equation}
		\input{dom_prod_marginal_3.tikz}
	\end{equation}
	over $m \geq n$, where we have used the copyability of $\pi_m$ as well as $f_1^m = f_1^{n} \otimes f_{n+1}^m$ to establish the equation.
	Using $\pi_n = \discard_{X_{n+1}^m} \comp \pi_m$ and $\sigma$-continuity, we thus obtain the first equation in
	\begin{equation}
		\begin{split}
			\pi_n \comp \bigotimes_{i \in \N} \dom{f_i} &= \pi_n \comp \bigwedge_{m \geq n} \dom{f_1^m \comp \pi_m} \\
				&= \pi_n \comp \bigwedge_{m \in \N} \dom{f_1^m \comp \pi_m},
		\end{split}
	\end{equation}
	where the second equation is by \eqref{eq:inf_tensor_3}.
	By the universal property of Kolmogorov products, we have therefore shown the second equality in \eqref{eq:inf_tensor_domain}.
\end{proof}

One can also state a version of \Cref{lem:KP_dom_cop} for the infinite parallel composite.
\begin{lem}\label{lem:inf_tensor_cop}
	Given the same assumptions as in \Cref{lem:inf_tensor}:
	\begin{enumerate}
		\item If each $f_i$ is copyable, then so is the infinite parallel composite $f_1^\N$.
		\item If each $f_i$ is total, then so is $f_1^\N$.
	\end{enumerate}
\end{lem}
\begin{proof}
	The first part follows from \Cref{lem:KP_dom_cop}~\ref{item:Kolmogorov_map_copyablity}.
	Indeed, each leg $f^n \comp \pi_n \comp \dom{f_1^\N}$ of the defining cone is a composite of copyable morphisms, where we make use of \Cref{lem:inf_tensor_domain}.
	
	The second follows from \Cref{lem:KP_dom_cop}~\ref{item:Kolmogorov_map_domain}.
	If each $f_i$ is total, then so is each $f^n \comp \pi_n$.
	In this case the meet $\bigwedge_{m \in \N} \dom{f^m \comp \pi_m}$ is the identity, so that their composite is total as well.
\end{proof}
\begin{rem}
	The converses of the two assertions in \Cref{lem:inf_tensor_cop} hold provided that all the morphisms $f_i$ are equal to each other.
	While the case of totality is clear from definition, the converse for copyability can be shown as a consequence of the upcoming \Cref{prop:IID_inf_copy} and quasi-totality.
\end{rem}

One of the most important properties that is ensured by $\sigma$-continuity is the functoriality of forming infinite parallel composites.

\begin{prop}[Functoriality of infinite parallel composites]\label{prop:inf_tensor_composition}
	Let $\cC$ be a positive and $\sigma$-continuous quasi-Markov category with countable Kolmogorov products.
	Then the infinite parallel composite construction from \Cref{lem:inf_tensor} defines a functor $\cC^{\N} \to \cC$.
\end{prop}
\begin{proof}
	The preservation of identities is straightforward to verify.
	Indeed, if $f_1^n$ is given by $\id_{X^n}$, then the morphism from \eqref{eq:finite_tensor_morphisms} equals $\pi_n$, i.e.\ we get the limit cone itself.
	
	For the preservation of composites, consider sequences of morphisms $f_i : X_i \to Y_i$ and $g_i : Y_i \to Z_i$ for $i \in \N$ in $\cC$, and define $h_i \coloneq g_i \comp f_i$.
	The task is to show that
	\begin{equation}
		h_1^\N = g_1^\N \comp f_1^\N.
	\end{equation}
	To obtain this, it suffices show that the projections of both sides onto $Z_1^n$ coincide.
	According to \Cref{lem:inf_tensor_proj}, applying $\pi_n$ gives
	\begin{equation}\label{eq:projection_1}
		\pi_n \comp h_1^\N = \bigwedge_{m \ge n} \left( \discard_{Z_{n+1}^m} \comp g_1^m \comp f_1^m \comp \pi_m \right)
	\end{equation}
	and
	\begin{equation}\label{eq:projection_2}
		\pi_n \comp g_1^\N \comp f_1^\N = \bigwedge_{m \ge n} \left( \discard_{Z_{n+1}^m} \comp g_1^m  \comp \pi_m  \comp f_1^\N \right),
	\end{equation}
	respectively, where we have used $\sigma$-continuity to pull $f_1^\N$ inside the meet in the latter case.
	Using the definition of the infinite parallel composite and \Cref{lem:inf_tensor_domain}, we also have
	\begin{equation}\label{eq:projection_3}
		\pi_m \comp f_1^\N = f_1^m \comp \pi_m \comp \dom{f_1^\N}.
	\end{equation}
	Combining \Cref{eq:projection_1,eq:projection_2,eq:projection_3}, we get
	\begin{equation}\label{eq:projection_4}
		\pi_n \comp g_1^\N \comp f_1^\N = \pi_n \comp h_1^\N \comp \dom{f_1^\N} =  \pi_n \comp h_1^\N \comp \dom{h_1^\N} \comp \dom{f_1^\N}
	\end{equation}
	where in the second step we use quasi-totality to replace $h_1^\N$ by $h_1^\N \comp \dom{h_1^\N}$.
	By \Cref{lem:meet_of_domains}, the composite $\dom{h_1^\N} \comp \dom{f_1^\N}$ equals the meet of these two morphisms, and so once we show $\dom{f_1^\N} \domext \dom{h_1^\N}$, we can drop $\dom{f_1^\N}$ from the right-hand side of \Cref{eq:projection_4} and obtain the desired result.
	
	To prove $\dom{f_1^\N} \domext \dom{h_1^\N}$, note that $\discard_{Y_{1}^m} \domext \discard_{Z_{1}^m} \comp g_1^m$ trivially holds and implies 
	\begin{equation}
		\dom{f_1^m \comp \pi_m} \domext \dom{g_1^m \comp f_1^m \comp \pi_m}
	\end{equation}
	by \Cref{prop:enrichment}.
	Thus, we also get 
	\begin{equation}
		\dom{f_1^\N} = \bigwedge_{m \in \N} \dom{f_1^m \comp \pi_m}  \domext \bigwedge_{m \in \N} \dom{g_1^m \comp f_1^m \comp \pi_m} = \dom{h_1^\N}
	\end{equation}
	by \Cref{lem:inf_tensor_domain} and the proof is complete.
\end{proof}

In order to make use of these results on $\sigma$-continuity for quasi-Markov categories in our path towards synthetic strong laws of large numbers, we must establish that this property holds in measure-theoretic probability.

\begin{prop}\label{prop:sigma_parborelstoch}
	$\Par{\BorelStoch}$ is $\sigma$-continuous.
\end{prop}
\begin{proof}
	We can verify this explicitly via monotone $\sigma$-continuity of measures.
	We consider a generic countable descending chain of morphisms ${\left(u_n : A \to X\right)}_{n \in \N}$ in $\Par{\BorelStoch}$
	with representatives
	\begin{equation}
		A \supseteq D_n \xrightarrow{f_n} X,
	\end{equation}
	where each $f_n$ is a Markov kernel.
	First, let us prove the existence of their meet.
	By $u_n \domext u_{n+1}$ and the characterization of the extension relation in measure-theoretic probability (\Cref{ex:domext_PBS}), we have
	\begin{equation}\label{eq:descending_kernels}
		D_n \supseteq D_{n+1} \qquad \text{and} \qquad f_{n+1} \left( E \,\vert\, a\right) = f_{n}\left(E \,\vert\, a\right)
	\end{equation}
	for all $n\in\N$, all $a \in D_{n+1}$, and all measurable $E \subseteq X$.
	The intersection of all the domains $D_f \coloneqq \bigcap_{n \in \N} D_n$ is measurable as a countable intersection of measurable sets. 
	Moreover, we can define a Markov kernel $f : D_f \to Y$ via $f\left(E \,\vert\, a\right) \coloneqq f_n\left(E \,\vert\, a\right)$ by picking any $n$, which is independent of $n$ since the $f_n$ agree on $D_f$ by \eqref{eq:descending_kernels}.
	This gives us a
	partial Markov kernel $u$ represented by 
	\begin{equation}
		A \supseteq D_f \xrightarrow{f} X.
	\end{equation}
	We argue that $u$ is the meet of the chain we started with.
	To this end, consider any partial Markov kernel $w$, represented by
	\[
		A \supseteq D_h \xrightarrow{h} X,
	\]
	which satisfies $u_n \domext w$ for all $n$.
	By \Cref{ex:domext_PBS}, we have $D_n \supseteq D_h$ for all $n$, which implies $D_f \supseteq D_h$.
	Furthermore, for any $a \in D_g$, we also have ${h \left(E \,\vert\, a\right) = f_n\left(E \,\vert\, a\right)} $, thus showing the desired $u \domext w$.
	This establishes the existence of countable directed meets in any hom-set of $\Par{\BorelStoch}$.

	To show that the meet $u = \bigwedge_n u_n$ is preserved by tensoring, consider an arbitrary partial Markov kernel $v$ represented by
	\[
		B \supseteq D_g \xrightarrow{g} Y.
	\]
	Then we have a descending chain with elements $u_n \otimes v$ given by
	\begin{equation}
		A \times B \supseteq D_n \times D_g \xrightarrow{f_n \otimes g} X \times Y.
	\end{equation}
	The meet $\bigwedge_{n \in \N} \left(u_n \otimes v \right)$ has domain
	\begin{equation}
		\bigcap_n \, \left(D_n \times D_g\right) = \biggl( \, \bigcap_n D_n \biggr) \times D_g = D_f \times D_g.
	\end{equation}
	On cylinder sets $E \times F \subseteq X \times Y$, the meet gives the transition probability
	\begin{equation}
		\bigwedge_n \, \left(u_n \otimes v\right) \bigl( E \otimes F \, \vert \, (a,b) \bigr) = f_n\left(E \,\vert\, a\right) \, g\left(F \,\vert\, b\right) = f\left(E \,\vert\, a\right) \, g\left(F \,\vert\, b\right)
	\end{equation}
	for all $\left(a,b\right) \in D_f \times D_g$.
	Comparing on cylinder sets is enough to ensure that this coincides with $f \otimes g$, and so we conclude $\bigwedge_{n \in \N} \left(u_n \otimes v \right) = u \otimes v$ as needed.
	
	Now we check compatibility of the above meet with post-composition by an arbitrary partial Markov kernel $w$ represented by
	\[
		X \supseteq D_g \xrightarrow{g} Y.
	\]
	%FIXED Split into two lines for margins
	Recalling \Cref{eq:dom_comp}, composites of the form $w \comp u_n$ have domain $E_n = \left\{ a \in D_n \quad \middle| \right. \vphantom{\}} $
	$\left. \vphantom{\{}f_n \left(D_g \, \vert \, a \right)= 1 \right\}$,
	and they form a descending chain.
	Per the above construction of meets, $\bigwedge_n \left(w \comp u_n\right)$ has domain
	\begin{equation}
		\begin{split}
			\bigcap_n E_n &= \Set{ a \in \bigcap_n D_n \given \forall n \in \N \; : \; f_n\left(D_g \, \vert \, a\right)= 1 } \\[2pt]
				&= \Set*[\Big]{ a \in D_f \given f \left(D_g \, \vert \, a \right)= 1 },
		\end{split}
	\end{equation}
	where the second equality follows because $f_n\left(D_g \, \vert \, a\right) = 1$ holds if and only if $f\left(D_g \, \vert \, a\right) = 1$ does for any $a \in \bigcap_n D_n = D_f$.
	This set is also the domain of $w \comp u$.
	Furthermore, on this domain, $f_n$ and $f$ are concentrated pointwise on $D_g$, and act the same.
	Hence the composites act the same as well, and we have $\bigwedge_n \left(w \comp u_n\right) = w \comp u$.

	Finally, we check compatibility with pre-composition by an arbitrary partial Markov kernel $v$ represented by
	\[
		B \supseteq D_g \xrightarrow{g} A.
	\]
	Each composite $u_n \comp v$ has domain $F_n \coloneqq \Set{b \in D_g \given g\left(D_n \, \vert \, b\right) = 1 }$, and these again remain a descending chain.
	Therefore the meet $\bigwedge_n \left(u_n \comp v\right)$ has domain
	\begin{equation}\label{eq:dom1}
		\bigcap_n F_n = \Set*[\Big]{ b \in D_g \given \forall n \in \N \; : \; g\left(D_n \, \vert \, b\right) = 1 }.
	\end{equation}
	On the other hand, the domain of $u \comp v$ is 
	\begin{equation}\label{eq:dom2}
		\Set{b \in D_g \given g\biggl( \,\bigcap_n D_n \, \Big\vert \, b \biggr) = 1 }.
	\end{equation}
	Since the $D_n$ form a descending chain of measurable sets, and for any $b \in D_g$ the probability measure $g \left(\ph \, \vert \, b\right)$ is $\sigma$-continuous from above~\cite[Theorem 2.1 (ii)]{billingsley}, meaning
	\begin{equation}
		g \biggl( \, \bigcap_n D_n \, \Big\vert \, b \biggr) = \inf_{n \in \N} g\left(D_n \,\vert\, b \right),
	\end{equation}
	we conclude that $g \left( D_n \, \vert \, b \right) = 1$ holds for all $n \in \N$ if and only if $g\left(\, \bigcap_n D_n \, \vert \, b \right) = 1$ does.
	Therefore the two domains from \eqref{eq:dom1} and \eqref{eq:dom2} coincide.
	Moreover, on this domain $g$ is concentrated entirely on $D_f$ by construction.
	Since each $f_n$ agrees pointwise with $f$ on $D_f$, the composites $u_n \comp v$ and $u \comp v$ coincide on their domain, hence establishing $\bigwedge_n \left(u_n \comp v\right) = u \comp v$.
\end{proof}

\subsubsection{Kolmogorov Powers and IID Morphisms}\label{sec:IID}

Suppose now that all $f_i : X_i \to Y_i$ are identical, which is the case of interest for us in the rest of this article.
Then we also write $X^\N$ for the Kolmogorov product $\bigotimes_{i \in \N} X$ and refer to it as a \newterm{Kolmogorov power}.
Similarly, we denote the infinite parallel composite of countably many instances of $f$ by ${f^\N : X^\N \to Y^\N}$.
%FIXED pushed quarter a line for margins
\vskip0.25\baselineskip
\noindent %FIXED indent
Given any morphism $f : A \to X$, we can also form a compatible family $\left( f^{(n)} : A \to X^n \right)_{n \in \N}$, in the sense that it provides a cone over the Kolmogorov diagram \eqref{eq:kolmogorov_diagram}, by defining
\begin{equation}\label{eq:notation_fN}
	f^{(n)} \coloneqq f^n \comp \cop,
\end{equation}
i.e.\ $f^{(n)}$ is the composite of the $n$-fold copy morphism $\cop : A \to A^n$ followed by the $n$-fold tensor power of $f$.
Thanks to the quasi-totality of $f$, this defines a cone over the Kolmogorov diagram~(\ref{eq:kolmogorov_diagram}).
If the Kolmogorov product $X^\N$ exists, we obtain a morphism
\begin{equation}\label{eq:IID}
	f^{(\N)} : A \longrightarrow X^\N,
\end{equation}
which commonly appears in categorical probability.
It is the unique morphism $A \to X^\N$ whose finite marginals are the \eqref{eq:notation_fN}.
Intuitively, it corresponds to a probabilistic process that produces infinitely many samples, which are conditionally independent given $A$ and identically distributed.

For string-diagrammatic manipulations with such \newterm{IID morphisms}, we use the \newterm{plate notation} \cite[Section 2.5]{chen2024aldoushoover}:
\begin{equation}
	\label{eq:plate}
	\input{plate_notation.tikz}
\end{equation}
For $f= \id_A$ the above construction gives us the $\N$\textbf{-fold copy} morphism $\id^{(\N)} : A \to A^\N$.

Under the additional assumption of $\sigma$-continuity, we can prove that IID morphisms indeed display conditional independence given $A$ in the following sense.

\begin{prop}[IID morphisms are independent products]\label{prop:IID_inf_copy}
	Let $\cC$ be a positive and $\sigma$-continuous quasi-Markov category with countable Kolmogorov products.
	Then for any morphism $f : A \to X$ in $\cC$, we have
	\begin{equation}\label{eq:IID_inf_copy}
		f^{\left(\N\right)} = f^{\N}\comp \id^{(\N)}.
	\end{equation}
\end{prop}

\begin{proof}
	Let us first show a useful intermediate result.
	Namely, our claim is that $\dom{f}$ is ``countably-copyable'' in the following sense:
	\begin{equation}\label{eq:dom_inf_copyable}
		\dom{f}^\N \comp \id^{(\N)} = \id^{(\N)} \comp \dom{f}.
	\end{equation}
	To show this, we use \Cref{lem:inf_tensor_domain} and $\sigma$-continuity to deduce
	\begin{equation}
		\begin{split}
			\dom{f}^\N \comp \id^{(\N)} &= \bigwedge_{n \in \N} \Bigl[ \dom{f^n \comp \pi_n} \comp \id^{(\N)} \Bigr] \\
				&= \bigwedge_{n \in \N} \Bigl[ \id^{(\N)} \comp \dom{f^n \comp \pi_n \comp \id^{(\N)}} \Bigr] \\
				&= \id^{(\N)} \comp  \bigwedge_{n \in \N} \dom{f^n \comp \id^{(n)}} = \id^{(\N)} \comp \dom{f},
		\end{split}
	\end{equation}
	where the second equation is by \eqref{eq:dom_comp_copyable}, the third one is by the definition of the $\N$-fold copy and $\sigma$-continuity, and the last one by copyability of domains which implies $\dom{f^{(n)}} = \dom{f}$.
	
	Turning our attention to the claim \eqref{eq:IID_inf_copy} now, it suffices to show that the finite projections of both sides agree.
	For the left-hand side, we have $\pi_n \comp f^{\left(\N\right)} = f^n \comp \id^{(n)}$ by construction.
	On the other hand, we can compute
	\begin{equation}
		\begin{split}
			\pi_n \comp f^{\N}\comp \id^{(\N)} &= f^n \comp \pi_n \comp \dom{f}^\N \comp \id^{(\N)} \\
				&= f^n \comp \pi_n \comp \id^{(\N)} \comp \dom{f} \\
				&= f^n \comp \id^{(n)} \comp \dom{f} = f^n \comp \id^{(n)},
		\end{split}
	\end{equation}
	where the first equation is by the definition of the infinite parallel composite $f^\N$ and \Cref{lem:inf_tensor_domain}, the second one uses \Cref{eq:dom_inf_copyable}, and the last two are as in the previous paragraph.
	Therefore the marginals agree for all $n \in \N$, and we conclude $f^{\left(\N\right)} = f^{\N}\comp \id^{(\N)}$.
\end{proof}

\begin{prop}\label{lem:IID_naturality}
	Let $\cC$ be a quasi-Markov category with countable Kolmogorov products and consider two morphisms $f \colon X \to Y$ and $g \colon Y \to Z$.
	
	If $\cC$ is positive and $\sigma$-continuous, then we have
	\begin{equation}\label{eq:IID_naturality}
		\input{IID_naturality1.tikz}
	\end{equation}
	
	If $f$ is copyable, then we have
	\begin{equation}\label{eq:IID_naturality2}
		\input{IID_naturality2.tikz}
	\end{equation}
\end{prop}
\begin{proof}
	The first equation is a direct consequence of \Cref{prop:inf_tensor_composition,prop:IID_inf_copy}.
	For the second, we consider the finite marginals to find
	\begin{equation}
		\input{IID_naturality_finite.tikz}
	\end{equation}
	The first and third equalities use the definition of $g^{\left(\N\right)}$ and $(g \circ f)^{\left(\N\right)}$, while the second one follows from the copyability of $f$.
\end{proof}

\begin{rem}[Copyable implies countably-copyable]
	In particular, we can now generalize \Cref{eq:dom_inf_copyable}: 
	Taking $g = \id_Y$ in \eqref{eq:IID_naturality2} and using \Cref{prop:IID_inf_copy} shows that if $f$ is copyable, then it is also ``countably-copyable'', in the sense of
	\begin{equation}\label{eq:countably-copyable}
		f^\N \comp \id^{(\N)} = \id^{(\N)} \comp f.
	\end{equation}
\end{rem}

\subsubsection{Exchangeability}

For any injective function $\sigma : \N \to \N$, we can apply the universal property to obtain an action on the Kolmogorov power $X^\N$, which we denote by
\begin{equation}
	X^\sigma : X^\N \longrightarrow X^\N.
\end{equation}
It is the unique morphism whose finite marginals $\pi_n \comp X^\sigma$ are given by $X^{\sigma{|_{n}}} \comp \pi_m$ where we define $X^{\sigma{|_{n}}} \colon X^{m} \to X^n$ with $m \coloneq \max_{i \leq n}\{\sigma(i)\}$ to be the composite of discarding maps and structure isomorphisms which maps the $i$-th input to the $\sigma^{-1}(i)$-th output whenever $i$ is in the forward image of $\{1,\ldots,n\}$ under $\sigma$ and discards the $i$-th input otherwise.\footnote{Our choice of $m$ merely ensures that the domain contains all the components in the forward image. Any larger integer would work as well, as the corresponding morphism would simply discard all additional inputs.}
See \cite[Section~5]{fritzrischel2019zeroone} and \cite[Section~2.4]{chen2024aldoushoover} for more details.
By \Cref{lem:KP_dom_cop} and the fact that every $X^{\sigma{|_{n}}} \comp \pi_m$ is a deterministic morphism, we can also conclude that $X^\sigma$ is deterministic.

It is not hard to show that the assignment $\sigma \mapsto X^\sigma$ is contravariant.
Indeed, if $\tau : \N \to \N$ is another injective function, we have
\begin{equation}\label{eq:perm_contravariance}
	\pi_n \comp X^\sigma \comp X^\tau = X^{\sigma{|_{n}}} \comp \pi_m \comp X^\tau = X^{\sigma{|_{n}}} \comp X^{\tau{|_{m}}} \comp \pi_k = X^{(\tau \comp \sigma){|_{n}}} \comp \pi_k = \pi_n \comp  X^{\tau \comp \sigma},
\end{equation}
where $k \coloneq \max_{j \leq m}\{\tau(j)\}$ is also equal to $\max_{i \leq n}\{\tau \comp \sigma(i)\}$.
Thus, we conclude $X^\sigma \comp X^\tau = X^{\tau \comp \sigma}$.

We mainly consider $X^\sigma$ in the case where $\sigma$ is a \emph{finite} permutation, which is a bijection that leaves all but finitely many elements of $\N$ fixed.
In this way, we obtain an action of the finite permutation group $S_\infty$ on $X^\N$.
The morphisms fixed by this action are particularly relevant to the present work, as will become clearer in \Cref{sec:dF_obs_rep,sec:es_def}.
\begin{defi}\label{def:exchangeable}
	Let $\cC$ be a quasi-Markov category such that the Kolmogorov power $X^{\N}$ exists for an object $X$.
	A morphism $f: A \to X^{\N}\otimes Y$ is \newterm{exchangeable in the first factor} if 
	\begin{equation}
		\input{exchangeable.tikz}
	\end{equation}
	holds.
	If $Y = I$, i.e.\ if $f$ is of type $A \to X^{\N}$, we simply say that $f$ is \newterm{exchangeable}.
\end{defi}

As one would expect, IID morphisms as defined in \Cref{sec:IID} are exchangeable:
\begin{lem}[IID morphisms are exchangeable]\label{lem:fN_exchangeable}
	Let $\cC$ be a quasi-Markov category and let $X^\N$ be a Kolmogorov power. Then, for any $f : A \to X$, the IID morphism $f^{(\N)} : A \to X^\N$ is exchangeable.
\end{lem}
\begin{proof}
	By the universal property, we need to show that
	\begin{equation}
		\pi_n \comp X^\sigma \comp f^{(\N)} = \pi_n \comp f^{(\N)}
	\end{equation}
	for any $n \in \N$ and any finite permutation $\sigma$.
	
	The right-hand side is given by $f^{(n)}$, which is also the right-hand side of \Cref{eq:notation_fN}.
	Let us show that the left-hand side is the same.
	To this end, we use the definition of the marginals of $X^\sigma$ with the notation $X^{\sigma|_{n}}$ and $m\coloneq \max_{i \leq n}\{\sigma(i)\}$ being just as in the paragraphs above \Cref{def:exchangeable}:
	\begin{equation}\label{eq:fN_exchangeable_LHS}
		\begin{split}
			\pi_n \comp X^\sigma \comp f^{(\N)} &=  X^{\sigma{|_{n}}} \comp \pi_m \comp f^{(\N)} \\
			&= X^{\sigma{|_{n}}} \comp f^{(m)} \\
			&= f^{(n)}
		\end{split}
	\end{equation}
	The second equality in \eqref{eq:fN_exchangeable_LHS} follows by the definition of $f^{(\N)}$. 
	The third one follows by the fact that $X^{\sigma{|_{n}}}$ essentially discards $m-n$ of the inputs, yielding $f^{(n)}$ by quasi-totality of $f$, and permutes the remaining $n$ of them, which preserves $f^{(n)}$.
\end{proof}

Just as finite tensor products, infinite parallel composites are also permutation-covariant in the following sense.
\begin{lem}[Parallel composites are permutation-covariant]\label{lem:perm-covariance}
	Let $\cC$ be a positive and $\sigma$-continuous quasi-Markov category with Kolmogorov powers $X^\N$ and $Y^\N$.
	For any morphism $f : X \to Y$ in $\cC$ and a finite permutation $\nu : \N \to \N$ we have
	\begin{equation}\label{eq:perm-covariance}
		Y^\nu \comp f^\N = f^\N \comp X^\nu.
	\end{equation}
\end{lem}
\begin{proof}
	As usual, it suffices to prove that \Cref{eq:perm-covariance} holds when post-composed with an arbitrary finite projection.
	Applying $\pi_n$ to the left-hand side yields (once again we use the notation $X^{\nu|_{n}}$ and $m\coloneq \max_{i \leq n}\{\nu(i)\}$ as in the paragraphs above \Cref{def:exchangeable})
	\begin{equation}\label{eq:perm-covariance_1}
		\input{perm-covariance_1.tikz}
	\end{equation}
	where use the definition of $Y^\nu$ (and $X^\nu$) in steps 1 and 4; the definition of $f^\N$ and \Cref{lem:inf_tensor_domain} in step 2; the definition of $Y^{\nu|_{n}}$ (and $X^{\nu|_{n}}$) and quasi-totality of $f$ in step 3; and copyability of $\pi_m$ in step 4.
	By \Cref{lem:inf_tensor_domain} and \Cref{eq:inf_tensor_4}, the right-hand side of \Cref{eq:perm-covariance_1} can be simplified to
	\begin{equation}
		f^n \comp \pi_n \comp X^\nu \comp \bigwedge_{k \in \N} \dom{f^k \comp \pi_k}.
	\end{equation}
	By the fact that the infinite meet is of a descending sequence (\Cref{eq:inf_tensor_3}), we can also write it as 
	\begin{equation}\label{eq:perm-covariance_3}
		f^n \comp \pi_n \comp X^\nu \comp \bigwedge_{k \geq w} \dom{f^k \comp \pi_k} = f^n \comp \pi_n \comp \bigwedge_{k \geq w} \Bigl[ X^\nu \comp \dom{f^k \comp \pi_k} \Bigr].
	\end{equation}
	for an arbitrary $w \in \N$.	
	
	Let us pick the $w$ such that $\nu(i) = i = \nu^{-1}(i)$ holds for all $i \geq w$, which is possible since $\nu$ is a finite permutation.
	We now argue that $X^\nu$ and $\dom{f^k \comp \pi_k}$ commute for any $k \geq w$.
	By \Cref{eq:perm_contravariance}, the former has inverse $X^{\nu^{-1}}$ and thus we can express the composite $X^\nu \comp \dom{f^k \comp \pi_k}$ as 
	\begin{equation}\label{eq:perm-covariance_2}
		\input{perm-covariance_2.tikz}
	\end{equation}
	where, thanks to our choice of $w$, we can use $k = \max_{i \leq k}\{\nu^{-1}(i)\}$ and the fact that $X^{\nu^{-1}|_{k}}$ merely permutes the tensor factors of $X^k$ (i.e.\ it involves no discarding maps).
	Thus, we get 
	\begin{equation}\label{eq:perm-covariance_4}
		X^\nu \comp \dom{f^k \comp \pi_k} = \dom{f^k \comp \pi_k} \comp X^\nu.
	\end{equation}
	
	Combining \Cref{eq:perm-covariance_1,eq:perm-covariance_3,eq:perm-covariance_4}, we thus obtain
	\begin{equation}
		\begin{split}
			\pi_n \comp Y^\nu \comp f^\N &= f^n \comp \pi_n \comp \bigwedge_{k \geq w} \Bigl[ \dom{f^k \comp \pi_k} \comp X^\nu \Bigr] \\
			&= f^n \comp \pi_n \comp \bigwedge_{k \geq w} \Bigl[ \dom{f^k \comp \pi_k} \Bigr] \comp X^\nu \\
			&= \pi_n \comp f^\N \comp X^\nu,
		\end{split}
	\end{equation}
	which completes the proof by the universal property of the Kolmogorov power $Y^\N$.
\end{proof}

\subsection[Representability]{Representable Quasi-Markov Categories}
\label{sec:representable}

One can think of a Markov kernel either as a channel with a random outcome or as a deterministic map whose output is a probability distribution. 
Formulating this idea for Markov categories in general leads to the definition of representable Markov category~\cite{fritz2023representable}.
We now extend this notion to quasi-Markov categories.

\begin{defi}
	Let $\cC$ be a quasi-Markov category.
	\begin{enumerate}
		\item For $X \in \cC$, a \newterm{distribution object} is $PX \in \cC$ together with bijections
	 		\begin{equation}\label{eq:representability}
			 	\cC(A,X) \, \cong \, \cC_{\rm cop}(A,PX)
			 \end{equation}
			 natural in $A \in \cC_{\rm cop}$.
	 
		\item A quasi-Markov category is \newterm{representable} if every object has a distribution object.
	\end{enumerate}
\end{defi}

\noindent %FIXED indent
If $f : X\to Y$ is a morphism in a representable quasi-Markov category $\cC$, then we denote its counterpart under~\Cref{eq:representability} by $f^\sharp : X\longrightarrow PY$.

The correspondence of \Cref{eq:representability} can be used to define a functor $P$ that is right adjoint to the inclusion $\cC_{\rm cop} \hookrightarrow \cC$. 
This adjunction induces a monad $(P,\mu,\delta)$ on $\cC_{\rm cop}$, which we also denote by $P$, and $\cC$ can be seen as the Kleisli category of $P$.
The unit of the adjunction at $X$ is a copyable morphism
\begin{equation}
	\delta_X : X\longrightarrow PX,
\end{equation}
while the counit of the adjunction at $X$ is a morphism
\begin{equation}
	\samp_X : PX \longrightarrow X,
\end{equation}
namely the counterpart of the copyable $\id_{PX}$ under \Cref{eq:representability}.
The sampling morphism gives us a concrete implementation of Bijection \eqref{eq:representability}, i.e.\ we have $f = \samp_X \comp f^\sharp$ for every morphism $f : A \to X$.
In the opposite direction, we can write $f^\sharp$ as $Pf \comp \delta_A$.
Furthermore, since every morphism in $\cC(A,I)$ is trivially copyable by quasi-totality, we can conclude that $\samp_I : PI\cong I$ is an isomorphism, making the monad $P$ an affine monad.

One of the triangle identities is that each $\delta_X$ is a section of $\samp_X$.
As a consequence, $\delta_X$ is monic, and is thus total by \Cref{lem:mono_total}.
We leave open the question as to whether $\samp_X$ is total in general.
At present, we only know how to ensure this fact under the stronger assumption of \emph{observational representability} introduced in \Cref{def:obs_rep} below.
Whenever the sampling maps are total, the bijection ${\cC(A,X) \cong \cC_{\rm cop}(A,PX)}$ restricts to a bijection 
\begin{equation}\label{eq:representability_res}
	{\cC_{\rm tot}(A,X) \,\cong\, \cC_{\rm \det}(A,PX)}.
\end{equation}
for all objects $A$ and $X$.
This also means that a morphism $f$ is total if and only if $f^\sharp$ is.

\begin{exa}[Representability in measure-theoretic probability]
	$\Par{\BorelStoch}$ is representable, with $PX$ being the measurable space of probability measures on $X$, which is a standard Borel space again~\cite{giry}.
	Indeed one obtains a natural bijection as in \Cref{eq:representability} by sending a partial Markov kernel $f : A \to X$ to the partial measurable function
	\begin{align*}
		f^\sharp : A & \longrightarrow PX \\
			x & \longmapsto f(\ph|x),
	\end{align*}
	where the second row applies whenever $x \in D_f$ and $f^\sharp(x)$ is undefined otherwise.
	Then the above discussion instantiates to the following:
 \begin{itemize}
  \item The unit $\delta:X\to PX$ is the measurable function assigning to each $x\in X$ the Dirac delta $\delta_x\in PX$;
  \item The counit $\samp:PX\to X$ is the Markov kernel that takes a measure $\mu \in PX$ as input and returns a sample from $\mu$ as output.
	Formally, for all measurable $T \subseteq X$, we have
  \begin{equation}
  	\samp(T|\mu) \;=\; \mu(T)
  \end{equation} %FIXED Margins
  for all measurable subsets $T\subseteq X$.
	\item In particular, since $\samp$ is total, the monad $P$ on $\Par{\BorelMeas}\!\!\cong\!\! \Par{\BorelStoch}_{\rm cop}$ restricts to the Giry monad on $\BorelMeas$, and we thus recover the standard fact that $\BorelStoch$ is its Kleisli category \cite{giry}.
 \end{itemize}
\end{exa}

\subsection{Observational Representability and de Finetti Objects}
\label{sec:dF_obs_rep}

In \cite[Section~8.2]{moss2022probability}, a strengthening of representability has been introduced, which abstracts the idea that distinct measures can be distinguished by iterated sampling.
More details about this axiom can be found in \cite[Section A.4]{fritz2023supports}.
We generalize it now to the quasi-Markov setting.

\begin{defi}
	\label{def:obs_rep}
	A representable quasi-Markov category $\cC$ is \newterm{observationally representable} if for each $X$, the collection of morphisms $(\samp^{(n)} : PX \to X^n)_{n \in \N}$ is jointly monic.\footnote{Here, $\samp^{(n)}$ denotes the $n$-fold copy of $PX$ followed by $\samp^{n}$, in analogy to the notation from \Cref{eq:notation_fN}.}
\end{defi}

If $\cC$ also has countable Kolmogorov products, then it is easy to see that observational representability is equivalent to $\samp^{(\N)} : PX \to X^\N$ being monic.

\begin{lem}[Compatibility of enrichment and representability]
	\label{lem:rep_ext}
	In an observationally representable quasi-Markov category $\cC$ with countable Kolmogorov products,\footnote{$\Par{\BorelStoch}$ is one such category, as we elaborate on in \Cref{sec:theorems} (specifically, \Cref{thm:representability}).} every sampling morphism $\samp : PX \to X$ is total.
	Furthermore, we have
	\begin{equation}
		\label{eq:rep_ext}
		f \domext g \; \iff \; f^\sharp \domext g^\sharp
	\end{equation}
	for all morphisms $f, g : A \to X$ in $\cC$.
\end{lem}
As highlighted in \Cref{sec:representable}, this also means that the restricted bijection~\eqref{eq:representability_res} holds, and that a morphism $f$ is total if and only if $f^{\sharp}$ is.
\begin{proof}
	As noted above, the morphism $\samp^{(\N)} : PX \to X^\N$ is a monomorphism by observational representability, and thus total by \Cref{lem:mono_total}.
	Since the cone components of a Kolmogorov product are assumed to be total, we get that $\samp$, which can be obtained as ${\pi_1 \comp \samp^{(\N)}}$, is itself total.
	
	Totality of sampling means that for every morphism $g : A \to X$, we have
	\begin{equation}
		\label{eq:dom_of_det}
		\input{dom_of_det.tikz}
	\end{equation}
	which can be also written as $\dom{g} = \dom{g^\sharp}$ by using the $\dom{g}$ notation to refer to the domain of $g$ (\Cref{def:domain}) and recalling that $g=\samp \comp g^{\sharp}$ holds.
	Since $\delta$ is natural with respect to copyable morphisms, we have 
	\begin{equation}
		\label{eq:dom_sharp}
		(\dom{g})^\sharp = P\bigl( \dom{g} \bigr) \comp \delta_A = \delta_X \comp \dom{g}.
	\end{equation}
	Thus we get the desired \Cref{eq:rep_ext} via
	\begin{equation}
		\begin{split}
			f \domext g  \; &\iff \; f \comp \dom{g} = g \\
			&\iff \; Pf \comp (\dom{g})^\sharp = g^\sharp  \\
			&\iff \; f^\sharp \comp \dom{g} = g^\sharp \\
			&\iff \; f^\sharp \comp \dom{g^\sharp} = g^\sharp \\
			&\iff \; f^\sharp \domext g^\sharp
		\end{split}
	\end{equation}
	where the first and last equivalences are by definition, the second one is by Bijection \eqref{eq:representability}, the third one by \Cref{eq:dom_sharp}, and finally the fourth one by \Cref{eq:dom_of_det}.
\end{proof}

Our aim in the rest of this section is to show that observational representability can be derived from the existence of abstract de Finetti representations of exchangeable measures (and suitable compatibility conditions thereof).

\begin{defi}\label{def:dF_object}
	Let $\cC$ be a quasi-Markov category $\cC$ with countable Kolmogorov products and $X \in \cC$.
	Then a \newterm{de Finetti object} for $X$ is $QX \in \cC$ together with a morphism $\ell : QX \to X$ such that:
	\begin{enumerate}
		\item The morphism
			\begin{equation}
				\ell^{(\N)} : QX \to X^{\N}
			\end{equation} 
			makes $QX$ into the equalizer of all finite permutations $\permact{\sigma}{X} : X^{\N} \to X^{\N}$.
		\item This limit is preserved by ${\ph} \otimes \id_Y$ for every object $Y$.
	\end{enumerate}
\end{defi}

\noindent %FIXED indent
In~\cite{fritz2021definetti}, we proved that if $\cC$ is an \as{}-compatibly representable Markov category with countable Kolmogorov products and conditionals, then $PX$ together with $\samp : PX \to X$ satisfies the existence part of the universal property required for a de Finetti object.
In other words, this means that every exchangeable morphism $p : A\to X^\N$ factors across $\samp^{(\N)}$.
If the representability assumption is strengthened to observational representability, then every $X$ has a de Finetti object given by $PX$~\cite[Section~8.2]{moss2022probability}.\footnote{While the preservation of the limit by ${\ph} \otimes \id_Y$ was not considered in~\cite{fritz2021definetti}, this can be shown by the same arguments. See also~\cite[Theorem~B.2]{chen2024aldoushoover} for a similar result.}${}^{,}$\footnote{See~\cite[Proposition A.4.4]{fritz2023supports} for the proof that observational representability implies \as{}-compatible representability.}
The universal property is illustrated by the diagram:
\begin{equation}
\begin{tikzcd}
	&&& X^\N \arrow[dd,bend left,"X^\sigma"] \\
	A \arrow[urrr,bend left=20,"p"]\arrow[drrr,bend right=20,swap,"p"] \ar{rr}[near end]{\mu} && PX \ar{ur}[swap,inner sep=0.5mm]{\samp^{(\N)}} \ar{dr}{\samp^{(\N)}} \\
	&&& X^\N
\end{tikzcd}
\end{equation}

Taking $A = I$ in $\BorelStoch$, this says that every exchangeable probability measure $p$ on $X^\N$ can be expressed as a unique convex mixture of IID measures, i.e.\ it is in the form
\begin{equation}
	p(T_1\times\dots\times T_n \times X \times \dots ) \;=\; \int_{q\in PX} q(T_1)\cdots q(T_n)\,\mu(\mathrm{d}q)
\end{equation}
for a unique probability measure $\mu$ on $PX$, the so-called de Finetti measure of $p$.
Moreover, the analogous statement is true for every exchangeable Markov kernel $p : A\to X^\N$ for general $A$, which expresses the universal property in general.

The following result now can be thought of as an approximate converse.

\begin{thm}[Distributions from de Finetti representations]\label{thm:defin_obs}
	Let $\cC$ be a quasi-Markov category with countable Kolmogorov products such that every object has a de Finetti object.
	Then $\cC$ is observationally representable with distribution objects $QX$ and sampling morphisms $\ell : QX \to X$.
\end{thm}
\noindent %FIXED indent
This theorem and its proof are variations on \cite[Corollary 3.66]{fritz2023involutive}, which is a similar result for involutive Markov categories.
We provide the complete proof here to emphasize the role of quasi-totality.

\begin{proof}
	We need to show that the morphism $\ell : QX \to X$ gives a representation of the functor ${\cC\left( \ph ,X\right): \cC_{\rm cop}^\op \to \mathsf{Set}}$. 
	This is sufficient to also conclude observational representability, since $\ell^{(\N)}$ is monic as the universal morphism of an equalizer. 
	Thus we need to show that for every $A \in \cC$, the mapping	
	\begin{equation}
		\ell \comp \ph \; : \; \cC_{\rm cop} (A, QX)\to \cC (A,X)
	\end{equation}
	is a bijection.
	
	Let us first prove the injectivity, so suppose that $\ell \comp f = \ell \comp g$ holds for some copyable $f,g : Z \to A$.
	Then we have
	\begin{equation}\label{eq:dF_rep0}
		\input{dF_rep0.tikz}
	\end{equation}
	where we use \Cref{lem:IID_naturality} twice.
	Since $\ell^{(\N)}$ is monic, we obtain $f = g$ as required.
	 
	To prove surjectivity, we consider an arbitrary morphism $f : A \to X$ and aim to show that it factors through a \emph{copyable} morphism of type $A \to QX$. 
	By the assumed properties of the de Finetti object $QX$ as a limit, we can factor the exchangeable morphism $f^{(\N)} : A \to X^{\N}$ as
	\begin{equation}\label{eq:exch_factorization}
		f^{(\N)} = \ell^{(\N)} \comp \mu
	\end{equation}
	for a suitable morphism $\mu : A \to QX$. 
	In the rest of the proof, we show that $\mu$ must be copyable.	
	
	By the universal property of $X^\N$, \Cref{eq:exch_factorization} implies $f = \ell \comp \mu$, and thus we can rewrite~\eqref{eq:exch_factorization} as
	\begin{equation}
		\input{dF_rep4.tikz}
	\end{equation}
	which upon composition with the finite projection $\pi_n$ implies the same equality with $\N$ replaced by any finite $n \in \N$.
	Thinking of $n$ as an $n$-element set, we can use these equations together with the definition of Kolmogorov products and associativity of copy to obtain
	\begin{equation}\label{eq:dF_rep2}
		\input{dF_rep2_finite.tikz}
	\end{equation}
	where we omit an implicit isomorphism $X^n \otimes X^n \cong X^{n \sqcup n}$, the choice of which is irrelevant by the exchangeability of the IID morphisms involved.
	Since \Cref{eq:dF_rep2} holds for any $n$ and Kolmogorov products are preserved by tensor products, we can use the universal property to get
	\begin{equation}
		\input{dF_rep5.tikz}
	\end{equation}
	and finally
	\begin{equation}\label{eq:dF_rep6}
		\input{dF_rep6.tikz}
	\end{equation}
	
	Since $\ell^{(\N)} \tensor \id$ is monic by virtue of the limit $QX$ being preserved by tensor products, the morphism $\ell^{(\N)} \tensor \ell^{(\N)}$ is also monic. 
	Therefore, \Cref{eq:dF_rep6} implies that $\mu$ is copyable, which proves the surjectivity of the map $\ell \comp \ph$.
\end{proof}

\section{Empirical Sampling Morphisms}
\label{sec:emp_samp}

The key idea of this paper is that the limiting behaviour of infinite sequences of observations can be encapsulated by morphisms satisfying a few properties expressible in the language of quasi-Markov categories  (see \Cref{sec:es_def}). 
We devote \Cref{sec:es_BorelStoch,sec:es_integration} to the construction of instances of such morphisms, i.e.\ of partial Markov kernels between standard Borel spaces.
It is non-trivial to show that they satisfy the requisite properties, and we defer the technical proofs to \Cref{sec:proofs}.
Nevertheless, these two sections are markedly different in style from the rest of the article in their heavier use of measure-theoretic rather than category-theoretic language.

\subsection{The Idea and the Categorical Definition}
\label{sec:es_def}

Here we give our axiomatization of empirical sampling morphisms.
As mentioned in the Introduction, in measure-theoretic probability we want to associate to suitable sequences $(x_i)_{i \in \N}$ of points in a measurable space $X$
a probability measure formed out of limiting relative frequencies
\begin{equation}
	\lim_{n \to \infty} \frac{|\{ i \le n \mid x_i \in T \}|}{n}.
\end{equation}
When such a limiting probability measure exists, we call it \emph{empirical measure}.
Sampling from this empirical measure gives a random element of $X$.
As mentioned in the Introduction, a number of subtleties arise in the precise constructions.
For example, the limits above do not exist for every sequence. 
Therefore our constructions of empirical sampling morphisms will be \emph{partial} morphisms in $\BorelStoch$, or equivalently morphisms in $\Par{\BorelStoch}$.
This is why we formulate the general definition in terms of quasi-Markov categories rather than Markov categories.

The following two key properties were motivated and informally discussed in the Introduction.

\begin{defi}
	\label{def:es}
	Let $\cC$ be a quasi-Markov category with countable Kolmogorov products and $X \in \cC$.
	An \newterm{empirical sampling morphism} for $X$ is a morphism
	\begin{equation}
		\es \: : \: X^\N \longrightarrow X
	\end{equation}
	satisfying the following conditions:
	\begin{enumerate}
		\item\label{it:es_invariant} \newterm{Permutation invariance:} 
			For every finite permutation $\sigma$ of $\N$, $\es$ is invariant under pre-composition by the corresponding permutation ${X^\sigma : X^\N \to X^\N}$, we have
			\begin{equation}
				\es \comp X^\sigma \;=\; \es .
			\end{equation}
		\item\label{it:es_invariance} \newterm{Empirical adequacy:} If a morphism $f : A \to X^\N \otimes Y$ is exchangeable in the first factor,\footnote{This means that $(X^\sigma \otimes \id_Y) \comp f = f$ holds for every finite permutation $\sigma$.} then we have
			\begin{equation}
				\label{eq:es_invariance}
				\input{empirical_distribution_axiom.tikz}
			\end{equation}
	\end{enumerate}
\end{defi}
\noindent %FIXED indent
These two conditions are the general categorical formulations of the corresponding properties of empirical sampling mentioned in the Introduction.

\begin{rem}
	The converse of Property \ref{it:es_invariance} holds automatically: 
	Any $f$ satisfying \Cref{eq:es_invariance} is exchangeable in the first factor, since the $\es^{(\N)}$ appearing in the diagram on the left is exchangeable by \Cref{lem:fN_exchangeable}.
\end{rem}

Constructing interesting examples of empirical sampling morphisms, which we devote \Cref{sec:es_BorelStoch,sec:es_integration} to, requires a fair amount of measure-theoretic details.
Before doing so, let us note that we cannot expect empirical sampling morphisms to be natural transformations in $X$, at least not in $\BorelStoch$ (see \Cref{rem:es_not_natural}).
Nevertheless, they can be transported along copyable retracts, which can be thought of as a weak form of naturality.

\begin{lem}[Retracts of empirical sampling morphisms]\label{lem:transfer_es}
	Let $\cC$ be a positive and $\sigma$-continuous quasi-Markov category with Kolmogorov powers $X^\N$ and $Y^\N$.
	Let $\es_X : X^\N \to X$ be an empirical sampling morphism for $X$.
	Then for any copyable morphism $\iota : Y \to X$ with left inverse $\pi : X \to Y$, the morphism
	\begin{equation}
		\input{transfer_es.tikz}
	\end{equation}
	is an empirical sampling morphism for $Y$.
\end{lem}

\begin{proof}
	The permutation invariance of $\es_Y$ is a direct consequence of the naturality of braiding and the permutation invariance of $\es_X$.
	
	Concerning empirical adequacy, let $f : A \to Y^{\N} \tensor Z$ be exchangeable in the first factor.
	Then we have
	\begin{equation}\label{eq:transfer_es_proof2_noplate}
		\input{transfer_es_proof2.tikz}
	\end{equation}
	where the first step uses that $\iota^\N$ is copyable (\Cref{lem:inf_tensor_cop}) as well as \Cref{lem:IID_naturality}, the second is by empirical adequacy of $\es_X$ (exchangeability of $\iota^\N \comp f$ in the first factor follows from \Cref{lem:perm-covariance}), and the last one by \Cref{prop:inf_tensor_composition}.
\end{proof}

\subsection{Empirical Sampling Morphisms for Standard Borel Spaces}
\label{sec:es_BorelStoch}

Our goal now is to construct an empirical sampling morphism for every standard Borel space (\Cref{cor:standard_borel_es}).
By Kuratowski's theorem, it is enough to show that every finite set as well as $\N$ and $\R$, with their usual Borel $\sigma$-algebras, have empirical sampling morphisms.
In fact, since every standard Borel space is a measurable retract of $\R$, by \Cref{lem:transfer_es} it would even be sufficient to construct an empirical sampling morphism for $\R$ and we would automatically get it for every standard Borel space.
However, in order to progressively build up to $\R$ as the most challenging case, we start with the finite case, followed by $\N$, and only then deal with all the nuances in the continuous case.

Proofs of the theorems stated here can be found in \Cref{sec:proofs}.

\begin{defi}\label{def:es_finite}
	Let $F$ be a finite set equipped with the discrete $\sigma$-algebra. 
	For a given sequence $(x_i) \in F^\N$ and a subset $T \subseteq F$ we define
	\begin{equation}
		\label{eq:es_finite}
		\es_F \bigl( T | (x_i) \bigr) \coloneq \lim_{n \to \infty} \frac{\abs{ \Set{i \le n \given x_i \in T }}}{n}
	\end{equation}
	whenever this limit exists for all $T$ and leave it undefined otherwise.
\end{defi}

\begin{thm}[Concrete empirical sampling, the finite case]
	\label{lem:finite_es}
	Given any finite set $F$, the kernel $\es_F : F^\N \to F$ defined as in \Cref{eq:es_finite} is an empirical sampling morphism for $F$ in $\Par{\BorelStoch}$.
\end{thm}
\noindent %FIXED indent
In the infinite case we would like to use the same ideas, but this meets further subtleties. 
For example, the sequence $(1,2,3,\dots)$ over $\N$ has well-defined relative frequencies in the limit, since every single number has limiting frequency zero. 
Therefore, even though all the relative frequencies have well-defined limits, $\sigma$-additivity fails and we do not get a probability measure on $\N$.
So already for countably infinite spaces, the definition of $\es$ is not as straightforward as in the finite case, and only requiring the limits to exist is not enough.
To avoid situations such as the above case, we require additionally the following equivalent conditions.

\begin{prop}[Equivalent ways to ensure sensible empirical measures]\label{prop:eq_uniform}
	For every sequence $(x_i)\in \N^\N$ for which the limiting relative frequencies of singletons
	\begin{equation}
		\lim_{n \to \infty} \frac{|\{ i \le n \mid x_i = t \}|}{n}
	\end{equation}
	exist, the following properties are equivalent:
	\begin{enumerate}
		\item\label{it:es_tight}
		The sequence of finite empirical measures is tight:\footnote{For more details about the concept of tight sequences, see \cite[Section~25]{billingsley}.} For every $\eps > 0$ there is a $t \in \N$ such that for all $n \gg 1$, we have
		\begin{equation}
			\label{eq:es_tight}
			\frac{|\{ i \le n \mid x_i \ge t \}|}{n} < \eps.
		\end{equation}
		\item\label{it:es_uniform}
		The limits
		\begin{equation}
			\label{eq:es_uniform}
			\lim_{n \to \infty} \frac{|\{ i \le n \mid x_i \le t \}|}{n}
		\end{equation}
		exist \emph{uniformly} in $t \in \N$.\footnote{This means that there is a function $f \colon \N \to [0,1]$ such that for every $\eps>0$ there is $N\in\N$ with
		\begin{equation}
			\left|  \frac{|\{ i \le n \mid x_i \le t \}|}{n} - f(t) \right| < \eps  \qquad \forall n \ge N, t \in \N.
		\end{equation}}
		\item\label{it:es_norm} The limiting relative frequencies of singletons sum up to $1$,
		\begin{equation}
			\label{eq:es_sum}
			\sum_{t \in \N} \lim_{n \to \infty} \frac{|\{ i \le n \mid x_i = t \}|}{n} = 1.
		\end{equation}
	\end{enumerate}
	Whenever these conditions are satisfied, then for every $T \subseteq \N$, its probability coincides with its limiting relative frequency,
	\begin{equation}
		\label{eq:es_additive}
		\sum_{t \in T} \lim_{n \to \infty} \frac{|\{ i \le n \mid x_i = t \}|}{n} = \lim_{n \to \infty} \frac{|\{ i \le n \mid x_i \in T \}|}{n}.
	\end{equation}
\end{prop}
		
Informally, one may think of the tightness property as saying that no mass should ``escape to infinity''.
All three properties fail for the sequence $(1,2,3,\dots)$.

\begin{proof}
	The tightness condition (Property \ref{it:es_tight}) implies that the limits
	\begin{equation}
		\lim_{n \to \infty} \frac{|\{ i \le n \mid x_i \ge t \}|}{n}
	\end{equation}
	are uniform in $t$: 
	For a given $\eps > 0$ and large $t$, the right-hand side is less than $\eps$ for every $n$, and therefore we are effectively dealing with a finite range of $t$, where uniform convergence is trivial.
	Considering the complementary events (characterized by $x_i < t$) instead now shows that~\ref{it:es_uniform} holds.
	The converse follows by similar arguments.
	
	The equivalence with~\ref{it:es_norm} follows by standard conditions on the interchange of limit and summation~\cite{hildebrandt1912interchange}, since the tightness is exactly the statement that the sums $\sum_{t \in \N} \frac{|\{ i \le n \mid x_i = t \}|}{n}$ converge uniformly in $n$.
	
	Let us now show \Cref{eq:es_additive}, i.e.\ that the probability of any $T \subseteq \N$ coincides with its limiting relative frequency.
	The inequality $\le$ follows easily by taking the supremum over all finite subsets of $T$.
	The equality then follows upon combining the same inequality applied to $\N \setminus T$ combined with the normalization condition given by \Cref{eq:es_sum}.
\end{proof}

\begin{defi}\label{def:es_N}
	Let $\N$ be equipped with the discrete $\sigma$-algebra. 
	We define $\es_\N : \N^\N \to \N$ as the partial Markov kernel where:
	\begin{itemize}
		\item A sequence $(x_i)$ belongs to $\dom{\es_\N}$ if and only if it satisfies Property \ref{it:es_uniform} above.
		
		\item In this case, we take 
		\begin{equation}
			\label{eq:es_countable}
			\es_\N \bigl(T \,|\, (x_n) \bigr) \coloneq \lim_{n \to \infty} \frac{ \abs{ \Set{ i \le n \given x_i \in T }}}{n} .
		\end{equation}
		for every $T \subseteq \N$.
	\end{itemize}
\end{defi}

\begin{thm}[Concrete empirical sampling, the countable case]
	\label{lem:countable_es}
	The partial Markov kernel $\es_\N : \N^\N \to \N$ from \Cref{def:es_N} is an empirical sampling morphism for $\N$ in $\Par{\BorelStoch}$.
\end{thm}
\noindent %FIXED indent
Let us now turn to the construction of an empirical sampling morphism for $\R$.
As explained in the Introduction, we now cannot even expect the equation
\begin{equation}
	\label{eq:es_rel_freq}
	\es_\R \bigl( T | (x_i) \bigr) = \lim_{n \to \infty} \frac{ \abs{\{ i \le n \mid x_i \in T\}}}{n}
\end{equation}
to hold for all measurable $T$.
Indeed, sequences over $\R$ with all elements mutually distinct are generic, and so the failure of $\sigma$-additivity mentioned for case of the very special sequence $(1,2,3,\ldots)$ in $\N$ would now be a generic feature of elements of $\R^\N$.

We therefore need to relax the requirement that \Cref{eq:es_rel_freq} holds for all measurable sets $T$ to a smaller class of events.
In the case of $\N$ above, we characterize the viable sequences by requiring the relative frequencies of events $T = \{1,2, \ldots, t\}$ to converge uniformly in $t$ (\Cref{prop:eq_uniform}).
Similarly, let us consider subsets of the form $T=(-\infty,t]$ here. 
In other words, we view probability measures through their \newterm{cumulative distribution functions (CDF)}, which are monotone right continuous functions $F : \R \to [0,1]$ satisfying
\begin{equation}
	\label{limit_cdf}
	\lim_{t \to -\infty} F(t) = 0, \qquad \lim_{t \to \infty} F(t) = 1. 
\end{equation}
It is a basic fact that the probability measures $\mu$ on $\R$ are in bijection with these functions,
according to the assignment $F(t) \coloneqq \mu((-\infty,t])$ for all $t \in \R$.
We can therefore define a probability measure $\es_\R(\ph | (x_i))$ on $\R$ by specifying its CDF, as long as we ensure that this function has the requisite properties.
Moreover, the right continuity determines such a function uniquely from its values on $\Q$, as every monotone and right continuous function on $\Q$ satisfying~\Cref{limit_cdf} extends uniquely to such a function on $\R$.
Therefore from now on, we consider the domain of CDFs to be $\Q$.
Its countability facilitates our arguments for measurability.

This discussion motivates the following construction of an empirical sampling morphism for $\R$.

\begin{defi}\label{def:es_R}
	Let $\R$ be equipped with its Borel $\sigma$-algebra. 
	We define $\es_\R : \R^\N \to \R$ as the partial Markov kernel where:
	\begin{itemize}
		\item A sequence $(x_i)$ belongs to $\dom{\es_\R}$ if and only if the limits
			\begin{equation}
				\label{eq:es_cdf}
				\lim_{n \to \infty} \frac{ \abs{\{ i \le n \mid x_i \le t\}}}{n}
			\end{equation}
			exist uniformly in $t \in \Q$.
		\item In this case, we define $\es_\R(\ph|(x_i))$ to be the unique probability measure on $\R$ satisfying
			\begin{equation}
				\label{eq:es_intervals}
				\es_\R\bigl(T \,|\, (x_i)\bigr) = \lim_{n \to \infty} \frac{ \abs{\{ i \le n \mid x_i \in T\}}}{n}
			\end{equation}
			for all intervals $T \subseteq \R$.
	\end{itemize}
\end{defi}
\noindent %FIXED indent
A priori it is not clear whether Equation~\eqref{eq:es_intervals} defines a probability measure at all, since it does not construct the measure explicitly.
While the uniqueness of the measure is clear from the fact that intervals form a generating $\pi$-system for Borel sets, the existence is most easily seen via standard properties of cumulative distribution functions.
We spell this out in the proof of the following result (which can be found in \Cref{sec:proofs}).
As detailed there, the uniformity of the limit plays a role and indeed guarantees that \eqref{eq:es_cdf} defines a CDF.

\begin{thm}[Concrete empirical sampling, the continuous case]
	\label{lem:real_es}
	The partial Markov kernel $\es_\R : \R^\N \to \R$ from \Cref{def:es_R} is an empirical sampling morphism for $\R$ in $\Par{\BorelStoch}$.
\end{thm}

\begin{exa}\label{ex:es_domain}
	We illustrate $\dom{\es_\R}$ with two simple examples, and compare it with the weak convergence of empirical measures as in \eqref{eq:finite_empirical_distribution}.
	\begin{enumerate}
		\item The sequence $(1, 1/2, 1/3, \ldots)$ has well-defined limiting relative frequencies given by
		\begin{equation}
		\lim_{n \to \infty} \frac{|\{ i \le n \mid x_i \le t\}|}{n} = \begin{cases}
			1 & \text{if $t > 0$}, \\
			0 & \text{if $t \le 0$}.
		\end{cases}
		\end{equation}
		But the convergence is not uniform in $t$.
		Indeed this limit function is not a CDF, and therefore does not correspond to a probability measure.
		The sequence of empirical measures, however, converges weakly to $\delta_0$, which gives empirical probability $1$ to the singleton event $\{0\}$, despite that fact that $0$ never appears in the sequence.
		
		\item The similar sequence $(-1, -1/2, -1/3, \ldots)$ also has well-defined limiting relative frequencies given by
		\begin{equation}
		\lim_{n \to \infty} \frac{|\{ i \le n \mid x_i \le t\}|}{n} = \begin{cases}
			1 & \text{if $t \ge 0$}, \\
			0 & \text{if $t < 0$},
		\end{cases}
		\end{equation}
		and again the finite empirical measures converge weakly to $\delta_0$.
		However, now the convergence of CDFs \emph{is} uniform in $t$, and therefore this sequence belongs to $\dom{\es_\R}$ and its empirical measure is $\delta_0$.
		Once again $0$ does not appear in the sequence.
	\end{enumerate}
\end{exa}
\noindent %FIXED indent
This shows that $\es_\R$ is not invariant under the isomorphism $x \mapsto -x$ of $\R$ and thus that empirical sampling morphisms are not unique.
In fact, we construct another inequivalent empirical sampling morphism for $\R$ in \Cref{sec:es_integration}.

It is conceivable that the domain of $\es_\R$ could be extended to include both of the above sequences, e.g.\ by requiring the weak convergence of finite empirical measures instead of the uniform convergence of the CDF.
The considerations of Austin and Panchenko~\cite[Theorem~3]{austin2014hierarchical} suggest that this could indeed be the case on $[0,1]$, but we have not explored this possibility further so far.

\begin{rem}
	In \emph{discrepancy theory}~\cite{aistleitner2011discrepancy}, an important topic is the study of sequences $(x_i)$ whose finite empirical measures converge quickly to the uniform distribution on $[0,1]$.
	In this context, one studies the speed of convergence also in terms of the sup-norm-distance between the corresponding CDFs.
	This matches our requirement of uniform convergence in~\Cref{eq:es_cdf}.
\end{rem}

\begin{cor}
	\label{cor:standard_borel_es}
	In $\Par{\BorelStoch}$, every standard Borel space $X$ admits an empirical sampling morphism.
\end{cor}

\begin{proof}
	By Kuratowski's theorem, every standard Borel space is either finite or isomorphic to $\N$ with the discrete $\sigma$-algebra or to $\R$ with the Borel $\sigma$-algebra. 
	By \Cref{lem:finite_es,lem:countable_es,lem:real_es}, these cases are covered.
	Moreover, by \Cref{lem:transfer_es}, we know that empirical sampling morphisms can be transferred along isomorphisms.
\end{proof}

\begin{rem}\label{rem:es_relation}
	By applying \Cref{lem:transfer_es} to the empirical sampling morphism $\es_\R$, we can construct an empirical sampling morphism on $\N \subseteq \R$ and on any finite subset $F \subseteq \R$.
	This exactly reproduces the empirical sampling morphisms $\es_F$ and $\es_\N$ from \Cref{def:es_finite,def:es_N} respectively.
	In particular, it is easy to check that any sequence in $\N$ belongs to the domain of $\es_\N$ if and only if it is in the domain of $\es_\R$.
\end{rem}

We leave open the question of whether empirical sampling morphisms exist for other kinds of measurable spaces, and what kind of structure is needed in order to construct them.

Before we conclude this section, let us also discuss why the naturality of empirical sampling morphisms cannot hold in $\Par{\BorelStoch}$.
We already saw that our construction of $\es_\R$ is not invariant under isomorphisms, since \Cref{ex:es_domain} shows that its domain is not invariant under the map $x \mapsto -x$.

\begin{rem}[Non-naturality of empirical sampling morphisms]
	\label{rem:es_not_natural}
	Given that an empirical sampling morphism exists on every object of $\BorelStoch$, one might ask whether these constitute can be constructed as a natural transformation, in the sense that the square
	\begin{equation}
		\label{eq:es_naturality}
		\begin{tikzcd}[row sep=small]
			{X^\N} && X \\
			& \\
			{Y^\N} && Y
			\arrow["{\es_Y}"', from=3-1, to=3-3]
			\arrow["{f^\N}"', from=1-1, to=3-1]
			\arrow["f", from=1-3, to=3-3]
			\arrow["{\es_X}", from=1-1, to=1-3]
		\end{tikzcd}
	\end{equation}
	commutes for every measurable map $f : X \to Y$.
	However, naturality is too much to ask for in general, because then taking $f \coloneqq \discard_X$ shows that $\es_X$ would have to be total.
	Although this could be achieved at least for finite $X$, e.g.~by using ultralimits in~\Cref{eq:es_finite}, this is \emph{not what we want}, since the very existence of the limit is an important part of results like the Glivenko--Cantelli theorem (\Cref{thm:glivenko}).

	More plausibly, one might hope for empirical sampling morphisms to satisfy \emph{lax naturality} instead, by which we mean
	\begin{equation}
		\label{eq:es_lax_naturality}
		\es_Y \comp f^\N \sqsupseteq f \comp \es_X 
	\end{equation}
	for every $f$, where $\sqsupseteq$ is the extension relation of \Cref{def:dom_ext}.
	In $\Par{\BorelStoch}$ and for total $f$, this states that whenever a sequence $(x_i)$ has a well-defined empirical measure, then the image sequence $(f(x_i))$ has an empirical measure as well, and it is given by the pushforward measure.
	\Cref{lem:relu_lax_naturality} gives one example where this holds.
	However, lax naturality must fail in general, since it entails the empirical sampling morphism on $\R$ to be total, and we already noted that this is not desirable.
	
	To prove this, take any diffuse probability measure $\mu$ on $\R$. 
	Sampling from $\mu$ produces a sequence of mutually distinct elements with probability $1$, 
	and therefore there must be at least one such sequence $(x_i)$ in $\dom{\es}$ by empirical adequacy.
	Take any other sequence $(y_i) \in \R^\N$.
	There exists a measurable function $f$ satisfying $f(x_i) = y_i$ for all $i$.
	Lax naturality thus implies $(y_i) \in \dom{\es}$ and consequently that $\es$ is total.
\end{rem}

\subsection[Empirical Averaging]{Empirical Averaging as Integration Against an Empirical Measure}
\label{sec:es_integration}

Even though we have now shown that every object in the quasi-Markov category $\Par{\BorelStoch}$ has an empirical sampling morphism, we have also seen how delicate is the choice of its domain.
We can liken this choice to the choice of a topology for a space{\,---\,}it specifies for which sequences do we consider their limiting empirical measure to be well-defined.
In intuitive terms, enlarging the domain makes it easier for statements involving $\es$ to hold with high probability (e.g.\ with probability $1$).
On the other hand, restricting the domain can encode additional properties or assumptions about the sequences considered.

Given a sequence $(x_i)\in \dom{\es}$ and any bounded measurable function $f : X \to \R$, one might hope that its integral is given by the empirical average of $f$ as in the right-hand side of
\begin{equation}
	\label{eq:es_integral}
	\int_{y \in X} f(y) \, \es \bigl(\mathrm{d}y \,|\, (x_i) \bigr) \stackrel{?}{=} \lim_{n \to \infty} \frac{1}{n} \sum_{i=1}^n f(x_i).
\end{equation}
The reason for why one might expect this is because its finite analogue trivially holds: 
The expectation value of $f$ with respect to the finite empirical measure  from \eqref{eq:finite_empirical_distribution} is the empirical average $\frac{1}{n} \sum_{i=1}^n f(x_i)$.

It is not hard to show that also the full version of \Cref{eq:es_integral} holds for any finite $X$.
However, it is too much to hope for in general: 
For instance, take $X = \R$ and any sequence whose empirical measure has no atoms, and let $f$ be the indicator function of the support $\{x_i\}_{i \in \N}$ of the sequence.
Then the left-hand side of \Cref{eq:es_integral} is $0$ while its right-hand side is $1$.

The aim of this section is to show that the domain of $\es_\R$ from \Cref{def:es_R} can be restricted in such a way that
\begin{enumerate}
	\item \Cref{eq:es_integral} holds for $f = \id$, and
	\item the permutation invariance and empirical adequacy axioms still hold.
\end{enumerate}
\noindent %FIXED indent
To begin our analysis, note the following weaker statement that is true also for the original $\es_\R$, and with the same proof also for $\es_\N$ (in which case the continuity is trivial).

\begin{prop}[Sufficient conditions for the empirical averaging property]
	\label{prop:es_integral}
	\Cref{eq:es_integral} holds for $\es_\R$ from \Cref{def:es_R} and for every function ${f : \R \to \R}$ for which the left and right limits
	\begin{equation}
		\lim_{s \nearrow t} f(s), \qquad \lim_{s \searrow t} f(s)
	\end{equation}
	exist and are finite for every $t \in \R$, as well as for $t = +\infty$ (left limit only) and $t = -\infty$ (right limit only).
\end{prop}
\noindent %FIXED indent
The following proof uses the fact that every such function is a uniform limit of step functions.
Thus, every such function is measurable and bounded.

\begin{proof}
	By definition of $\es_\R$, the desired~\Cref{eq:es_integral} holds for $f = 1_{(-\infty,t]}$ for all $t \in \R$.
	We show first that it also holds for indicator functions of the form $f = 1_{(-\infty,t)}$.
	This is because by the assumed uniform convergence, the left limits of the finite empirical CDFs converge to the left limit of the limiting CDF.
	By combining these two cases for $f$ and using linearity, it follows that~\Cref{eq:es_integral} holds for all step functions.
	
	Moreover, clearly~\Cref{eq:es_integral} is preserved under uniform limits.
	Therefore it holds for all functions that can be written as uniform limits of step functions.
	On compact intervals, it is a standard fact that such functions are exactly those with left and right limits~\cite[Section~VII.6]{dieudonne1969analysis}.
	The same arguments also prove the relevant statement in the case of $\R$ (although only the converse inclusion is needed).
\end{proof}

As a special case of \Cref{prop:es_integral}, note that every continuous function $f$ which converges at both $\pm \infty$ necessarily satisfies \Cref{eq:es_integral}.
However, we would like the equation to hold also for suitable unbounded $f$, and in particular for $f = \id_\R$.
The reason is that the right-hand side of \Cref{eq:es_integral} then becomes the empirical average appearing in the statement of the strong law of large numbers.
By relating it to the mean of the empirical measure, we will show that our synthetic strong law of large numbers (\Cref{cor:genericLLN}) indeed recovers the standard result.

As the following example shows, this is not true for the empirical sampling morphisms from \Cref{sec:es_BorelStoch}.

\begin{exa}[A sequence with no empirical average]
	\label{ex:escaping_sequence}
	Consider the sequence $(x_n)$ over $\N$ given by
	\begin{equation}
		x_n \coloneqq \begin{cases}
			n + 1 & \text{if $n$ is a power of $2$}, \\ 
			1 & \text{otherwise}.
		\end{cases}
	\end{equation}
	Then the limit
	\begin{equation}
		\lim_{n \to \infty} \frac{|\{ i \le n \mid x_i \le t\}|}{n} 
	\end{equation}
	of relative frequencies is equal to $1$ for every $t \in \N$,\footnote{Recall that $\N$, by convention, does not contain $0$ (\Cref{foot:N}).} since almost all elements of the sequence are one.
	The convergence is uniform by monotonicity in $t$.
	Therefore this sequence lies in $\dom{\es_\N}$ and we have
	\begin{equation}
		\es_\N \bigl( \ph \,|\, (x_n) \bigr) = \delta_1.
	\end{equation}
	By \Cref{rem:es_relation}, the same holds when considering $(x_n)$ as a sequence in $\R$, and for the map $\es_\R$.
	
	Consider now the condition \Cref{eq:es_integral} for $f=\id$, which gives
	\begin{equation}\label{eq:emp_av}
		\int_{y \in \R} y \, \es_\R \bigl(\mathrm{d}y \,|\, (x_n) \bigr) \stackrel{?}{=} \lim_{n \to \infty} \frac{1}{n} \sum_{i=1}^n x_i.
	\end{equation}
	The left-hand side evaluates to $1$, while the limit on the right-hand side does not even exist: 
	The sequence of averages oscillates between $1$ and $3/2$ indefinitely.
\end{exa}

Let us rectify this issue by constructing a revised empirical sampling morphism $\esav_{\R}$.
It coincides with $\es_\R$, except in that its domain is restricted such that \Cref{eq:emp_av} is guaranteed to hold.

\begin{defi}
	\label{def:es_av}
	The partial Markov kernel $\esav_{\R} : \R^\N \to \R$ is defined just as $\es_\R$ from \Cref{def:es_R} but for the fact that its domain is restricted to those sequences for which additionally the sequence $\left(\frac{1}{n} \sum_{i=1}^n |x_i|\right)$ converges in $\R \cup \{\infty\}$ and satisfies
	\begin{equation}\label{eq:es_integral_Eabs}
		\int_\R \abs{y} \, \esav_{\R} \bigl( \mathrm{d}y \mid (x_n) \bigr) = \lim_{n \to \infty} \frac{1}{n} \sum_{i=1}^n |x_i|	 
	\end{equation}
	and, moreover, such that if both sides of \Cref{eq:es_integral_Eabs} are finite, then the sequence $\left( \frac{1}{n} \sum_{i=1}^n x_i \right)$ converges and satisfies
	\begin{equation}\label{eq:es_integral_E}
		\int_\R y \, \esav_{\R} \bigl( \mathrm{d}y \mid (x_n) \bigr) = \lim_{n \to \infty} \frac{1}{n} \sum_{i=1}^n x_i.		
	\end{equation}
\end{defi}
We leave it understood that we include sequences for which both sides of \Cref{eq:es_integral_Eabs} are infinite.
The sequence from \Cref{ex:escaping_sequence} is excluded because the limit on the right-hand side does not exist.

\begin{thm}[Concrete empirical sampling with well-behaved empirical average]
	\label{prop:tightened_es}
	The partial Markov kernel $\esav_{\R}$ from \Cref{def:es_av} is an empirical sampling morphism in $\Par{\BorelStoch}$.
\end{thm}
The proof can again be found in \Cref{sec:proofs}.

\section{Synthetic Laws of Large Numbers}
\label{sec:theorems}

We now present our main synthetic results that follow from the existence of empirical sampling morphisms plus two auxiliary assumptions.
Remarkably, these conditions imply representability of the category and the de Finetti theorem (\Cref{cor:dF}).
We then prove a synthetic Glivenko--Cantelli theorem (\Cref{thm:lln}).
One of its immediate consequences is a synthetic strong law of large numbers (\Cref{cor:genericLLN}).

\subsection{Our Assumptions}
\label{sec:assumptions}

Throughout this section, we assume the following.
\begin{asm}
	\label{assumptions}
	Let $\cC$ be a quasi-Markov category such that 
	\begin{itemize}
		\item $\cC$ is positive and $\sigma$-continuous.
		\item Countable Kolmogorov products exist in $\cC$.
		\item All balanced idempotents \cite[Definition 4.1.1]{fritz2023supports} split in $\cC$.
		\item Every object $X$ has an empirical sampling morphism $\es : X^\N \to X$.
	\end{itemize}
\end{asm}
\noindent %FIXED indent
As the reader may observe, balanced idempotents are the only concept whose details we defer to the literature. 
This is because they play a tangential role in the present work and appear only once, in the proof of \Cref{thm:representability}, to ensure representability.

\begin{rem}
	\label{rem:borelstoch_assumptions}
	\Cref{assumptions} holds in $\Par{\BorelStoch}$.
	Indeed the existence of empirical sampling morphisms is \Cref{cor:standard_borel_es}.
	The $\sigma$-continuity was shown in \Cref{prop:sigma_parborelstoch}.
	For the other properties, we employ the abstract arguments in \cite{shahmohammed2025partial}, which show that if a given Markov category $\cD$ satisfies one of these properties then so does its partialization \cite{shahmohammed2025partial}.
	The fact that these properties hold in $\BorelStoch$ has been shown in the past: positivity in \cite[Example 11.25]{fritz2019synthetic}, existence of countable Kolmogorov products in \cite[Example 3.6]{fritzrischel2019zeroone}, and the splitting of all idempotents in \cite[Corollary 4.4.5]{fritz2023supports}.
\end{rem}

\begin{rem}
	The splitting of balanced idempotents has been proven synthetically for any Markov category in \cite[Theorem 4.4.3]{fritz2023supports} under the assumptions of positivity, observational representability, and the equalizer principle defined therein.
	If one generalized this result to quasi-Markov categories, one would obtain alternative synthetic results without splitting of idempotents among its assumptions.
	However, we would find this a bit less satisfactory, since one \emph{conclusion} of our synthetic results in this section is the observational representability of $\cC$ (as a consequence of \Cref{thm:representability}).
\end{rem}

\subsection{Representability and the de Finetti Theorem}
\label{sec:dF_synthetic}

Already our categorical formulation of the Glivenko--Cantelli theorem will involve a distribution object $PX$.
Even though its existence is not part of \Cref{assumptions}, it is in fact ensured.
In essence, we identify the distribution object $PX$ via the splitting of a suitable idempotent.

As we show for our measure-theoretic construction of $\es_\R$ in \Cref{eq:es_random_permutation}, the \newterm{resampling} morphism ${\es^{(\N)} : X^\N \to X^\N}$ may be seen as implementing a uniformly random permutation.
Based on this, we might expect that it factors through the equalizer of all permutations on $X^\N$ in general, which by definition is a de Finetti object (\Cref{def:dF_object}) and therefore also a distribution object (\Cref{thm:defin_obs}).
Let us now make this argument precise.

\begin{thm}[Representability]\label{thm:representability}
	Given \Cref{assumptions}, $\cC$ is an observationally representable quasi-Markov category.
	In particular, for every object $X$ we have a distribution object $PX$ and a sampling morphism ${\samp : PX \to X}$.
\end{thm}
\begin{proof}
	First, for any finite permutation $\sigma$ we have 
	\begin{equation}\label{eq:resamp_balanced}
		\input{resamp_balanced.tikz}
	\end{equation}
	where the first equality holds because $X^\sigma$ is deterministic with inverse $X^{\sigma^{-1}}$, the second one follows by permutation invariance of $\es$ together with \Cref{lem:IID_naturality}, and the last one is just the exchangeability of $\es^{(\N)}$, which is an instance of \Cref{lem:fN_exchangeable}.
	Therefore, the morphism on the right-hand side of \Cref{eq:resamp_balanced} is exchangeable in its first $X^\N$ output, and therefore empirical adequacy of $\es$ implies 
	\begin{equation}\label{eq:resamp_idempotent}
		\input{resamp_idempotent.tikz}
	\end{equation}
	This is the defining equation of a balanced idempotent.
	By assumption, we obtain a splitting, i.e.\ a factorization as $\es^{(\N)} = \iota \comp \pi$ satisfying $\pi \comp \iota = \id_{E}$ for some object $E$.
	
	Let us now show that $\iota$ is IID, i.e.\ that it is equal to $\ell^{(\N)}$ for some $\ell : E \to X$.
	Firstly, we introduce a diagrammatic notation for the isomorphism $X^\N \cong X \otimes X^\N$ given in terms of the compatible family 
	\begin{equation}\label{eq:bifurcation}
		\input{bifurcation.tikz}
	\end{equation}
	of morphisms of type $X^\N \to X \otimes X^{n}$ for any $n \in \N$ and $\pi_n$ defined as in \Cref{def:kolmogorov}.
	For any morphism $f : A \to X$, we then have 
	\begin{equation}\label{eq:extract_copy}
		\input{extract_copy.tikz}
	\end{equation}
	which follows from
	\begin{equation}\label{eq:extract_copy_2}
		\input{extract_copy_2.tikz}
	\end{equation}
	that holds for any $n \in \N$.
	Then we have
	\begin{equation}\label{eq:iota_IID}
		\input{iota_IID.tikz}
	\end{equation}
	where the first and third equality follow from the splitting of $\es^{(\N)}$, second one is a version of \Cref{eq:extract_copy} for $f = \es$, and the last one is an application of the positivity assumption to the deterministic morphism $\pi \comp \iota = \id_{E}$. 
	We can now apply the same argument repeatedly to obtain
	\begin{equation}
		\input{det_displaycondind5.tikz}
	\end{equation}
	for any $n \in \N$.
	By the universal property of the Kolmogorov product, we thus get $\iota = \ell^{(\N)}$ for ${\ell \coloneqq \es \comp \iota}$.
	
	Let us now show that $E$ is a de Finetti object with universal arrow $\iota : E \to X^\N$.
	By $\iota = \ell^{(\N)}$, we already know that $\iota$ is invariant under finite permutations on $X^\N$.
	Consider now a generic morphism $f : A \to X^\N \tensor Y$ that is invariant under all finite permutations of $X^\N$, i.e.\ one that is exchangeable.
	By the empirical adequacy of $\es$ and the above splitting of $\es^{(\N)}$, we have
	\begin{equation}
		\label{eq:es_defin}
		\input{es_defin.tikz}
	\end{equation}
	i.e.\ $f$ factors through $\iota \otimes \id_Y$ as required.
	
	Observational representability now follows from \Cref{thm:defin_obs}.
	In particular, we can consider $E$ as the distribution object $PX$ of $X$ with sampling morphism $\ell = \es \comp \iota : E \to X$.
\end{proof}

In addition to being a distribution object, the above proof in fact show that $E$ is a de Finetti object.
This fact may be interpreted as a synthetic de Finetti theorem.

\begin{cor}[Synthetic de Finetti theorem]
	\label{cor:dF}
	Given \Cref{assumptions}, the distribution object of any $X \in \cC$ is also a de Finetti object.
\end{cor}
\noindent %FIXED indent
However, this differs significantly from the synthetic de Finetti theorems for Markov categories from earlier works, namely \cite[Theorem 4.4]{fritz2021definetti} and \cite[Theorem 2.11]{chen2024aldoushoover}.
These two establish conditional independence of exchangeable morphisms, while here we have implicitly \emph{assumed} this conditional independence in the empirical adequacy axiom.
The non-trivial statement established by \Cref{cor:dF} is the fact that exchangeable morphisms must factor through the distribution object $PX$, and the universal property by virtue of the uniqueness of this decomposition.

\subsection{The Synthetic Glivenko--Cantelli Theorem}
\label{sec:lln_synthetic}

With the existence of distribution objects ensured, we can now aim for our main theorem, the synthetic Glivenko--Cantelli theorem.
An important stepping stone towards it is the following result, which shows that the resampling morphism ${\es^{(\N)} : X^\N \to X^\N}$ is a split idempotent given by infinitely many samples from the empirical measure.

\begin{prop}[Splitting of resampling]
	\label{lem:pre_lln}
	Given \Cref{assumptions}, the resampling morphism $\es^{(\N)}$ splits as\footnote{Recall that $f^{\sharp} : A \to PX$ denotes the counterpart of $f : A \to X$ via Bijection \eqref{eq:representability}.}
	\begin{equation}\label{eq:es_splitting}
		\es^{(\N)} = \samp^{(\N)} \comp \es^\sharp, \qquad
		\es^\sharp \comp \samp^{(\N)} = \id_{PX}.
	\end{equation}
\end{prop}

\begin{proof}
		By the representability proven in \Cref{thm:representability}, we have $\es = \samp \comp \es^{\sharp}$.
			Since the morphism ${\es^{\sharp} : X^\N \to PX}$ is copyable by definition, we get 
			\begin{equation}
				\samp^{(\N)} \comp \es^{\sharp} = \es^{(\N)},
			\end{equation}
			which is the first equation in \eqref{eq:es_splitting}.
			
			To prove the second equation in \eqref{eq:es_splitting}, consider
			\begin{equation} 
				\label{eq:lln_proof_1}
				\samp^{(\N)} \comp \es^{\sharp} \comp \samp^{(\N)} = \es^{(\N)} \comp \samp^{(\N)} = \samp^{(\N)},
			\end{equation}
			where the second equality follows by empirical adequacy (\Cref{eq:es_invariance}), since $\samp^{(\N)}$ is IID and hence exchangeable.
			Since $\cC$ is observationally representable (\Cref{thm:representability}), $\samp^{(\N)}$ is monic, and therefore we get $\es^\sharp \comp \samp^{(\N)} = \id$ from \Cref{eq:lln_proof_1}.
\end{proof}

Plugging in the splitting of the resampling morphism into the empirical adequacy axiom of the empirical sampling morphism, i.e.\ into \Cref{eq:es_invariance}, we obtain the following corollary of \Cref{lem:pre_lln}:
Every morphism $f : A \to X^\N \otimes Y$ that is exchangeable in $X^\N$ satisfies
\begin{equation}
	\label{eq:de_finetti_rep}
	\input{de_finetti_rep.tikz}
\end{equation}
In particular, for $Y = I$, \Cref{eq:de_finetti_rep} says that the de Finetti measure of $f$ is given by $\es^\sharp \comp f$.

Yet another corollary is that the resampling morphism is not just a balanced idempotent as shown in the proof of \Cref{thm:representability}, but in fact a \newterm{strong idempotent}, which means that it satisfies 
\begin{equation}
	\input{resamp_strong_idempotent.tikz}
\end{equation}
This follows from the splitting proven above and \cite[Theorem~4.2.1(ii)]{fritz2023supports},\footnote{That result is stated for Markov categories, but the proof works in quasi-Markov categories just as well.} since the projection of the splitting is $\es^\sharp$ by \Cref{lem:pre_lln} and thus copyable. 

Our Glivenko--Cantelli theorem for quasi-Markov categories is also not hard to derive once we have \Cref{lem:pre_lln}.
Its name is justified by the measure-theoretic Glivenko--Cantelli theorem, which we recover as \Cref{thm:glivenko}.

\begin{thm}[Synthetic Glivenko--Cantelli theorem]\label{thm:lln}
	Consider an arbitrary exchangeable morphism ${f : A \to X^\N}$ and an arbitrary morphism ${p : A \to X}$ in a quasi-Markov category.
	Given \Cref{assumptions}, they satisfy
	\begin{equation}
		\label{eq:lln_1}
		\input{lln_3.tikz} \qquad\quad \text{and} \qquad\quad \input{lln_1.tikz}
	\end{equation}
	respectively.
\end{thm}

It is worth noting that the first equation establishes that $\es^\sharp$ is a sufficient statistic for $f$ in the sense of~\cite[Definition~14.3]{fritz2019synthetic}.

\begin{proof}
	Let us start with the first equation in \eqref{eq:lln_1}, since the second one can be then derived as a consequence.
	To this end, we calculate
	\begin{equation}
		\label{eq:es_splitting_1}
		\input{es_splitting_1.tikz}
	\end{equation}
	where the first equality is by \Cref{eq:de_finetti_rep}, and the second by the positivity axiom applied to the morphism $\es^\sharp \comp \samp^{(\N)}$, which is equal to $\id_{PX}$ by \Cref{lem:pre_lln}.
			
	Let us now move to the second equation in \eqref{eq:lln_1}, applicable to the specific exchangeable morphism given by $p^{(\N)}$.
	Since we can write $p = \samp \comp p^\sharp$ by representability and $p^\sharp$ is copyable, we have the first equality in
	\begin{equation}
		\label{eq:lln_proof_2}
		\es^\sharp \comp p^{(\N)} = \es^\sharp \comp \samp^{(\N)} \comp p^\sharp = p^\sharp,
	\end{equation}
	while the second one follows from \Cref{lem:pre_lln}.
	In particular, $\es^\sharp \comp p^{(\N)}$ is copyable and we can apply the positivity axiom to obtain the desired equality.
\end{proof}

	Instead of assuming the existence of empirical sampling morphisms for each object in $\cC$, we could also just consider a fixed object $X$ and an empirical sampling morphism for it.
	In this case, the resulting synthetic Glivenko--Cantelli theorem still holds just as above with the same proof, but of course only for morphisms $f$ and $p$ for this fixed object $X$.

To recover the standard measure-theoretic version of the theorem \cite[Theorem 20.6]{billingsley}, we employ the empirical sampling morphism $\es_\R$ constructed in \Cref{def:es_R}.

\begin{cor}[Standard Glivenko--Cantelli theorem]
	\label{thm:glivenko}	
	Let $(x_i)_{i \in \N}$ be a sequence of real-valued random variables with an IID law $p^\N$.
	Then we have
	\begin{equation}
		\lim_{n \to \infty} \frac{| \{ i \le n \mid x_i \le t \}|}{n} = \P{x_1 \le t}
	\end{equation}
	$p^\N$-almost surely and uniformly in $t \in \R$.
\end{cor}

\begin{proof}
	By \Cref{rem:borelstoch_assumptions}, the category of partial Markov kernels between standard Borel spaces, $\Par{\BorelStoch}$, satisfies \Cref{assumptions}.
	Let us interpret the probability measure $p$ as a morphism $p : I \to \R$ in this category and consider the empirical sampling morphism $\es_\R : P \R \to \R$ from \Cref{def:es_R}.
	With this, $\es_\R^\sharp : \R^\N \to P\R$ is the measurable map that sends a sequence of real numbers to its empirical measure in the sense of \Cref{def:es_R}, whenever it is defined.

	By \Cref{eq:lln_1}, we thus have 
	\begin{equation}\label{eq:glivenko_proof_1}
		\input{glivenko_proof_1.tikz}
	\end{equation}
	where $p^\N$ is the IID measure on $X^\N$.
	As shown in \cite[Example 13.3]{fritz2019synthetic}, this is the string-diagrammatic way to express measure-theoretic almost sure equality.
	In this case \Cref{eq:glivenko_proof_1} says that, $p^\N$-almost surely, the empirical measure $\es^\sharp$ does not depend on the sequence and is simply equal to the measure $p^\sharp$:
	\begin{equation}
		\label{eq:glivenko_proof_2}
		\es^\sharp \bigl( (x_i) \bigr) \;\ase{p^\N}\; p^\sharp.
	\end{equation}
	Using the definition of $\es_\R$ shows that \Cref{eq:glivenko_proof_2} is exactly the desired statement.
\end{proof}

\begin{rem}
	A well-known generalization of the Glivenko--Cantelli theorem is the \emph{Vapnik--Chervo\-nenkis theorem} \cite{vapnik1971uniform}, which replaces the collections of intervals $(-\infty, t]$ in $\R$ by an arbitrary class of measurable sets satisfying a suitable ``finite-dimensionality'' condition and concludes almost sure uniform convergence on this class.
	It seems plausible to us that this could be incorporated into our framework by constructing an empirical sampling morphism for each such class of sets, but we have not yet worked out the details.
\end{rem}

We turn to a few more consequences that can be formulated and proven synthetically.
Comparing the empirical adequacy of empirical sampling with the statement of one version of the synthetic de Finetti theorem \cite[Theorem 2.11]{chen2024aldoushoover} suggests that $\es$ could play the role of the \emph{tail conditional}.
Intuitively, both morphisms can be used to generate the individual outputs $X$ of an exchangeable morphism given the full sequence in $X^\N$.
The following result makes this precise.

\begin{cor}[Tail conditionals]\label{prop:es_tail_cond}
	Given \Cref{assumptions}, let $f : A \to X^{\N}\otimes Y$ be exchangeable in the first factor.
	Then the morphism
	\begin{equation}
		\discard_A \otimes \es \otimes \discard_Y \; : \; A \otimes X^\N \otimes Y \longrightarrow X
	\end{equation}
	is a conditional\footnote{See \cite[Definition 11.5]{fritz2019synthetic} for a definition of conditionals.} of $f$ with respect to $X^{\N\setminus\{1\}} \otimes Y$.
\end{cor}

\begin{proof}
	By the exchangeability of $f$, we can restrict to the first $X$ output without loss of generality. 
	We then get\footnote{We write $X^\N$ instead of $X^{\N\setminus\{1\}}$ in the string diagrams to save space.}
	\begin{equation}\label{eq:es_tail_conditional}
		\input{es_tail_conditional.tikz}
	\end{equation}
	where we use the empirical adequacy of $\es$ in the first and last equality, while the second one is a consequence of \Cref{eq:lln_1} with $p$ given by $\es$.
	The rightmost diagram includes a `splitting' of the double wire, which stands for the isomorphism from \Cref{eq:bifurcation}.
	\Cref{eq:extract_copy} is also used in the last equality in \eqref{eq:es_tail_conditional}.
	We can then calculate
	\begin{equation}
		\input{es_tail_conditional_2.tikz}
	\end{equation}
	where the first equation is by the spreadability lemma~\cite[Lemma 5.1]{fritz2021definetti}, which is a consequence of exchangeability\footnote{Although~\cite[Lemma~5.1]{fritz2021definetti} was proven in the context of Markov categories, the proof is the same also for quasi-Markov categories.} and counitality of copying.
	This equation witnesses that the dashed box morphism $\discard_A \otimes \es \otimes \discard_Y$ is a conditional $f_{|X^{\N\setminus\{1\}} \otimes Y}$.
\end{proof}

\subsection{Synthetic Strong Law of Large Numbers}
\label{sec:synthetic_lln}

The strong law of large numbers for bounded variables is an immediate consequence of the measure-theoretic Glivenko--Cantelli theorem.
In the categorical formulation, our synthetic Glivenko--Cantelli theorem also specializes immediately to a synthetic strong law of large numbers.

\begin{cor}[Synthetic strong LLN]\label{cor:genericLLN}
	Given \Cref{assumptions}, we have\footnote{Recall the symbol $\ase{p^\N}$ refers to equality $p^\N$-almost surely (\Cref{def:ase}).}
	\begin{equation}\label{eq:genericLLN2}
		\input{genericLLN2.tikz}
	\end{equation}
	for arbitrary morphisms ${p : I \to X}$ and ${m : PX \to X}$.
\end{cor}

\begin{proof}
	This follows directly from \Cref{eq:lln_1}, with $A = I$ as in \Cref{eq:glivenko_proof_1}, upon post-composing with $m \otimes \id_{X^\N}$.
\end{proof}

To recover a measure-theoretic law of large numbers, the idea is to take $X = Y = \R$ and $m : P \R \to \R$ to be the expectated value map which assigns to every probability measure its mean,
\begin{equation}
	E : \mu \longmapsto \int_\R x \, \mu(\mathrm{d} x).
\end{equation}
Of course, once again this is only a partial map $P \R \to \R$, since the integral generally does not converge.
For this reason, we define the domain of $E$ to be the set of all probability measures with finite first moment:
\begin{equation}
	\dom{E} \coloneqq \left\{ \mu \in P \R \;\bigg|\; \int_\R |x| \, \mu(\mathrm{d} x) < \infty \right\}.
\end{equation}
The subtle ways in which $E$ can be thought of as a partial $P$-algebra are discussed in \cite{shahmohammed2025partial}.

The following observation is key for interpreting the left-hand side of \Cref{eq:genericLLN2}.

\begin{lem}
	In $\Par{\BorelStoch}$ and with the empirical sampling morphism $\esav_{\R} : \R^\N \to P \R$ from \Cref{def:es_av}, every sequence $(x_i) \in \dom{\esav_\R}$ satisfies:
	\begin{enumerate}
		\item\label{it:empirical_averaging_domain}
			$E \comp \esav_{\R}^\sharp$ is defined on $(x_i)$ if and only if the empirical measure $\esav_{\R}^\sharp((x_i))$ has finite first moment.
		\item\label{it:empirical_averaging_recover}
			In this case, we have
			\begin{equation}
				\label{eq:empirical_averaging_recover}
				\left( E \comp \esav_{\R}^\sharp \right) \bigl( (x_i) \bigr)
				= \lim_{n \to \infty} \frac{1}{n} \sum_{i=1}^n x_i.
			\end{equation} 
	\end{enumerate}
\end{lem}

\begin{proof}
	Property \ref{it:empirical_averaging_domain} holds by the definition of $E$ and Property \ref{it:empirical_averaging_recover} by the definition of $\esav_{\R}$, and in particular by \Cref{eq:es_integral_E}.
\end{proof}

Based on this, we can now instantiate \Cref{cor:genericLLN} to obtain Kolmogorov's strong law of large numbers.

\begin{cor}[Strong law of large numbers]
	\label{cor:strong_lln}
	Let $(x_i)$ be a sequence of real-valued IID random variables with finite first moment.
	Then
	\begin{equation}
		\lim_{n \to \infty} \frac{1}{n} \sum_{i=1}^n x_i = \E{x_1}
	\end{equation}
	almost surely.
\end{cor}

\begin{proof}
	Take $p$ to be the law of the $x_i$ in \Cref{cor:genericLLN}.
	Then the left-hand side of \Cref{eq:genericLLN2} becomes the empirical average
	\begin{equation}
		\lim_{n \to \infty} \frac{1}{n} \sum_{i=1}^n x_i
	\end{equation}
	by \Cref{eq:empirical_averaging_recover}.
	The right-hand side of \Cref{eq:genericLLN2} is $E \comp p^\sharp$, which is the expectation value of the law $p$, also written as $\E{x_1}$.
\end{proof}

	Note that it is important that we use the modified empirical sampling morphism $\esav_{\R}$ instead of the original $\es_{\R}$ from \Cref{def:es_R}.
	Although it is possible to instantiate \Cref{cor:genericLLN} with $\es_{\R}$ and $m = E$, this does not yield the strong law of large numbers in its usual form, since the left-hand side need not be the empirical average, as we saw in \Cref{ex:escaping_sequence}.

\bibliographystyle{alphaurl}
\bibliography{markov}

\clearpage
\appendix

\section{Measure-Theoretic Preliminaries}
\label{sec:measure_theory}

Here, we develop some preliminary results towards proving empirical adequacy for $\es_\R$ and $\esav_{\R}$ introduced in \Cref{def:es_R,def:es_av} respectively.
Specifically, the bounds that we derive in this section are used in \Cref{sec:proofs} to show that the domain of the respective empirical sampling morphism has probability $1$ with respect to any exchangeable measure.
In other words, they ensure that the empirical measure of an exchangeable sequence of random variables exists almost surely.

The ingredients we need here are: Markov's inequality \cite[Section~5.6]{billingsley} the Borel--Cantelli lemma \cite[Section~4.5]{billingsley} and elementary properties of cumulants \cite[Section~9.1]{billingsley}.
It is important to note that our proof of the axioms of empirical sampling morphisms do not use results such as the Glivenko--Cantelli theorem.
Otherwise, we could not claim to prove the standard version thereof (\Cref{thm:glivenko}) from basic principles, since the proof would have been circular.

As is customary in measure-theoretic probability, we denote random variables by uppercase letters in this section. 
This is in contrast to the rest of this document, where uppercase letters generally denote objects in a category.

\begin{lem}
	\label{lem:exchangeable_convergence}
	Let $(Z_i)_{i \in \N}$ be an exchangeable sequence of $[0,1]$-valued random variables and consider
	\begin{equation}
		S_n \coloneqq \sum_{i=1}^n Z_i.
	\end{equation}
	Then:
	\begin{enumerate}
		\item\label{it:exchangeable_concentration}
			For $m, n \in \N$ and $\eps > 0$, we have
			\begin{equation}
				\label{eq:exchangeable_concentration}
				\P*{ \, \vphantom{\rule{0mm}{6mm}} \abs{\frac{S_n}{n} - \frac{S_m}{m}}  > \varepsilon \, } \;\le\; \frac{C}{\varepsilon^6 \min(m, n)^3}
			\end{equation}
			for some universal constant $C$.
		\item\label{it:exchangeable_convergence}
			The sequence $\left(\frac{S_n}{n} \right)_{n \in \N}$ converges almost surely.
	\end{enumerate}
\end{lem}
\noindent %FIXED indent
While our proof bears similarity to a standard proof of the law of large numbers for bounded random variables~\cite[Section~1.6]{billingsley}, our bound is tighter and more general, as it also applies in the exchangeable case. 
Due to this tighter bound, the proof is harder and employs a less straightforward argument based on cumulants.
On the other hand, we make no statement about what the limit of the sequence $\left(\frac{S_n}{n} \right)_{n \in \N}$ is.

\begin{proof}
	For \ref{it:exchangeable_concentration}, we assume $m \ge n$ without loss of generality.
	Then Markov's inequality gives
	\begin{align}
		\P*{ \, \vphantom{\rule{0mm}{6mm}} \abs{\frac{S_n}{n} - \frac{S_m}{m} } > \varepsilon \, } & = \P*{ \left( m S_n - n S_m \right)^6 > \varepsilon^6 m^6 n^6 } \nonumber\\
			& \le \frac{1}{\varepsilon^6 m^6 n^6} \cdot \E*{(m S_n - n S_m)^6} \nonumber\\
			& = \frac{1}{\varepsilon^6 m^6 n^6} \cdot \E*{\left( (m - n) \sum_{i=1}^n Z_i - n \sum_{i=n+1}^m Z_i \right)^6}.
			\label{eq:6th_moment_expression}
	\end{align}
	Thus we need to control the sixth moment of the ``discrepancy''
	\begin{equation}
		\label{def_D}
		D \coloneqq (m - n) \sum_{i=1}^n Z_i - n \sum_{i=n+1}^m Z_i.
	\end{equation}
	In the following, we will finish the proof by showing that this can be bounded by $O(n^3 m^6)$.

	When expanding out the sixth power and taking expectations, each term ends up being an expectation of a product of six of the $Z_i$'s.
	When all six factors are distinct, then we can use exchangeability to rewrite this expectation as $\E{Z_1 \cdots Z_6}$.
	Similarly by exchangeability, all other terms are multiples of one of
	\begin{equation}
		\E*{Z_1^2 Z_2 \cdots Z_6}, \quad \ldots, \quad \E*{Z_1^3 Z_2^2 Z_3}, \quad \ldots, \quad \E*{Z_1^6},
	\end{equation}
	which formally means that the terms are indexed by the partitions of $6$.
	The coefficient of each term is a polynomial in $n$ and $m$.

	Although these coefficient polynomials can be computed explicitly by performing this expansion, actually doing so to the relevant order is cumbersome, and we take a different route through the cumulants.
	Since the polynomials do not depend on the $Z_i$, we can assume without loss of generality that the $Z_i$ are independent.
	In this case, the additivity of cumulants under sums of independent variables together with the fact that the $j$-th cumulant is homogeneous of degree $j$ lets us write the $j$-th cumulant of~\Cref{def_D} as
	\begin{equation}
		\label{eq:cumulants_bound}
		\kappa_j(D) = (m - n)^j n \, \kappa_j(Z_1) + (-n)^j m \, \kappa_j(Z_1) = O(n m^j).
	\end{equation}
	The second step uses $m \ge n$, as well as the assumption that the $Z_i$'s are $[0,1]$-valued so that the cumulants can be bounded by a universal constant.
	Using now the standard conversion formula\footnote{While the general formula is given e.g.~at \cite[Eq.~(6)]{mccullagh1984cumulants}, it is worth noting that for our purposes we only need to know that this formula corresponds to the sum over partitions mentioned above and that $\kappa_1(D) = 0$.} expressing the sixth moment in terms of cumulants and specializing it to the case $\kappa_1(D) = \E{D} = 0$, we get
	\begin{equation}
		\E*{D^6} = \kappa_6(D) + 15 \kappa_4(D) \kappa_2(D) + 10 \kappa_3(D)^2  + 15 \kappa_2(D)^3.
	\end{equation}
	By~\Cref{eq:cumulants_bound}, this can be bounded by
	\begin{equation}
		\E*{D^6} = O(n m^6) + O(n^2 m^6) + O(n^2 m^6) + O(n^3 m^6) = O(n^3 m^6),
	\end{equation}
	as was to be shown.

	Concerning claim~\ref{it:exchangeable_convergence}, the Borel--Cantelli lemma combined with Inequality \eqref{eq:exchangeable_concentration} shows that the sequence $\left(\frac{S_n}{n} \right)_{n \in \N}$ is almost surely Cauchy, and hence converges almost surely.
\end{proof}

Of course, using similar arguments one can establish analogous bounds for any exponent of $\min(m, n)$ in the denominator of~\Cref{eq:exchangeable_concentration}, where larger exponents require more work.
An exponent of $2$ would be sufficient to conclude the strong law of large numbers for bounded random variables.
The reason we use an exponent of $3$ is that this is sufficient to make the following simple proof of a version of the Glivenko--Cantelli theorem work, while an exponent of $2$ would not be.

\begin{prop}
	\label{prop:exchangeable_ecdf}
	Let $(W_i)_{i \in \N}$ be an exchangeable sequence of $\R$-valued random variables.
	Then the empirical cumulative distribution functions $F_n : \R \to [0,1]$ defined by
	\begin{equation}
		F_n(t) \coloneqq \frac{ \abs{ \Set{ i \le n \given W_i \le t } }}{n}
	\end{equation}
	satisfy:	
	\begin{enumerate}
		\item\label{it:exchangeable_ecdf_concentration}
			For $m, n \in \N$ and $\eps > 0$, we have
			\begin{equation}
				\label{eq:exchangeable_ecdf_concentration}
				\P*{ \sup_{t \in \R} \, \abs{F_n(t) - F_m(t)}  > \eps } \le \frac{C}{\eps^6 \min(m, n)^2}
			\end{equation}
			for some universal constant $C$.
		\item\label{it:exchangeable_ecdf_convergence}
			The sequence $(F_n)_{n \in \N}$ converges almost surely uniformly to the cumulative distribution function of a probability measure on $\R$.
	\end{enumerate}
\end{prop}
\noindent %FIXED indent
Our proof is similar in spirit to a standard proof of the Glivenko--Cantelli theorem from the strong law of large numbers~\cite[Theorem~19.1]{vandervaart1998asymptotic}.

\begin{proof}
	For \ref{it:exchangeable_ecdf_concentration}, suppose again $m \ge n$ without loss of generality.
	Then the function $F_n$ jumps up at the points $t = W_1, \ldots, W_n$ and is otherwise constant.
	Since $F_m$ is also monotonically nondecreasing, this implies that the supremum is attained at one of the points $W_i$ (if $F_n$ is above $F_m$ at the supremum) or just before (if $F_n$ is below $F_m$ at the supremum).
	Hence there are $2n$ possibilities for where the supremum can be attained, and we can write the event on the left-hand side of~\Cref{eq:exchangeable_ecdf_concentration} as
	\begin{align*}
		\exists \, j \le n \: : \: \bigg( & \frac{|\Set{ i \le n \given W_i \le W_j }|}{n} - \frac{|\Set{ i \le m \given W_i \le W_j }|}{m} > \varepsilon \\
			& \lor \quad \frac{|\Set{ i \le n \given W_i < W_j }|}{n} - \frac{|\Set{ i \le m \given W_i < W_j }|}{m} < -\varepsilon \bigg)
	\end{align*}
	By \Cref{lem:exchangeable_convergence}, we have the bounds\footnote{Technically the case $i = j$ needs to be considered separately, but this does not affect the conclusion, since a change of $\frac{1}{n}$ in the relative frequency is negligible for $n \gg \eps^{-1}$, and this is the relevant case as we consider asymptotics in $n$ for fixed $\eps$.}
	\begin{align*}
		\P*{ \frac{|\Set{ i \le n \given W_i \le W_j }|}{n} - \frac{|\Set{ i \le m \given W_i \le W_j }|}{m} > \varepsilon } \le \frac{C}{\varepsilon^6 n^3}, \\[1pt]
		\P*{ \frac{|\Set{ i \le n \given W_i < W_j }|}{n} - \frac{|\Set{ i \le m \given W_i < W_j }|}{m} < -\varepsilon } \le \frac{C}{\varepsilon^6 n^3},
	\end{align*}
	for every $j$ and for some universal constant $C$.
	The claimed~\Cref{eq:exchangeable_ecdf_concentration} now follows by a union bound over the $2n$ possibilities.

	Finally, let us prove \ref{it:exchangeable_ecdf_convergence}.
	Since we have a square in the denominator of~\Cref{eq:exchangeable_ecdf_concentration} and $\sum_{n=1}^\infty \frac{1}{n^2}$ converges, we can again apply the Borel--Cantelli lemma together with \ref{it:exchangeable_ecdf_concentration} to conclude that the sequence $(F_n)$ is almost surely uniformly Cauchy. 
	Furthermore, part \ref{it:exchangeable_convergence} of \Cref{lem:exchangeable_convergence} implies that it almost surely converges pointwise.
	Combining these two things with the fact that the uniform limit of a sequence of cumulative distribution functions is a cumulative distribution function, we obtain the desired result.
\end{proof}

To deal with averages of unbounded random variables, we need some control over tails as provided by the following lemma.
Our proof is essentially Garsia's simple proof of the maximal ergodic theorem~\cite{garsia1965ergodic}, although our statement is only the special case for exchangeable variables.

\begin{lem}
	\label{lem:maximal_ergodic}
	Suppose that $(Y_n)_{n \in \N}$ is an exchangeable sequence of nonnegative random variables with finite expectation.
	Then for every $r > 0$, we have
	\begin{equation}\label{eq:maximal_ergodic}
		\P*{ \sup_n \frac{1}{n} \sum_{i=1}^n Y_i > r } \le r^{-1} \E{Y_1}.
	\end{equation}
\end{lem}
\noindent %FIXED indent
Note that the right-hand side is precisely the bound that one would get from Markov's inequality for fixed $n$ (if the supremum was not there).

\begin{proof}
	Consider the sequence of measurable subsets of the underlying probability space $\Omega$ given by
	\begin{equation}
		E_n \coloneqq \Set{ \omega \in \Omega  \given  \left( \max_{k=1,\dots,n}  \frac{1}{k} \sum_{i=1}^k Y_i(\omega) \right) > r }
	\end{equation}
	so that \Cref{eq:maximal_ergodic} can be expressed as $\P{E_\infty} \le r^{-1} \E{Y_1}$, where $E_\infty$ is $\bigcup_n E_n$.
	
	Consider now the shifted version $Z_n \coloneqq Y_n - r$
	of the sequence so that we can write 
	\begin{equation}
		E_n = \left\{ \max_{k=0,\dots,n} \sum_{i=1}^k Z_i > 0 \right\} 
	\end{equation}
	using the more standard implicit notation for random variables.
	We use the convention  $\sum_{i=1}^0 Z_i = 0$, so that we can add the extra case of $k = 0$, which is redundant as $0 > 0$ never holds.
	By
	\begin{equation}
		\label{eq:maximal_ergodic_proof_2}
		\max_{k=0,\dots,n} \sum_{i=1}^k Z_i  \;\le\; \max_{k=0,\dots,n+1} \sum_{i=1}^k Z_i , 
	\end{equation}
	we have $E_n\subseteq E_{n+1}$ and also the pointwise inequality
	\begin{equation}
		\label{eq:maximal_ergodic_proof}
		\max_{k=0,\dots,n} \sum_{i=1}^k Z_i \le 1_{E_n} Z_1 + \max_{k=0,\dots,n} \sum_{i=1}^{k} Z_{i+1}
	\end{equation}
	for each $\omega \in \Omega$.
	Inequality \eqref{eq:maximal_ergodic_proof} can be shown by distinguishing cases:
	\begin{itemize}
		\item If $E_n$ holds, then the maximum on the left is attained for some $k \geq 1$.
			But then the very same terms appear on the right-hand side which can be rewritten as the right-hand side of Inequality \eqref{eq:maximal_ergodic_proof_2} in this case.
		\item If $E_n$ does not hold, then the maximum on the left is attained at $k = 0$, and the inequality simply states its right-hand side is nonnegative, which is trivial.
	\end{itemize}
	If we now take expectations on both sides of Inequality \eqref{eq:maximal_ergodic_proof}, then the two maxima coincide by exchangeability.
	Therefore we obtain
	\begin{equation}
		\E{1_{E_n} Z_1} \ge 0.
	\end{equation}
	Since the sequence of events $E_n$ is increasing in $n$, we obtain the same inequality for $E_\infty$.
	Plugging in the definition of $Z_1$ into $\E{1_{E_\infty} Z_1} \ge 0$
	then gives the second step in
	\begin{equation}
		\E*{Y_1} \ge \E*{1_{E_\infty} Y_1} \ge  \E*{1_{E_\infty}\, r} = r \, \P*{ \sup_n \frac{1}{n} \sum_{i=1}^n Y_i > r },
	\end{equation}
	which amounts to the desired result.
\end{proof}

\section{Omitted Proofs}
\label{sec:proofs}

This section is devoted to proving the theorems in \Cref{sec:es_BorelStoch} plus \Cref{prop:tightened_es} in \Cref{sec:es_integration} based on the results of \Cref{sec:measure_theory}.

\begin{proof}[\textbf{Proof of \Cref{lem:finite_es}}]
	First of all, to show that $\es_F$ is indeed a partial Markov kernel, we need to prove that 
	\begin{enumerate}
		\item \label{it:esF_domain} the domain $\dom{\es_F}$ of $\es_F$ (i.e.\ the set where the limit \eqref{eq:es_finite} exists) is a measurable subset of $F^\N$,
		\item \label{it:esF_measurable} $\es_F(T | \ph) : \dom{\es_F} \to [0,1]$ is a measurable map for each $T \subseteq F$, and that
		\item \label{it:esF_distribution} $\es_F \bigl( \ph | (x_i) \bigr) : 2^F \to [0,1]$ is a probability measure for every $(x_i)$ in the domain.
	\end{enumerate}
	\noindent %FIXED indent
	To prove Property \ref{it:esF_domain}, it suffices to show that the set of sequences for which the limit exists for fixed $T$ is measurable, since $\dom{\es_F}$ is a finite intersection of these. 
	In order to show this, we can replace the existence of the limit by the equivalent condition that the sequence of relative frequencies is Cauchy. 
	Namely, this set of sequences is characterized by requiring that for all $\ell \in \N$, there exists $N \in \N$ such that for all integers $m,n \ge N$, we have
	\begin{equation}
		\abs{ \frac{ \abs{\Set{ i \le n \given x_i \in T} } }{n} - \frac{ \abs{\Set{ i \le m \given x_i \in T} } }{m} } < \ell^{-1}.
	\end{equation}
	This shows that the set of sequences with well-defined empirical measure appears at level $\mathbf{\Pi}^0_4$ of the Borel hierarchy~\cite{kechris}, and is therefore measurable.
	The measurability of $\es_F(T | \ph)$ on this domain, which is Property~\ref{it:esF_measurable}, follows by a similar argument.
	Finally, the additivity of $\es_F( \ph | (x_i))$ follows from the finite additivity of limits, and normalization holds because the relative frequency of $T = F$ is $1$ for every $n$.
	In conclusion, $\es_F$ is a partial Markov kernel $F^\N \to F$, i.e.\ a morphism of $\Par{\BorelStoch}$, as required.
	
	Let us show now that the axioms of \Cref{def:es} are satisfied.
	Permutation invariance holds because both the existence and the values of the limits of relative frequencies are invariant under any finite permutation of the sequence $(x_i)$.
	
	Proving empirical adequacy takes more work.
	We do so by first showing that the kernel
	\begin{equation}
		\es_F^{(\N)} : F^\N \to F^\N
	\end{equation}
	can be computed in terms of a limit of averaging over the finite permutation groups $S_n$.
	Namely, we claim that for any subsets $T_1, \ldots, T_m \subseteq F$ and any $m \in \N$, we have
	\begin{equation}
		\label{eq:es_random_permutation}
		\es_F^{(\N)} \bigl( T_1 \times \cdots \times T_m \times F \times \dots \,|\, (x_n) \bigr) = \lim_{n \to \infty} \frac{1}{n!} \sum_{\sigma \in \perms{n}} T_1(x_{\sigma(1)}) \, \cdots \, T_m(x_{\sigma(m)}),
	\end{equation}
	where $T_i(\ph)$ denotes the indicator function of $T_i$, and we leave it understood that the limit is over $n \ge m$ only.
	This formula suggests that we may think of $\es_F^{(\N)}$ as applying a uniformly random permutation to the given sequence $(x_i)$, whenever this makes sense (i.e.\ when the limits exist).
	As part of the claim that \Cref{eq:es_random_permutation} coincides with \eqref{eq:es_finite}, 
	we claim that the existence of the relevant limits holds on precisely the same set of sequences $(x_i)$.

	To prove \Cref{eq:es_random_permutation}, suppose first that $(x_i) \in \dom{\es_F}$.
	We then need to show that the limits on the right-hand side exist for all $m$ and that the equation holds.
	By definition of $\es_F^{(\N)}$, we can compute
	\begin{align*} %FIXED Reduced Spacing for Margins
			\es_F^{(\N)} \bigl(T_1 \times \cdots \times T_m \times F \times \dots \,|\, (x_j) \bigr) & = \es_F \bigl(T_1 | (x_j) \bigr) \, \cdots \, \es_F \bigl(T_m | (x_j) \bigr) \\[3pt]
			&= \lim_{n \to \infty} \frac{\abs{\Set{ i_1 \le n \!\given\! x_{i_1} \in T_1 } } }{n} \, \cdots  \frac{\abs{\Set{ i_m \le n \!\given\! x_{i_m} \in T_m } } }{n} \\
			&= \lim_{n \to \infty} \frac{1}{n^m} \sum_{i_1, \dots, i_m = 1}^n T_1(x_{i_1}) \,\cdots\, T_m(x_{i_m}).
	\end{align*}
	Thinking of $k \mapsto i_k$ as a map $\{1, \ldots, m\} \to \{1, \ldots, n\}$ shows that this amounts to averaging over all such maps, of which there are $n^m$ many.
	Since the fraction of injections among these maps goes to $1$ as $n \to \infty$, we obtain the same limit by averaging over all injections.
	But every such injection can be extended to a permutation of $\{1,\ldots,n\}$, there are $(n - m)!$ many ways to do so independently of what the injection is, and all of these are distinct as the map varies.
	Therefore, we can instead average over permutations, and this is exactly the claimed \Cref{eq:es_random_permutation}.
	
	Conversely, suppose that $(x_j)$ is such that the limits in \Cref{eq:es_random_permutation} exist for all $m$ and all $T_1, \ldots, T_m$.
	Then this holds in particular for $m = 1$, in which case the right-hand side coincides with that of \Cref{eq:es_finite} and the claim follows.
	In conclusion, $\es_F^{(\N)}$ and the formula from \Cref{eq:es_random_permutation} have the same domain and, on this domain, they coincide.
	
	With this formula for $\es_F^{(\N)}$ at hand, we can prove empirical adequacy for $\es_F$.
	To this end, consider an arbitrary standard Borel space $Y$ as in \Cref{eq:es_invariance}.
	Then for a Markov kernel $f : A \to F^{\N} \otimes Y$ that is exchangeable in its first factor, let us first show that
	\begin{equation}
		\label{eq:exchangeable_es_dom}
		f \bigl(\dom{\es_F} \times Y \, | \, a \bigr) = 1 \qquad \forall a \in A.
	\end{equation}
	This amounts to showing that if $(x_n)_{n \in \N}$ and $y$ are random variables taking values in the respective spaces, then the limit
	\begin{equation}
		\lim_{n \to \infty} \frac{| \{ i \le n \mid x_i \in T \}|}{n} = \lim_{n \to \infty} \frac{1}{n} \sum_{i=1}^n T(x_i)
	\end{equation}
	exists almost surely with respect to the joint distribution $f( \ph | a)$ of $F^\N \otimes Y$.
	This follows from \Cref{lem:exchangeable_convergence}.
	
	We can now calculate that, for arbitrary subsets $T_1, \ldots, T_m \subseteq F$ and any measurable $U \subseteq Y$,
	\begin{equation}
		\label{es_proof_finite}
		\begin{split}
			\left( \bigl( \es_F^{(\N)} \otimes \id_Y \bigr)  \comp f \right)& \bigl( T_1 \times \cdots \times T_m \times F \times \dots \times U \,|\, a \bigr)  \\
			& = \int_{x \in \dom{\es_F}} \es_F^{(\N)} \bigl( T_1 \times \cdots \times T_m \times F \times \dots | x \bigr) \, f(\mathrm{d} x \times U \,|\, a) \\
			& = \int_{x \in \dom{\es_F}} \lim_{n \to \infty} \frac{1}{n!} \sum_{\sigma \in \perms{n}} T_1(x_{\sigma(1)}) \cdots T_m(x_{\sigma(m)}) \, f(\mathrm{d} x \times U \,|\, a) \\
			& = \lim_{n \to \infty} \frac{1}{n!} \sum_{\sigma \in \perms{n}} \int_{x \in \dom{\es_F}} T_1(x_{\sigma(1)}) \cdots T_m(x_{\sigma(m)}) \, f(\mathrm{d} x \times U \,|\, a) \\[1pt]
			& = f \bigl( T_1 \times \cdots \times T_m \times F \times \dots \times U \,|\, a \bigr), 
		\end{split}
	\end{equation}
	where the first equation is plugging in definitions, the second one uses \Cref{eq:es_random_permutation}, 
	the commutation of the limit and the integral in the third is an instance of the dominated convergence theorem (the measure is finite and the integrand is bounded by $1$), and the last step holds by the exchangeability of $f$ and by \Cref{eq:exchangeable_es_dom}.
	This concludes the proof of empirical adequacy for $\es_F$.
\end{proof}

\begin{rem}
	Writing $\es_F$ as in averaging over permutations as in \Cref{eq:es_random_permutation} is strongly reminiscent of group averaging in ergodic theory~\cite{lindenstrauss2001pointwise}.
	We expect this to be connected to the fact that the permutation group $S_\infty = \bigcup_n S_n$ is amenable, and that similar constructions can be made for other amenable groups.
\end{rem}

\begin{proof}[\textbf{Proof of \Cref{lem:countable_es}}]
This works similarly as the construction from the finite case with some additional subtleties.
Recall from \Cref{eq:es_uniform} that we require the limits of~\Cref{eq:es_countable} to exist \emph{uniformly} for the sets $T = \{1, \ldots, t\}$, but otherwise the proof of measurability of $\dom{\es_\N}$ and of $\es_\N$ itself is the same as in the finite case.
The fact that $\es_\N(\ph | (x_i))$ is a probability measure for every $(x_i) \in \dom{\es_\N}$ holds because of \Cref{eq:es_additive}.

The rest of the proof is the same as in the finite case, but with all limits being suitably uniform; in particular, for fixed $m$ the limit in~\Cref{eq:es_random_permutation} must exist uniformly for sets of the form $T_i = \{1, \ldots, t_i\}$.
The analogue of \Cref{eq:exchangeable_es_dom} holds because the relevant limits are all uniform in $S$, as is immediate from the bound \Cref{eq:exchangeable_concentration} in \Cref{lem:exchangeable_convergence}.
\end{proof}

\begin{proof}[\textbf{Proof of \Cref{lem:real_es}}]
	Let us first show that \Cref{def:es_R} indeed gives a unique Markov kernel as claimed.
	To this end, let $F : \Q \to \R$ be the function given by 
	\begin{equation}
		F(t) \coloneq \lim_{n \to \infty} \frac{ \abs{\{ i \le n \mid x_i \le t\}}}{n}
	\end{equation}
	for any sequence in the domain of $\es_\R$.
	Since the function
	\begin{equation}
		t \mapsto \frac{|\{ i \le n \mid x_i \le t\}|}{n}
	\end{equation}
	is right continuous for every fixed $n$, and right continuous functions are stable under uniform limits, it follows that $F$ is right continuous as well. 
	The limits in~\Cref{limit_cdf} likewise hold thanks to the uniform convergence, while monotonicity is already a consequence of the pointwise convergence. Therefore $F$ is a CDF. 
	The probability measure that corresponds to it satisfies \Cref{eq:es_intervals} by finite additivity together with the commutation of limits and finite sums.
	Therefore, a measure $\es_\R(\ph|(x_i))$ with the properties specified in \Cref{def:es_R} indeed exists.

	To see that $\es_\R$ is a partial Markov kernel, we still need to verify that $\dom{\es_\R}$ is measurable and that on this domain, $\es_\R(T | \ph)$ is measurable for every measurable $T \subseteq \R$.
	The former follows by analogous arguments as before, as $\dom{\es_\R}$ consists of all sequences $(x_n)$ such that
	\begin{equation}
		\begin{split}
			&\forall \ell \in \mathbb{N}, \; \exists N \in \mathbb{N} \; : \\
			&\quad \quad \forall m, n \geq N, \; \forall t \in \mathbb{Q} \; : \\
			&\quad \quad \quad \quad \left| \frac{|\{ i \leq n \mid x_i \leq t\}|}{n} 
			- \frac{|\{i \leq m \mid x_i \leq t\}|}{m} \right| < \ell^{-1}.
		\end{split}
	\end{equation}
	which shows that $\dom{\es_\R}$ is a $\mathbf{\Pi}^0_4$ set.
	The measurability of $\es_\R(T | \ph)$ also holds by the same token for the sets of the form $T = (-\infty, t]$ with $t \in \Q$.
	Since these generate the Borel $\sigma$-algebra on $\R$, the general case follows by the $\pi$-$\lambda$ theorem.
	
	The remainder of the proof is analogous to the countable case, with all measurable sets being of the form $(-\infty, t)$ and all limits existing uniformly.
	The corresponding version of \Cref{eq:exchangeable_es_dom} holds thanks to \Cref{prop:exchangeable_ecdf}.
\end{proof}

Before proving \Cref{prop:tightened_es} in full, we first treat the restricted case of $\R_+$. 
To this end, we consider the positive and negative parts of a real number $x \in \R$,
\begin{equation}
	x^+ \coloneqq \max\{x, 0\} \qquad \text{and} \qquad x^- \coloneqq \max\{-x, 0\},
\end{equation}
so that we have $x = x^+ - x^-$ and $\abs{x} = x^+ + x^-$.
We also write 
\begin{align*}
	\pi : \R &\longrightarrow \R_+ \\
	x &\longmapsto x^+
\end{align*}
for the positive part function and $\iota : \R_+ \hookrightarrow \R$ for the inclusion map.
Then
\[
	\es_{\R_+} \coloneqq \pi \comp \es_\R \comp \iota^\N
\]
is an empirical sampling morphism for $\R_+$ by \Cref{lem:transfer_es}.

\begin{lem}
	\label{lem:relu_lax_naturality}
	The empirical sampling morphisms $\es_\R$ and $\es_{\R_+}$ satisfy lax naturality with respect to $\pi$, i.e.
	\begin{equation}
		\label{eq:relu_lax_naturality}
		\es_{\R_+} \comp \pi^\N \domext \pi \comp \es_{\R}.
	\end{equation}
\end{lem}

\begin{proof}
	First let us show that $\dom{\pi \comp \es_{\R}}$ is a subset of $\dom{\es_{\R_+} \comp \pi^\N}$.
	To this end, consider an arbitrary sequence $(x_i) \in \dom{\es_{\R}}$.
	The finite empirical CDFs of the sequence $(x_i^+)$ coincide with those of the original sequence $(x_i)$ on the positive axis (including $0$) while vanishing on the negative axis.
	Therefore the uniform convergence is preserved.
	
	It is also clear that the limiting CDF is the one of the pushforward measure, which likewise arises by truncation of the limiting CDF of the original sequence.
	Therefore, $\es_{\R_+} \comp \pi^\N$ coincides with $\pi \comp \es_{\R}$ when restricted to the domain of the latter.
	This completes the proof of the lemma.
\end{proof}

Similar to \Cref{prop:tightened_es}, we now consider a ``tightened'' version of $\es_{\R_+}$ that we denote by $\esav_{\R_+}$.
This is defined like $\es_{\R_+}$, but with the additional restriction that a sequence $(x_i)$ is in the domain of $\esav_{\R_+}$ if and only if it is in the domain of $\es_{\R_+}$ and additionally
\begin{equation}
	\label{eq:es_integral_Eabs_2}
	\int_{\R_+} y \, \es_\R(\mathrm{d}y | (x_i)) = \lim_{n \to \infty} \frac{1}{n} \sum_{i=1}^n x_i
\end{equation}
holds, where both sides may be infinite.
This defines a partial Markov kernel $\esav_{\R_+} : \R_+^\N \to \R_+$.

\begin{lem}
	\label{lem:es_nonnegative}
	The partial Markov kernel $\esav_{\R_+}$, defined as above, is an empirical sampling morphism for $\R_+$ in $\Par{\BorelStoch}$,
\end{lem}

\begin{proof}
	The proof that $\esav_{\R_+}$ is a partial Markov kernel is analogous to the previous cases.
	While permutation invariance is clear, empirical adequacy is where the main work lies.

	We are essentially using the same empirical sampling as in \Cref{lem:real_es}, the only change being a more restricted domain.
	Therefore, the only additional step to show is that the left-hand side of \Cref{eq:es_invariance} is a total morphism for every total $f$.
	This is what we focus on in the remainder of the proof.	
	To simplify notation, we assume without loss of generality that $f$ has no input, which makes it into a probability measure $p : I \to \R_+^\N \otimes Y$.
	
	Assuming that such $p$ is exchangeable in the first factor, we need to show 
	\begin{equation}
		p \bigl( \dom{\esav_{\R_+}} \times Y \bigr) = 1. 
	\end{equation}
	By what we have already shown in the previous proofs of empirical adequacy,
	this is equivalent to saying that \Cref{eq:es_integral_Eabs_2} holds $p$-almost surely whenever the $(x_i)$ is a random sequence with exchangeable distribution.
	
	For any given $r \ge 0$, we decompose each $x_i$ as $x_i = x_i^{\le r} + x_i^{> r}$, where
	\begin{equation}
	x_i^{\le r} \coloneqq x_i\, 1_{[0,r]}(x_i),
	\qquad x_i^{> r} \coloneqq x_i\,1_{(r,\infty)}(x_i),
	\end{equation}
	so that $x_i$ is equal to exactly one of them.
	We conduct the proof by showing that the second of the two inequalities in
	\begin{equation}
		\label{eq:supinf_ineqs}
		\limsup_{n \to \infty} \frac{1}{n} \sum_{i=1}^n x_i \le \int_{\R_+} y \, \es_\R(\mathrm{d}y \mid (x_i)) \le \liminf_{n \to \infty} \frac{1}{n} \sum_{i=1}^n x_i
	\end{equation}
	holds for every sequence in $\dom{\es_\R}$ and subsequently showing that the first one holds $p$-almost surely.
	These two facts imply that \Cref{eq:es_integral_Eabs_2} holds $p$-almost surely.
	
	The second inequality in \Cref{eq:supinf_ineqs} can be shown by taking the $r \to \infty$ limit of
	\begin{align}
		\label{eq:tightened_es_proof_1} \int_{[0,r]} y \, \es_\R(\mathrm{d}y | (x_i)) &\leq \int_{[0,r]} y \, \es_\R(\mathrm{d}y | (x_i)) + \liminf_{n \to \infty} \frac{1}{n} \sum_{i=1}^n x_i^{> r} \\
		&= \lim_{n \to \infty} \frac{1}{n} \sum_{i=1}^n x_i^{\le r} + \liminf_{n \to \infty} \frac{1}{n} \sum_{i=1}^n x_i^{> r} \\
		&\leq \liminf_{n \to \infty} \frac{1}{n} \sum_{i=1}^n x_i,
	\end{align}
	where the first inequality is by $r \ge 0$, the equality is an instance \Cref{prop:es_integral}, and the last inequality is the superadditivity of the limit inferior.
	
	In order to prove the first inequality in \Cref{eq:supinf_ineqs}, we can assume without loss of generality that the integral on its right-hand side is finite.
	Then, we have
	\begin{align}
		\limsup_{n \to \infty} \frac{1}{n} \sum_{i=1}^n x_i &\leq \lim_{n \to \infty} \frac{1}{n} \sum_{i=1}^n x_i^{\le r} + \limsup_{n \to \infty} \frac{1}{n} \sum_{i=1}^n x_i^{> r} \\
		&= \int_{[0,r]} y \, \es_\R(\mathrm{d}y | (x_i)) + \limsup_{n \to \infty} \frac{1}{n} \sum_{i=1}^n x_i^{> r} \\
		&\leq \int_{[0,r]} y \, \es_\R(\mathrm{d}y | (x_i)) + \sup_{n} \frac{1}{n} \sum_{i=1}^n x_i^{> r},
	\end{align}
	where the first inequality is by the subadditivity of limit superior and the equality is again by \Cref{prop:es_integral}.
	Taking the $r \to \infty$ limit gives
	\begin{equation}
		\limsup_{n \to \infty} \frac{1}{n} \sum_{i=1}^n x_i \leq  \int_{\R_+} y \, \es_\R(\mathrm{d}y | (x_i)) + \lim_{r \to \infty} \sup_{n} \frac{1}{n} \sum_{i=1}^n x_i^{> r}.
	\end{equation}
	It now suffices to argue that the second term on the right-hand side vanishes $p$-almost surely.
	To this end, for a fixed $\eps > 0$, an application of \Cref{lem:maximal_ergodic} gives
	\begin{equation}
		\P*{\sup_n \frac{1}{n} \sum_{i=1}^n x_i^{>r} > \eps} \le \eps^{-1} \E*{x_1^{>r}}.
	\end{equation}
	Since the functions $x \mapsto x^{>r}$ converge to zero pointwise as $r \to \infty$, the monotone convergence theorem implies that $\E{x^{>r}}$ converges to $0$.
	We thus conclude that 
	\begin{equation}
		\lim_{r \to \infty} \sup_n \frac{1}{n} \sum_{i=1}^n x_i^{>r} \le \eps
	\end{equation}
	holds $p$-almost surely for any $\eps > 0$, as was to be shown.
\end{proof}

\begin{proof}[\textbf{Proof of \Cref{prop:tightened_es}}]
	We can generally proceed by the same arguments as in the proof of \Cref{lem:es_nonnegative}.
	However, we still need to show
	\begin{equation}
		\label{eq:dom_esav}
		p \bigl( \dom{\esav_{\R}} \times Y \bigr) = 1
	\end{equation}
	for any exchangeable $p : I \to \R^\N \otimes Y$, i.e.\ that the conditions in \Cref{def:es_av} that are additionally imposed, on top of those in \Cref{def:es_R}, hold almost surely for a random sequence $(x_i)$ with an exchangeable distribution.
	
	To this end, let $L_+ \subseteq \R^\N$ be the set of those sequences, whose positive part is an element of the domain of $\esav_{\R_+}$ and similarly for $L^-$:
	\begin{align}
		L^+ &\coloneq \Set{ (x_i) \in \dom{\es_{\R}}  \given  \int_\R y^+ \, \es_\R\bigl(\mathrm{d}y \,|\, (x_i) \bigr) = \lim_{n \to \infty} \frac{1}{n} \sum_{i=1}^n x_i^+ } \\
		L^- &\coloneq \Set{ (x_i) \in \dom{\es_{\R}}  \given  \int_\R y^- \, \es_\R \bigl(\mathrm{d}y \,|\, (x_i) \bigr) = \lim_{n \to \infty} \frac{1}{n} \sum_{i=1}^n x_i^- } .
	\end{align}
	Note that, for any $(x_i) \in \dom{\es_{\R}}$, we also have
	\begin{equation}
		\label{eq:int_positive_part}
		\int_{\R} y^+ \, \es_\R \bigl( \mathrm{d}y \,|\, (x_i) \bigr) = \int_{\R_+} y \,\, (\pi \comp \es_{\R}) \bigl( \mathrm{d}y \,|\, (x_i) \bigr) = \int_{\R_+} y \, \es_{\R_+} \bigl( \mathrm{d}y \,|\, (x_i^+) \bigr),
	\end{equation}
	where the first equation is by change of variables and the second by the lax naturality from \Cref{lem:relu_lax_naturality}.
	The other integral in \Cref{eq:dom_esav} can be treated analogously.
	
	Since both $(x_i^+)$ and $(x_i^-)$ have an exchangeable distribution, we thus obtain
	\begin{equation}
		p(L^+ \otimes Y) = 1 = p(L^- \otimes Y)
	\end{equation}
	by \Cref{lem:es_nonnegative}, so that $(L^+ \cap L^-) \otimes Y$ also has full measure.
	
	We complete the proof by showing $L^+ \cap L^- \subseteq \dom{\esav_{\R}}$, from which \Cref{eq:dom_esav} then follows.
	We distinguish several distinct cases for the values of 
	\begin{equation}
		\lim_{n \to \infty} \frac{1}{n} \sum_{i=1}^n x_i^+ \qquad \text{and} \qquad \lim_{n \to \infty} \frac{1}{n} \sum_{i=1}^n x_i^- .
	\end{equation}
	\begin{itemize}
		\item If both of the limits are finite, we obtain \Cref{eq:es_integral_Eabs} by adding the two defining equations of $L^+$ and $L^-$, respectively, and we obtain \Cref{eq:es_integral_Eabs_2} by subtracting them.
		
		\item If at least one is infinite, then this limit is necessarily $+\infty$ since both sequences are non-negative, and we get \Cref{eq:es_integral_Eabs} by adding the two equations and getting $+\infty$ too.
	\end{itemize}
	In both cases, $(x_i) \in L^+ \cap L^-$ implies $(x_i) \in \dom{\esav_{\R}}$, and thus the proof is complete.
\end{proof}
\end{document}